\pdfoutput=1
\documentclass[11pt]{article}
\usepackage[a4paper,margin=24mm]{geometry}
\usepackage{amsmath,amssymb,amsthm,mathtools}
\usepackage{booktabs,longtable,array}
\usepackage{enumitem}
\usepackage{xcolor}
\usepackage{placeins}
\usepackage[colorlinks=true,hypertexnames=false,linkcolor=blue!45!black,citecolor=blue!45!black,urlcolor=blue!45!black]{hyperref}

\newcommand{\Z}{\mathbb{Z}}
\newcommand{\Q}{\mathbb{Q}}

\newcommand{\C}{\mathbb{C}}
\newcommand{\F}{\mathbb{F}}
\newcommand{\e}{\mathbf{e}}
\newcommand{\fd}{\mathfrak{d}}
\newcommand{\Mp}{\mathrm{Mp}_2(\Z)}
\newcommand{\Dp}{\mathcal{D}}
\newcommand{\Hh}{\mathcal{H}}
\DeclareMathOperator{\sig}{sig}
\DeclareMathOperator{\tr}{tr}
\DeclareMathOperator{\Irr}{Irr}
\DeclareMathOperator{\SL}{SL}
\DeclareMathOperator{\Hom}{Hom}

\newcommand{\gC}{\textup{[C]}}
\newcommand{\Gs}{\mathcal{G}}
\newcommand{\Bb}{\overline{B}}
\DeclareMathOperator{\Inv}{Inv}
\DeclareMathOperator{\Saw}{Saw}
\DeclareMathOperator{\GL}{GL}
\DeclareMathOperator{\PGL}{PGL}
\DeclareMathOperator{\SO}{SO}
\DeclareMathOperator{\rr}{rr}
\DeclareMathOperator{\RR}{RR}
\newcommand{\LS}{\mathrm{LS}}
\newcommand{\LP}{\mathrm{LP}}
\newcommand{\LO}{\mathrm{LO}}
\newcommand{\cp}{\mathrm{cap}}

\newtheorem{theorem}{Theorem}[section]
\newtheorem{lemma}[theorem]{Lemma}
\newtheorem{proposition}[theorem]{Proposition}
\newtheorem{corollary}[theorem]{Corollary}
\theoremstyle{definition}
\newtheorem{definition}[theorem]{Definition}
\newtheorem{remark}[theorem]{Remark}
\theoremstyle{plain}
\newtheorem*{thmA}{Theorem A}
\newtheorem*{thmOne}{Theorem 1}
\newtheorem*{thmTwo}{Theorem 2}

\newenvironment{cited}[1]{\par\smallskip\begingroup\small\leftskip=7mm\rightskip=5mm\noindent\textbf{#1.}\ \itshape\ignorespaces}%
  {\par\endgroup\smallskip}
\newcommand{\why}[1]{{\small\noindent\emph{Reason for quoting.} #1\par}}

\begin{document}

\begin{center}
{\Large\bfseries Petersson-rigid lattices of signature $(2,1)$\\[1mm] and their harmonic weak Maass forms}\\[5mm]
{Eungang Cho}\\[1mm]
{\small Department of Mathematics, Chung-Ang University}\\[3mm]
{\footnotesize 2026-10-30 \quad (paper III of a series: [I] arXiv:2608.19706, [II] arXiv:2609.02843)}\\
{\footnotesize Reproduction archive: Zenodo, \href{https://doi.org/10.5281/zenodo.23770095}{doi:10.5281/zenodo.23770095}.}
\end{center}

\vspace{2mm}
\begin{center}{\bfseries Abstract}\end{center}
{\small
An even lattice $L$ of signature $(2,1)$ is \emph{Petersson-rigid} if the orthogonal group of its discriminant form acts absolutely irreducibly on the
space of cusp forms of weight $\frac32$ for its Weil representation, a space of dimension at least two. We prove, with computer assistance, that
there are exactly twelve such lattices up to isometry, and we determine their root bases. Four of them are theta lattices, whose cusp forms are lifts
of unary theta series. On each, the space is $\C F_0\otimes\sigma_0$ for one form $F_0$ and
an irreducible $\sigma_0$, so the Petersson norm is $\|F\|^2=t_\mu\|c_F(\mu)\|^2$ for all cusp forms $F$, with one $t_\mu>0$, $\mu$ the least exponent of $F_0$
and $c_F(\mu)$ the coefficient of $q^\mu$ in $F$. We prove that the harmonic weak Maass forms of weight $\frac12$ in the $\sigma_0$-isotypic part with a pole of least possible order are
exactly the forms $\Hh_F$ with principal part $q^{-\mu}\overline{c_F(\mu)}$, and $\xi_{1/2}\Hh_F=F/t_\mu$.

For the four theta lattices it is a theorem that $t_\mu=\frac{s-1}{12}\sqrt s$, with $s$ the square-free part of $|\det L|/2$, and the $\Hh_F$ come from
umbral moonshine. For the eight others we compute $F_0$, $t_\mu$ and $\Hh_{F_0}$: the Hecke eigenvalues are integers and, at the primes computed, the
$a_p$ of an elliptic curve, and for three of the eight $F_0$ is explicit as an eta--theta product.

We prove, for every prime $p>3$, a closed formula in $p$ and one class number for the difference $g(X_0(4p)/w_p)-g(X_0(2p)/w_p)$ of the
genera of two Atkin--Lehner quotients of modular curves. In a family $\Dp_p$ of discriminant forms, $p$ an odd prime,
containing three of the eight, comparison with that classification shows that $\Dp_p$ is rigid exactly when $X_0(4p)/w_p$ has genus $1$ and
$X_0(2p)/w_p$ has genus $0$.
\par}

\vspace{3mm}
{\footnotesize\noindent\textbf{Conventions.} \gC\ marks a statement that rests on a computation; every such computation is reproduced by a script
and a log (Appendix~\ref{app:comp}), and where floating point enters, the margin is stated. A theorem quoted from the literature and used as a step of a
proof is stated in place with the reason for quoting it, and classical facts (Schur's lemma, Deuring's criterion, the Wiener--Ikehara theorem, the genus
formula for $\Gamma_0(N)$, \dots) are named in place; Appendix~\ref{app:cite} lists all of them. Numerical values are labelled as such, and none of them
is used as a step of a proof.\par}

\vspace{2mm}
\tableofcontents
\newpage
\section{Introduction}\label{sec:intro}

For an untwisted affine Kac--Moody algebra the automorphy of the denominator is a theorem. By the Weyl--Kac denominator identity, which specializes to
the Macdonald identities, the denominator is a Weyl-alternating sum of theta functions of level $h^\vee$, and these transform under $\Mp$ by the Weil
representation of the discriminant form of the coroot lattice scaled by $h^\vee$ [Mac72], [Ka90, Chs.~12--13], [KP84, \S\S3.1, 4.1--4.2]. In the
hyperbolic case this can fail. The denominator of the Feingold--Frenkel algebra [FF83] is not an automorphic form [KL13]; Gritsenko and Nikulin found
an automorphic correction, a generalized Kac--Moody algebra with the same real simple roots and further imaginary ones, whose denominator is the Igusa
cusp form $\Delta_{35}$ of weight $35$, with an infinite product expansion [GN96], [KL13, Thm.~4.2].

The root lattice of a symmetrizable hyperbolic Cartan matrix of rank three, in its minimal symmetrization, is an even lattice of signature $(2,1)$.
For an even lattice $L$ of signature $(2,1)$ the finite quadratic module $D(L)=L'/L$ carries the Weil representation $\rho_L$ of the metaplectic group
$\Mp$, and its isometry group $O(D(L))$ acts on every space of vector-valued modular forms for $\rho_L$, commuting with $\Mp$. Papers [I] and [II]
studied two such spaces. At weight $\frac52$, $S_{5/2}(\rho_L)$ is the obstruction space of Borcherds products on $L(-1)\oplus U$ [Bo99, Thm.~3.1], and
the shadow of an $O(D(L))$-invariant prescription lies in its invariant part [I, Thm.~7.6]. At weight $\frac32$, [II] called $L$
\emph{Petersson-rigid} if $O(D(L))$ acts absolutely irreducibly on $S_{3/2}(\rho_L)$ and $\dim S_{3/2}(\rho_L)\ge2$. We determine these lattices.

\begin{thmA}[{proved in \S\ref{sec:proofA}}]
There are exactly twelve Petersson-rigid lattices up to isometry, of determinants $-16,-40,-54,-64,-72,-88,-96,-120,-184,-250,-324,-648$; each of
their genera has class number one. For $L_{16},L_{54},L_{64},L_{250}$ the space $S_{3/2}(\rho_L)$ is spanned by lifts of unary theta series; for the
other eight it contains none.
\end{thmA}

Theorem~A says which lattices are rigid; it does not say what rigidity fixes. This paper is about the eight lattices of Theorem~A without
theta series,
\[
B_{40},\ L_{72},\ B_{88},\ L_{96},\ B_{120},\ B_{184},\ B_{324},\ B_{648},
\]
and asks which numbers and which harmonic weak Maass form rigidity determines for them. The main results are two theorems. Theorem~1 (\S\ref{subsec:nochoice}, proved in \S\ref{sec:canon}): once the pole is required to have the
least possible order, rigidity leaves no normalization to choose; it fixes one number $t_\mu$ and attaches to each $F\in S_{3/2}(\rho_L)$ one form
$\Hh_F$, with shadow $F/t_\mu$. Theorem~2 (\S\ref{subsec:results},
proved in \S\ref{sec:count}): a closed formula, for every prime $p>3$, for the difference of the genera of two Atkin--Lehner quotients of modular
curves, whose values single out the rigid members of one family of discriminant forms by comparison with the list of Theorem~A. Between them we compute the forms (\S\ref{sec:eight})
and their Hecke eigenvalues (\S\ref{sec:arith}).

\begin{table}[h]
\centering\footnotesize\setlength{\tabcolsep}{2.6pt}
\caption{The twelve Petersson-rigid lattices as lattices (Theorem~A, Proposition~\ref{prop:roots}): a model, the discriminant form $D(L)$ (notation of
\S\ref{subsec:fqm}), the quaternion discriminant $\fd$ of $L\otimes\Q$ (\S\ref{subsec:lat}; its Shimura curve is compact exactly when $\fd>1$), whether $L$ is
maximal among the even lattices in $L\otimes\Q$ (\S\ref{subsec:lat}), the hyperbolic root basis, written $M(k)$ with $M$ a row of the census of [II]
(C compact, N non-compact, --- none), the span $\Theta(L)$ of the lifts of unary theta series in $S(D)=S_{3/2}(\rho_L)$, $s=\dim S(D)$ and $|O(D)|$. The index in
the name is $|\det L|$, and every genus has class number one. Above the rule the four theta lattices, below it the eight lattices of this paper.
\gC}
\label{tab:twelve}
\begin{tabular}{lllcclccc}
\toprule
$L$ & model & $D(L)$ & $\fd$ & max. & root basis & $\Theta(L)$ & $s$ & $|O(D)|$\\
\midrule
$L_{16}$ & $U(2)\oplus\langle4\rangle$ & $U(2)\oplus\langle 1/8\rangle$ & $1$ & no & $FF(2)$, N & $S(D)$ & $2$ & $12$\\
$L_{54}$ & $U(3)\oplus\langle6\rangle$ & $\langle 3/4\rangle\oplus 3^{(1,1,1)}[1,1,2]$ & $1$ & no & $FF(3)$, N & $S(D)$ & $3$ & $48$\\
$L_{64}$ & $M_8(2)$ & $V(2)\oplus\langle 5/32\rangle$ & $1$ & no & $M_8(2)$, N & $S(D)$ & $2$ & $48$\\
$L_{250}$ & $U(5)\oplus\langle10\rangle$ & $\langle 1/4\rangle\oplus 5^{(1,1,1)}[1,1,1]$ & $1$ & no & $FF(5)$, N & $S(D)$ & $5$ & $240$\\
\midrule
$B_{40}$ & $\langle2\rangle\oplus\bigl(\begin{smallmatrix}4&2\\2&-4\end{smallmatrix}\bigr)$ & $\langle 1/4\rangle\oplus V(2)\oplus 5^{(1)}[2]$ & $10$ & yes & --- (rim row) & $0$ & $2$ & $12$\\
$L_{72}$ & $A_2(2)\oplus\langle-6\rangle$ & $\langle 1/4\rangle\oplus V(2)\oplus 3^{(1,1)}[1,1]$ & $6$ & yes & $L_{72}(1)$, C & $0$ & $4$ & $48$\\
$B_{88}$ & $\langle2\rangle\oplus\bigl(\begin{smallmatrix}-4&2\\2&10\end{smallmatrix}\bigr)$ & $\langle 3/4\rangle\oplus V(2)\oplus 11^{(1)}[2]$ & $22$ & yes & --- (rim row) & $0$ & $2$ & $12$\\
$L_{96}$ & $A_2(2)\oplus\langle-8\rangle$ & $V(2)\oplus\langle 7/16\rangle\oplus 3^{(1)}[2]$ & $6$ & no & $M_{12}(2)$, C & $0$ & $3$ & $96$\\
$B_{120}$ & $\langle2\rangle\oplus\langle-6\rangle\oplus\langle10\rangle$ & $\langle 3/4\rangle\oplus V(2)\oplus 3^{(1)}[2]\oplus5^{(1)}[1]$ & $10$ & yes & --- & $0$ & $2$ & $24$\\
$B_{184}$ & $\langle2\rangle\oplus\langle2\rangle\oplus\langle-46\rangle$ & $\langle 3/4\rangle\oplus V(2)\oplus 23^{(1)}[5]$ & $46$ & yes & --- & $0$ & $2$ & $12$\\
$B_{324}$ & $A_2\oplus\langle-108\rangle$ & $\langle 5/8\rangle\oplus 3^{(3,1)}[1,1]$ & $6$ & no & --- & $0$ & $2$ & $24$\\
$B_{648}$ & $A_2\oplus\langle-216\rangle$ & $\langle 5/16\rangle\oplus 3^{(3,1)}[2,1]$ & $1$ & no & --- & $0$ & $2$ & $24$\\
\bottomrule
\end{tabular}
\end{table}

\noindent Table~\ref{tab:twelve} lists the twelve with their invariants as lattices; their modular data are collected in
Table~\ref{tab:canon} at the beginning of \S\ref{sec:eight}, once they have been defined. On a first reading, what matters in Table~\ref{tab:twelve}
is the rule that separates the four theta lattices from the eight; the other columns are invariants of the lattices (\S\ref{subsec:lat}), and the
column of root bases is explained below. The list is finite because, for a rigid $D=D(L)$, two bounds enclose $\dim S_{3/2}(\rho_L)$: a lower bound from the dimension formula, which
grows like $|D|/48$, and an upper bound for the degree of an irreducible representation of $O(D)$ that can occur in $S_{3/2}(\rho_L)$, which grows like
$|D|^{1/2+o(1)}$. The proof of Theorem~A in \S\ref{sec:proofA} makes this comparison effective: no module with $|D|>10^4$ is rigid
(\S\ref{subsec:large}); the two bounds leave $258$ modules with $|D|\le10^4$ (\S\ref{subsec:small}); the character of $O(D)$ on $S_{3/2}(\rho_D)$ leaves $23$
rigid ones (\S\ref{subsec:chi}); and twelve of them are discriminant forms of even lattices of signature $(2,1)$, each genus with one class
(\S\S\ref{subsec:nik}--\ref{subsec:classno}). The steps that rest on computation are marked \gC; their scripts and logs are in the reproduction archive
(Appendix~\ref{app:comp}). The proof comes after the forms, in \S\ref{sec:proofA}; the other sections use only the statement of Theorem~A.

A \emph{root} of $L$ is a primitive $r\in L$ with $(r,r)>0$ and $2(r,x)/(r,r)\in\Z$ for all $x\in L$, and a \emph{root basis} is a $\Z$-basis of roots
with pairwise non-positive inner products. Its Cartan matrix $a_{ij}=2(r_i,r_j)/(r_i,r_i)$ is a symmetrizable generalized Cartan matrix; as $L$ has
signature $(2,1)$, it is \emph{hyperbolic} in the sense of Kac--Moody theory exactly when its Dynkin diagram is connected and $a_{ij}a_{ji}\le4$ for
$i\ne j$, and \emph{compact} when moreover $a_{ij}a_{ji}\le3$. The census of [II] lists the $44$ classes of symmetrizable hyperbolic Cartan matrices of rank three (its \emph{core}) and $56$
matrices one edit away from them (its \emph{rim}) [II, \S3.1] (\S\ref{subsec:roots}). In the names of the twelve (Table~\ref{tab:twelve}) the letter is $L$ when the
lattice has a Kac--Moody hyperbolic root basis and $B$ when it has none, and the index is $|\det L|$ (Proposition~\ref{prop:roots}); [I] and [II] write $L_4$, $B_{10}$,
$B_{22}$ for $L_{72}$, $B_{40}$, $B_{88}$. The six lattices with a hyperbolic root basis are rescaled root lattices: $L_{16},L_{54},L_{64},L_{72},L_{96},L_{250}$
are $M(k)$, the root lattice $M$ of a hyperbolic Cartan matrix of the census of [II] in its minimal symmetrization with its form multiplied by $k$, for
\[
(M,k)=(FF,2),\ (FF,3),\ (M_8,2),\ (L_{72},1),\ (M_{12},2),\ (FF,5)
\]
respectively, where $FF=U\oplus\langle2\rangle$ is the root lattice of the Feingold--Frenkel matrix [FF83], $M_{12}=A_2\oplus\langle-4\rangle$, and $M_8$ is the
census row with Gram matrix $\bigl(\begin{smallmatrix}2&-2&-1\\-2&2&-1\\-1&-1&2\end{smallmatrix}\bigr)$ (Proposition~\ref{prop:roots}(b)). Hyperbolicity is a property of
the Cartan matrix, which does not see $k$; rigidity is a property of the discriminant form, which does: $FF(k)$ is rigid for $k=2,3,5$ and not for
$k=1,4,6$ (Proposition~\ref{prop:roots}(c)), and for these $k$ the forms that \S\ref{sec:canon} attaches to $FF(k)$ come from umbral moonshine of lambency $k$
(Remark~\ref{rem:theta4}). Of the eight
lattices of this paper only $L_{72}$ and $L_{96}$ have a hyperbolic root basis, and theirs are compact; $B_{40}$ and $B_{88}$ are rows of the census of [II]
(its rim rows), and $B_{120},B_{184},B_{324},B_{648}$ are not of the form $M(k)$ for any row $M$ of the census (Proposition~\ref{prop:roots}(b)).

\subsection{The canonical forms}\label{subsec:nochoice}
In general a harmonic weak Maass form of weight $\frac12$ for the dual Weil representation is fixed by its principal part only up to holomorphic
modular forms, and to tie it to a given cusp form as its shadow one has to choose a normalization of that cusp form and of the Petersson product. On a
Petersson-rigid lattice none of this is a choice, once the pole is required to have the least possible order. Write $S(D)=S_{3/2}(\rho_L)$, $D=D(L)$. Rigidity means
\[
S(D)=\C F_0\otimes\sigma_0
\]
for one irreducible representation $\sigma_0$ of $O(D)$ and one form $F_0$ with values in the multiplicity space $\pi_{\sigma_0}$ of $\sigma_0$
(\S\ref{subsec:mult}). By Schur's lemma the Petersson product of $S(D)$ is unique up to a positive scalar, and the scalar is fixed by any one Fourier
coefficient: for every $n>0$ with $c_{F_0}(n)\ne0$ the coefficient pairing $B_n(F,G)=\sum_\gamma c_F(\gamma,n)\overline{c_G(\gamma,n)}$ is
$t_n^{-1}$ times the Petersson product, for one number $t_n>0$ (Proposition~\ref{prop:Bn}). Conversely, on a reducible space the Petersson product is
not determined by the symmetry (Remark~\ref{rem:converse}). Let $\mu\in(0,1]$ be the least $r$ for which $\e(r)$ is an eigenvalue of $\rho_L(T)$ on
$\pi_{\sigma_0}$, and let $V_{\sigma_0}=\pi_{\sigma_0}\otimes\sigma_0\subset\C[D]$ be the $\sigma_0$-isotypic part (\S\ref{subsec:mult}).

\begin{thmOne}[{Theorem~\ref{thm:harm}}]
Let $L$ be one of the twelve lattices of Theorem~A, and consider harmonic weak Maass forms of weight $\frac12$ for $\rho_L^*$ with values in $V_{\sigma_0}$.
\begin{enumerate}[label=(\roman*),leftmargin=7mm,itemsep=1pt]
\item Such a form with non-zero shadow has a pole of order at least $\mu$.
\item For each $F\in S(D)$ exactly one such form $\Hh_F$ has principal part $q^{-\mu}\,\overline{c_F(\mu)}$, and its shadow is
\[
\xi_{1/2}\Hh_F=F/t_\mu .
\]
\item The map $F\mapsto\Hh_F$ is an antilinear $O(D)$-isomorphism from $S(D)$ onto the space of such forms whose principal part is $q^{-\mu}x$ for a
vector $x$.
\end{enumerate}
\end{thmOne}

\noindent The proof is in \S\ref{sec:canon}; it uses of rigidity only that $\sigma_0$ occurs once in $S(D)$, and the hypotheses under which
Theorem~\ref{thm:harm} is stated hold for the twelve lattices \gC\ (\S\ref{subsec:mult}). So each lattice carries one harmonic weak
Maass form $\Hh=\Hh_{F_0}$, with principal part $q^{-\mu}\,\overline{c_{F_0}(\mu)}$ and shadow $F_0/t_\mu$ ($F_0$ realized in $S(D)$ as in
\S\ref{subsec:mult}), and one number $t_\mu$, the Petersson norm of $F_0$ normalized by its leading coefficient.
For the four theta lattices $t_\mu=\frac{s_L-1}{12}\sqrt{s_L}$, with $s_L$ the square-free part of $|\det L|/2$;
for them the same construction gives the forms of umbral moonshine, and the eigenvalues of $T(p^2)$ are $p+1$, those of an Eisenstein series
(Remark~\ref{rem:theta4}). For the eight, $t_\mu$ is computed numerically (\S\ref{subsec:table}), and at the primes computed the eigenvalues of
$T(p^2)$ are the coefficients $a_p(E)$ of the elliptic curve $E$ of Table~\ref{tab:canon} \gC\ (\S\ref{sec:arith}). Table~\ref{tab:canon} lists these data for all
twelve rigid lattices; the coefficients of $\Hh$ are in Appendix~\ref{app:tabH}.

\subsection{Results}\label{subsec:results}
\emph{The forms.} For the eight, the forms $\Hh$ of Theorem~1 are computed by
Hejhal's method from $F_0$ and $t_\mu$ (\S\ref{sec:eight}, Appendix~\ref{app:tabH}) \gC. For three of the eight, $F_0$ is explicit as an eta--theta product (Proposition~\ref{prop:etatheta}):
for $B_{40}$, $L_{72}$, $L_{96}$ the coordinate functions of $F_0$ are $\eta^2$ times theta functions of the lattices $\langle30\rangle$,
$\langle6\rangle$, $\langle24\rangle$ (differences of two for $\langle30\rangle$ and $\langle24\rangle$), and for the other five no identity of this shape
holds \gC\ (given the computed exponents).

\emph{Arithmetic} (\S\ref{sec:arith}). The Hecke operators $T(p^2)$, $p\nmid2|D|$, act on $S(D)$ by scalars, which are integers
(Propositions~\ref{prop:hecke}, \ref{prop:cm}(a)), and at the primes computed these are the coefficients $a_p(E)$ of the curves of
Table~\ref{tab:canon} \gC; among them only $36a1$ has complex multiplication, and only $L_{72}$ can carry an eigensystem with complex multiplication at
all (Proposition~\ref{prop:cm}) \gC. That the eigensystem is that of a newform of weight $2$ is not proved for vector-valued forms on these
discriminant forms (Problem~\ref{prob:hecke}); granting it, rationality and strong multiplicity one make it the eigensystem of an elliptic curve
over $\Q$ (Proposition~\ref{prop:cm}(b)).

\emph{Rigidity by counting} (\S\ref{sec:count}). Three of the eight, $B_{40},B_{88},B_{184}$, belong to the family
$\Dp_p=(\Z/2)^3_{\rm an}\oplus(\Z/p,ux^2/p)$, $p$ an odd prime.

\begin{thmTwo}[{Lemma~\ref{lem:three}, Corollary~\ref{cor:geom}}]
For every prime $p>3$,
\[
g\bigl(X_0(4p)/w_p\bigr)-g\bigl(X_0(2p)/w_p\bigr)=
\begin{cases}\frac{p-1}8+\frac{h(-4p)}4,&p\equiv1\pmod4,\\[2pt] \frac{p-3}8,&p\equiv3\pmod8,\\[2pt] \frac{p-3}8-\frac{h(-p)}2,&p\equiv7\pmod8.\end{cases}
\]
This number is $0,1,2,3$ exactly for $p\in\{7\},\{5,11,23\},\{13,19,31\},\{17,47\}$ among $p>3$, and, by comparison with the list of Theorem~A,
$\Dp_p$ is rigid exactly when $X_0(2p)/w_p$ has genus $0$ and $X_0(4p)/w_p$ genus $1$.
\end{thmTwo}

\noindent The proof uses the class numbers below $10^5$ \gC. The difference of genera is the number of newforms of weight $2$ and level $p$,
$2p$ or $4p$ with $w_p=+1$ (Lemma~\ref{lem:three}), and the count extends to every even weight (Proposition~\ref{prop:allw}). For $3\le p\le103$ and
$k\le\frac92$, $\dim S_k(\rho_{\Dp_p})$ is an explicit combination of newform counts of weight $2k-1$ \gC\ (Remark~\ref{rem:dims}); whether this holds
for every $p$ is open (Problem~\ref{prob:dict}), and nothing in this paper uses it. \S\ref{sec:conj} lists open problems.

\subsection{Relation to the literature}\label{subsec:litrel}
Vector-valued harmonic weak Maass forms of weight $\frac12$ (the term of [BO10]; weak Maass forms in [BF04], harmonic Maass forms in [Br13]) whose shadows are the
forms attached to weight-$2$ cusp forms by the Shimura correspondence are constructed in [BO10], for prime level, and in [Br13] for the Weil
representation of $\Z/2N$, whose curves are $X_0(N)$; there the coefficients at fundamental discriminants are algebraic exactly when a twisted
derivative $L'(G,\chi_\Delta,1)$ vanishes [BO10, Thm.~1.1], [Br13, Thm.~3.2], and [AGOR15] treats the Weierstrass mock modular forms of elliptic curves.
Hejhal's method for these forms is [BSt12], with tables for $\Z/2N$, $N=11,37$. The forms of this paper live on non-cyclic
discriminant forms; seven of the eight lattices have $\fd>1$, so $L\otimes\Q$ is anisotropic, its even Clifford algebra is a division algebra and the
Shimura curve of an order in it is compact [Vo21, Prop.~38.1.3, Main Thm.~38.4.3]. Hecke operators for the Weil representation in odd signature are [BS10]; their
eigenvalues are rational on a rational basis by McGraw's theorem [McG03], and on the rigid spaces here they are integers
(Proposition~\ref{prop:cm}). The counts of Atkin--Lehner signs behind \S\ref{sec:count} are classical:
the fixed points of $w_Q$ on $X_0(N)$ are counted in [Og74], traces of Atkin--Lehner operators in [Ya73], and refined dimensions with signs in
[Ma18]. The traces at the levels $p$, $2p$, $4p$ of Theorem~2 are given in closed form in [Ma25, \S4.1], and the groups $\Gamma_0(2n)+n$ of genus $0$
are listed in [CDH14b, \S2.3]; both are compared with our results in \S\ref{sec:count}. For the four theta lattices the forms $\Hh_F$ come from the
completed forms of umbral moonshine [ORTL15, Thm.~4.1], whose shadows are unary theta series, and the Petersson norms of these theta series are
computed in [BSc17, Example~5.2] (Remark~\ref{rem:theta4}).

\subsection{What is not in this paper}\label{subsec:notin}
Weight $\frac52$ enters only through the dimensions of Table~\ref{tab:canon}. The obstructions to Borcherds products, the corrections to the numbers of
[I] and [II] at weight $\frac52$, the shadows of the menus of [II] and the further weight-$\frac52$ computations (the level-$23$ block of $B_{184}$, the
$3$-adic descent of $B_{324}$ and $B_{648}$) are not treated; draft v0.9.8 of this paper, kept in the reproduction archive (folder \texttt{paper\_v098}),
discusses them, and nothing here depends on it. The discriminant-kernel curves of the five maximal lattices are not treated. No result of this paper depends on [I] or [II]; the comparison with their numbers is
Remark~\ref{rem:I-II}.
\section{Setting}\label{sec:setting}

\subsection{Finite quadratic modules}\label{subsec:fqm}
A \emph{finite quadratic module} is a finite abelian group $D$ with a map $q:D\to\Q/\Z$ such that $q(nx)=n^2q(x)$ and
$b(x,y):=q(x+y)-q(x)-q(y)$ is bilinear and nondegenerate. We write $\e(x)=e^{2\pi ix}$, $D_p$ for the $p$-primary part of $D$ and $D[m]=\{x:mx=0\}$.
The \emph{level} $N_D$ of $D$ is the least $N\ge1$ with $Nq(x)=0$ for all $x\in D$, and for $k\in\Z$ prime to $|D|$, $D(k)$ is $D$ with the form $kq$; in
particular $D(-1)=(D,-q)$.
The \emph{signature} $\sig D\in\Z/8$ is defined by $\sum_{x\in D}\e(q(x))=\sqrt{|D|}\,\e(\sig D/8)$ (Milgram's formula [MH73, App.~4] shows that the sum has
absolute value $\sqrt{|D|}$ and argument in $\frac{2\pi}8\Z$). $D$ is \emph{anisotropic} if $q(x)\ne0$ for $x\ne0$. For components we write, with $a$ odd,
$k\ge1$, odd $p$ and units $u_i$,
\[
\langle a/2^{k+1}\rangle=\bigl(\Z/2^k,\tfrac{ax^2}{2^{k+1}}\bigr),\qquad
V(2^k)=\bigl((\Z/2^k)^2,\tfrac{x^2+xy+y^2}{2^k}\bigr),
\]
\[
p^{(k_1,k_2,\dots)}[u_1,u_2,\dots]=\bigoplus_i\bigl(\Z/p^{k_i},\tfrac{u_ix^2}{p^{k_i}}\bigr),
\]
and $U(2^k)=((\Z/2^k)^2,xy/2^k)$. $(\Z/2)^3_{\rm an}$ denotes an anisotropic quadratic module of order $8$ and exponent $2$. Its values lie in
$\frac14\Z/\Z$, and it is $\langle a/4\rangle^3$ for one $a\in\{1,3\}$, the discriminant form of $\langle2\rangle^3$ resp.\ $\langle-2\rangle^3$: on
the seven non-zero elements $q$ takes the value $\frac a4$ three times, $\frac12$ three times and $\frac{3a}4$ once. (A sum
$\langle a/4\rangle\oplus\langle b/4\rangle\oplus\langle c/4\rangle$ is anisotropic iff $a\equiv b\equiv c\pmod4$, and $\langle a/4\rangle\oplus V(2)\cong\langle3a/4\rangle^3$:
the three elements $f+v$, $v\in V(2)\setminus0$, $f$ the generator of $\langle a/4\rangle$, are pairwise orthogonal of norm $\frac{3a}4$ and generate.)

\subsection{The Weil representation}\label{subsec:weil}
Let $\Mp$ be the metaplectic double cover of $\SL_2(\Z)$, generated by $T=\bigl(\bigl(\begin{smallmatrix}1&1\\0&1\end{smallmatrix}\bigr),1\bigr)$ and
$S=\bigl(\bigl(\begin{smallmatrix}0&-1\\1&0\end{smallmatrix}\bigr),\sqrt\tau\bigr)$, with $Z=S^2=(-I,i)$. The Weil representation $\rho_D$ of $\Mp$ on
$\C[D]=\bigoplus_x\C e_x$ is
\begin{equation}\label{eq:weil}
\rho_D(T)e_x=\e(q(x))\,e_x,\qquad \rho_D(S)e_x=\frac{\e(-\sig D/8)}{\sqrt{|D|}}\sum_{y\in D}\e(-b(x,y))\,e_y .
\end{equation}
It is unitary, and its dual $\rho_D^*$, which is its complex conjugate, is $\rho_{D(-1)}$; $\sig D(-1)=-\sig D$.
This is the normalization of Borcherds [Bo98] and Bruinier [Br02]; for $D=D(L)$ it is the representation computed by \texttt{WeilRep}$(-G)$ in the package
\texttt{weilrep} [Wi], and $S_{3/2}(\rho_{D(L)})$ is the space written $S_{3/2}(\rho_L)$ in [I] and [II]. [BEF16, (8)--(9)] use the dual normalization,
$\rho_A(T)e_\gamma=\e(-q(\gamma))e_\gamma$, so their $S_{k,A}$ is our $S_k(\rho_A^*)$. In particular [BEF16, Table~4] lists $D(B_{40})$ and $D(L_{72})$
(as $2_5^{+3}5^{+1}$ and $2_5^{+3}3^{+2}$, both of signature $1$ \gC) among the modules with $S_{3/2,A}=0$: this says $S_{3/2}(\rho_D^*)=0$, whereas
$S_{3/2}(\rho_D)$ has dimension $2$ and $4$ (Table~\ref{tab:twelve}). For $k\in\frac12\Z$ and a representation $\rho$ of
$\Mp$ on a finite-dimensional space $V$, $M_k(\rho)$ is the space of holomorphic $F:\mathbb H\to V$ with $F(M\tau)=\varphi(\tau)^{2k}\rho(M,\varphi)F(\tau)$
for $(M,\varphi)\in\Mp$, bounded at $i\infty$; $S_k(\rho)$ is the subspace of cusp forms. For $F$ with values in $\C[D]$ we write
$F=\sum_xF_xe_x=\sum_{x,n}c_F(x,n)\,q^ne_x$, $q=\e(\tau)$, and $c_F(n)=\sum_xc_F(x,n)e_x$.

The group $O(D)$ of isometries of $(D,q)$ acts on $\C[D]$ by $g\,e_x=e_{gx}$. It commutes with $\rho_D(T)$ and $\rho_D(S)$, hence with $\rho_D$, and
acts on $M_k(\rho_D)$ and $S_k(\rho_D)$ through the values. Let $\iota=-1\in O(D)$ and $V^\pm=\{v\in\C[D]:\iota v=\pm v\}$.

\begin{lemma}[parity]\label{lem:parity}
Every $F\in M_k(\rho_D)$ satisfies $F_{-x}=i^{\,\sig D-2k}F_x$. If $\sig D\equiv1\pmod8$, then $M_{3/2}(\rho_D)$ and $M_{1/2}(\rho_D^*)$ take values in
$V^-$.
\end{lemma}
\begin{proof}
Since $\sum_y\e(-b(x+z,y))=|D|\,[z=-x]$, squaring \eqref{eq:weil} gives $\rho_D(Z)e_x=i^{-\sig D}e_{-x}$. As $Z$ fixes $\tau$ with automorphy factor $i$, the transformation law at $Z$ reads $F=i^{2k}\rho_D(Z)F$, i.e.\ $F_x=i^{2k-\sig D}F_{-x}$. For $(k,\sig)=(\frac32,1)$ the factor is
$-1$; and $\rho_D^*=\rho_{D(-1)}$ has signature $-1$, so for $k=\frac12$ it is $i^{2}=-1$.
\end{proof}
\noindent Applying the lemma twice gives $F=(-1)^{\sig D-2k}F$; so $M_k(\rho_D)=0$ unless $2k\equiv\sig D\pmod2$.

\begin{lemma}[reality]\label{lem:real}
For $F\in M_k(\rho_D)$ put $F^c(\tau)=\overline{F(-\bar\tau)}$, the bar acting on the coordinates. Then $F\mapsto F^c$ is an antilinear involution of
$M_k(\rho_D)$, preserving $S_k(\rho_D)$ and commuting with $O(D)$, and $c_{F^c}(x,n)=\overline{c_F(x,n)}$. Hence the character of $O(D)$ on
$S_k(\rho_D)$ is real-valued.
\end{lemma}
\begin{proof}
$F^c$ is holomorphic, with the stated coefficients; so it has the same exponents as $F$, which gives the law at $T$, and it is bounded, resp.\ cuspidal,
with $F$. At $S$: for $\tau\in\mathbb H$ put $w=-\bar\tau\in\mathbb H$, so that $-1/w=1/\bar\tau$; the law of $F$ at $w$ reads
$F(1/\bar\tau)=w^{k}\rho_D(S)F(w)$, principal branch. Conjugating, with $\overline{w^k}=(-\tau)^k=\e(-\frac k2)\tau^k$ and
$\overline{\rho_D(S)}=\rho_D(S)\rho_D(Z)^{-1}$ (from \eqref{eq:weil} and $\rho_D(Z)e_x=i^{-\sig D}e_{-x}$), gives
$F^c(-1/\tau)=\e(-\frac k2)\tau^k\rho_D(S)\rho_D(Z)^{-1}F^c(\tau)$. By Lemma~\ref{lem:parity} for $F$, $F^c_{-x}=i^{2k-\sig D}F^c_x$, hence
$\rho_D(Z)F^c=i^{2k-2\sig D}F^c$, and the factor $\e(-\frac k2)\,i^{2\sig D-2k}=i^{2\sig D-4k}$ is $1$ because $2k\equiv\sig D\pmod2$ when $F\ne0$. So $F^c$
satisfies the law at the generators $T$, $S$, hence at all of $\Mp$. $O(D)$ acts by permutation matrices, which commute with the conjugation. Finally, an
antilinear bijection $c$ commuting with $g$ gives $\tr(g)=\tr(cgc^{-1})=\overline{\tr(g)}$.
\end{proof}

\subsection{Lattices}\label{subsec:lat}
$L$ is an even lattice with Gram matrix $G$, $(x,y)=x^{\mathsf T}Gy$, of signature $(2,1)$. Its discriminant form is $D(L)=L'/L$ with
$q(\gamma)=(\gamma,\gamma)/2\bmod1$; $|D(L)|=|\det G|$, and $\sig D(L)=1$ by Milgram's formula. $L(k)$ is $L$ with its form multiplied by $k$; $U$ is the
hyperbolic plane, $A_2$ the root lattice with Gram matrix $\bigl(\begin{smallmatrix}2&-1\\-1&2\end{smallmatrix}\bigr)$, $\langle n\rangle$ the lattice $\Z v$
with $(v,v)=n$.

\begin{cited}{Nikulin [Ni79, Cor.~1.9.4]}
The genus of an even lattice is determined by its signature and its discriminant form.
\end{cited}
\why{the passage from discriminant forms to genera. Consequence: Petersson-rigidity, defined through $D(L)$, is a property of the genus, and for the
twelve, whose genera have class number one (Theorem~A, \S\ref{subsec:classno}), the discriminant form determines the lattice.}

\noindent\textbf{The quaternion algebra.} For $L$ of rank three the even Clifford algebra $C^0(L\otimes\Q)$ is a quaternion algebra over $\Q$, ramified at
a place $v$ if and only if $L\otimes\Q_v$ is anisotropic [Vo21, Prop.~5.1.2, 5.3.23]. The \emph{quaternion discriminant} $\fd(L)$ is the
product of the ramified primes; [II] writes $D$ for it. As $L$ is indefinite, $\fd(L)$ is a product of an even number of primes, $\fd(L)=1$ if and
only if $L$ is isotropic over $\Q$, and $\fd(L(k))=\fd(L)$. A prime not dividing $2|D(L)|$ does not divide $\fd(L)$: there $L\otimes\Z_\ell$ is
unimodular of rank three, hence isotropic.

\noindent\textbf{Maximal lattices.} If $H\subseteq D(L)$ is a subgroup with $q(H)=0$, its preimage in $L'$ is an even lattice containing $L$, and every
even lattice $\tilde L\supseteq L$ in $L\otimes\Q$ arises in this way. So $L$ is maximal among even lattices in $L\otimes\Q$ if and only if $D(L)$ is
anisotropic. Among the eight, $B_{40},L_{72},B_{88},B_{120},B_{184}$ are maximal; $L_{96},B_{324},B_{648}$ are not.

\noindent\textbf{The lattice $L_{72}$.} Let $B=(\frac{5,6}{\Q})$ be the quaternion algebra of discriminant $6$, with $\alpha^2=5$, $\beta^2=6$,
$\beta\alpha=-\alpha\beta$, and $\mathcal O=\Z+\frac{1+\alpha}2\Z+\frac{\alpha+\alpha\beta}5\Z+\frac{5+\alpha+5\beta+\alpha\beta}{10}\Z$ the maximal order of [Err11, (7.1)]. The lattice
$\mathcal O^0=\{x\in\mathcal O:\operatorname{tr}x=0\}$ with the form $\operatorname{nrd}$ is the lattice of [Err11, \S3.1], and $|D(\mathcal O^0)|=8\cdot3^2=72$
[Err11, Prop.~3.1]. A change of basis takes $\mathcal O^0(-1)$ to $A_2(2)\oplus\langle-6\rangle$, so $L_{72}\cong\mathcal O^0(-1)$ \gC. [Err11, \S7]
uses this lattice for Borcherds products on the Shimura curve of discriminant $6$.

\subsection{Rigidity and the form \texorpdfstring{$F_0$}{F0}}\label{subsec:mult}
We write $S(D)=S_{3/2}(\rho_D)$ and $s=\dim S(D)$.

\begin{definition}[{[II, \S5.1]}]\label{def:rigid}
A finite quadratic module $D$ with $\sig D\equiv1\pmod8$ is \emph{rigid} if $O(D)$ acts absolutely irreducibly on $S(D)$ and $s\ge2$. An even lattice $L$
of signature $(2,1)$ is \emph{Petersson-rigid} if $D(L)$ is rigid.
\end{definition}

\noindent\textbf{The multiplicity spaces.} For each irreducible representation $\sigma$ of $O(D)$ the space
$\pi_\sigma=\Hom_{O(D)}(\sigma,V^-)$ is a representation of $\Mp$, and $V^-=\bigoplus_\sigma\pi_\sigma\otimes\sigma$. A function with values in $\pi_\sigma\otimes\sigma$ that transforms under $\Mp$
by $\pi_\sigma\otimes1$ is an element of $S_{3/2}(\pi_\sigma)\otimes\sigma$, so by Lemma~\ref{lem:parity}, as $O(D)$-modules,
\[
S(D)\cong\bigoplus_\sigma S_{3/2}(\pi_\sigma)\otimes\sigma .
\]
Hence $D$ is rigid if and only if $S(D)=\C F_0\otimes\sigma_0$ for one $\sigma_0$ of degree at least two and one form $F_0\in S_{3/2}(\pi_{\sigma_0})$, unique up
to a scalar. Beyond $\sigma_0$ and the inclusion $\pi_{\sigma_0}\otimes\sigma_0\subset V^-$, which are fixed by $D$, the space $S(D)$ is described
by the single form $F_0$: the Fourier coefficients of $F_0\otimes v$, $v\in\sigma_0$, are $c_{F_0}(n)\otimes v$ with $c_{F_0}(n)\in\pi_{\sigma_0}$.
We write $\mu$ for the least $T$-exponent of $\pi_{\sigma_0}$, taken in $(0,1]$; $\pi_{\sigma_0}(r)$ is the eigenspace of $\rho(T)$ with eigenvalue
$\e(r)$. We write $\chi_S$ for the character of $O(D)$ on $S(D)$; for rigid $D$ it is $\chi_{\sigma_0}$, which is real by Lemma~\ref{lem:real}.

\noindent\textbf{Realization.} Fix $v_0\in\sigma_0$, $v_0\ne0$. As $\sigma_0$ is irreducible, $\phi\mapsto\phi(v_0)$ is an injective $\Mp$-map
$\pi_{\sigma_0}\to V^-$. We identify $\pi_{\sigma_0}$ with its image $\pi_{\sigma_0}\otimes v_0\subset\C[D]$ and $F_0$ with $F_0\otimes v_0\in S(D)$; then the
coefficients $c_{F_0}(n)$ are vectors of $\C[D]$, and the inner product of $\C[D]$ restricts to an $\Mp$-invariant inner product on $\pi_{\sigma_0}$. Nothing
below depends on $v_0$: the numbers $t_n$ of \S\ref{subsec:onenumber} are ratios of such inner products, and the coordinates of the form $\Hh$ are
independent of $v_0$ (\S\ref{subsec:harm}).

\noindent\textbf{The twelve} \gC. Among the finite quadratic modules with $\sig D\equiv1\pmod8$ and at most three generators in each $p$-primary part,
twenty-three are rigid (\S\ref{subsec:chi}); twelve of them are the discriminant forms of the twelve lattices of Theorem~A. On
all twenty-three, $\chi_S$ is rational-valued, and $\sigma_0$ occurs neither in the odd part of $\C[\{x\in D:q(x)=0\}]$ nor in $M_{1/2}(\rho_D^*)$; we write
this $d_0(\sigma_0)=m_{1/2}(\sigma_0)=0$. For the twelve, $\dim\pi_{\sigma_0}$ and $\mu$ are as in Table~\ref{tab:canon}; $\mu$ has multiplicity one in
$\pi_{\sigma_0}$, and $c_{F_0}(\mu)\ne0$. The span of the coefficients $c_{F_0}(n)$ is $\Mp$-stable (it is the span of the values of $F_0$), and it is all
of $\pi_{\sigma_0}$ \gC.

\noindent\textbf{Coordinates.} For each class $r\bmod1$ of $T$-exponents of $\pi_{\sigma_0}$ let $e_{r,1}<e_{r,2}<\cdots$ be the exponents $e\equiv r$ at which
$c_{F_0}(e)\in\pi_{\sigma_0}(r)$ is not a linear combination of the $c_{F_0}(e')$, $e'\equiv r$, $e'<e$, and put $u_{r,i}=c_{F_0}(e_{r,i})$. As the coefficients
of $F_0$ span $\pi_{\sigma_0}$, the $u_{r,i}$ are a basis of $\pi_{\sigma_0}$, and
\[
F_0=\sum_{r,i}f_{r,i}\,u_{r,i},\qquad f_{r,i}=q^{e_{r,i}}\bigl(1+O(q)\bigr);
\]
the $e_{r,i}$ are distinct, the vectors $u_{r,i}$ of different classes are orthogonal, and $\mu=e_{r_0,1}$ for one class $r_0$. Replacing $F_0$ by
$\lambda F_0$ replaces $u_{r,i}$ by $\lambda u_{r,i}$ and leaves the $f_{r,i}$ unchanged.

\subsection{Conventions on newforms}\label{subsec:newforms}
Let $w\ge2$ be even and $N\ge1$. A \emph{newform} of weight $w$ and level $N$ is a normalized ($a_1=1$) eigenform of the Hecke operators $T_q$, $q\nmid N$,
in $S_w(\Gamma_0(N))$ with trivial character that is orthogonal to the old forms $g(d\tau)$, $g\in S_w(\Gamma_0(M))$, $M\mid N$, $M<N$, $dM\mid N$ [AL70].
\begin{enumerate}[label=(\arabic*),leftmargin=7mm,itemsep=1pt]
\item Its \emph{eigensystem} is $(a_q(f))_{q\nmid N}$. By strong multiplicity one a newform is determined by its eigensystem, even by the $a_q(f)$ for
all but finitely many $q$ [AL70]; so newforms of different levels have different eigensystems. Galois-conjugate newforms are different newforms
and are counted separately.
\item A newform of weight $2$ with rational eigenvalues is the form $f_E$ of an isogeny class of elliptic curves $E$ over $\Q$ of conductor $N$, with
$a_q(f_E)=a_q(E)$ (Eichler--Shimura; every isogeny class arises in this way [BCDT01]). We name $E$ by its label in [Cr97], the first curve of its isogeny class ($20a1$), and newforms by their labels in
[LMFDB] ($20.2.\mathrm{a.a}$ is $f_E$ for $E=20a1$).
\item On $S_k(\rho_D)$, $k\in\frac12+\Z$, the Hecke operators $T(q^2)$, $q\nmid2|D|$, are normalized as in \S\ref{subsec:hecke} (where the prime is
written $p$), so that the eigenvalue on the image of $f$ under the Shimura correspondence, where that is known, is $a_q(f)$ for $f$ of weight $2k-1$.
\end{enumerate}
We write \emph{eigensystem} where only Hecke eigenvalues matter, \emph{newform} where the level or a sign matters, and never ``modular form of
level $N$'' for a newform. Atkin--Lehner signs are used only in \S\ref{sec:count} and are fixed there (\S\ref{subsec:fam}).

\subsection{Hecke operators}\label{subsec:hecke}
\begin{cited}{Bruinier--Stein [BS10, Thm.~4.10]}
Let $A$ be a finite quadratic module of odd signature and level $N_A$, $k\in\frac12+\Z$, and $p$ a prime coprime to $N_A$. There is a Hecke operator
$T(p^2)$ on $M_k(\rho_A)$, preserving $S_k(\rho_A)$, which maps $F=\sum_{\lambda,n}c(\lambda,n)\e(n\tau)e_\lambda$ to $\sum_{\lambda,n}b(\lambda,n)\e(n\tau)e_\lambda$
with
\[
b(\lambda,n)=c(p\lambda,p^2n)+\kappa_{\rm BS}(p,n)\,c(\lambda,n)+p^{2k-2}c(\lambda/p,n/p^2),
\]
where $\kappa_{\rm BS}(p,n)=\epsilon_p^{\sig(A)+(\frac{-1}{|A|})}\bigl(\frac p{|A|\,2^{\sig(A)}}\bigr)p^{k-3/2}\bigl(\frac{-n}p\bigr)$ with $\epsilon_p=1$ or $i$ as
$p\equiv1$ or $3\pmod4$, $\lambda/p$ is the preimage of $\lambda$ under multiplication by $p$, and $c(\lambda/p,n/p^2)=0$ if $p^2\nmid n$.
\end{cited}
\why{the construction of Hecke operators for the Weil representation in odd signature. The normalization of $\rho_A$ in [BS10, (2.3)--(2.4)] is that of
\eqref{eq:weil}; the symbols are read as Kronecker symbols, as in [BS10, Rem.~4.11]. [BS10] writes $\kappa(p,n)$ for $\kappa_{\rm BS}(p,n)$.}

\begin{proposition}\label{prop:hecke}
Let $D$ be rigid and $p\nmid2|D|$ prime. Then $T(p^2)$ commutes with $O(D)$ and acts on $S(D)$ by a scalar $\lambda_p$, so $S(D)$ carries a single
eigensystem $(\lambda_p)_{p\nmid2|D|}$. If $\sig D\equiv1\pmod8$ and $k=\frac32$, then $\kappa_{\rm BS}(p,n)=\bigl(\frac{2|D|n}p\bigr)$, i.e.
\begin{equation}\label{eq:T32}
\bigl(T(p^2)F\bigr)(\lambda,n)=c(p\lambda,p^2n)+\Bigl(\frac{2|D|n}p\Bigr)c(\lambda,n)+p\,c(\lambda/p,n/p^2).
\end{equation}
\end{proposition}
\begin{proof}
Let $g\in O(D)$. Then $c_{gF}(\lambda,n)=c_F(g^{-1}\lambda,n)$, and as $g$ is a group automorphism, $g^{-1}(p\lambda)=p\,g^{-1}\lambda$ and
$g^{-1}(\lambda/p)=(g^{-1}\lambda)/p$; so $b_{gF}(\lambda,n)=b_F(g^{-1}\lambda,n)$, i.e.\ $T(p^2)(gF)=g\,T(p^2)F$. As $S(D)$ is absolutely irreducible,
$T(p^2)$ is a scalar on it by Schur's lemma [Se77, \S2.2, Prop.~4].

For the formula let $D_{\rm odd}$ be the odd part, so that $|D|=|D_2|\,|D_{\rm odd}|$, and read the symbols as Kronecker symbols, so that
$(\frac{-1}{|D|})=(\frac{-1}{|D_{\rm odd}|})$. Then $\epsilon_p^{1+(\frac{-1}{|D|})}$ is $(\frac{-1}p)$ if $|D_{\rm odd}|\equiv1$ and $1$ if
$|D_{\rm odd}|\equiv3\pmod4$, while $(\frac p{2|D|})=(\frac{2|D_2|}p)(\frac p{|D_{\rm odd}|})$ and
$(\frac p{|D_{\rm odd}|})=(\frac{-1}p)^{(|D_{\rm odd}|-1)/2}(\frac{|D_{\rm odd}|}p)$. In both cases the product is $(\frac{2|D_2|}p)(\frac{-|D_{\rm odd}|}p)$, and
$\kappa_{\rm BS}(p,n)=(\frac{2|D_2|}p)(\frac{-|D_{\rm odd}|}p)(\frac{-n}p)=(\frac{2|D|n}p)$.
\end{proof}
\noindent For the eight lattices all relations \eqref{eq:T32} computed agree with this value of $\kappa_{\rm BS}$, and each of the symbols $(\frac{-2|D|n}p)$,
$(\frac{\pm|D|n}p)$ fails on some of them \gC.
\section{The canonical scalar and the form \texorpdfstring{$\mathcal H$}{H}}\label{sec:canon}

Throughout this section $D$ is rigid, $S(D)=\C F_0\otimes\sigma_0$, $V_{\sigma_0}=\pi_{\sigma_0}\otimes\sigma_0\subset V^-$ is the $\sigma_0$-isotypic part,
and $P(F,G)=\int_{\SL_2(\Z)\backslash\mathbb H}\langle F(\tau),G(\tau)\rangle\,y^{3/2}\,\frac{dx\,dy}{y^2}$ is the Petersson product on $S(D)$, with
$\langle\ ,\ \rangle$ the standard inner product of $\C[D]$, which $\rho_D$ and $O(D)$ preserve.

\subsection{One number}\label{subsec:onenumber}
The inner product of $\C[D]$, restricted to $V_{\sigma_0}$, is invariant under $\Mp\times O(D)$. As $\sigma_0$ is irreducible, it is
$\langle\ ,\ \rangle_\pi\otimes\langle\ ,\ \rangle_\sigma$ with an $\Mp$-invariant inner product on $\pi_{\sigma_0}$ and an $O(D)$-invariant one on
$\sigma_0$, the pair being unique up to $(c\langle\ ,\ \rangle_\pi,c^{-1}\langle\ ,\ \rangle_\sigma)$, $c>0$. Hence
\[
P(F_0\otimes v,F_0\otimes w)=\|F_0\|_\pi^2\,\langle v,w\rangle_\sigma,\qquad \|F_0\|_\pi^2=\int_{\SL_2(\Z)\backslash\mathbb H}\|F_0(\tau)\|_\pi^2\,y^{3/2}\,\frac{dx\,dy}{y^2}.
\]
For $n>0$ put $B_n(F,G)=\sum_{\gamma\in D}c_F(\gamma,n)\overline{c_G(\gamma,n)}$.

\begin{proposition}\label{prop:Bn}
For every $n>0$ with $c_{F_0}(n)\ne0$ there is a number $t_n>0$ with $P=t_n\,B_n$ on $S(D)$, namely
\[
t_n=\|F_0\|_\pi^2\big/\|c_{F_0}(n)\|_\pi^2 ,
\]
the Petersson norm of $F_0$ normalized by $\|c_{F_0}(n)\|_\pi=1$. It depends neither on the normalization of $F_0$ nor on the choice of
$\langle\ ,\ \rangle_\pi$.
\end{proposition}
\begin{proof}
$B_n$ is $O(D)$-invariant, as $c_{gF}(\gamma,n)=c_F(g^{-1}\gamma,n)$ and $g$ permutes the positions, and positive semi-definite. Its kernel is $O(D)$-stable,
hence $0$ if $B_n\ne0$; so by Schur's lemma $B_n$ is a positive multiple of $P$. The $n$-th coefficient of $F_0\otimes v$ is $c_{F_0}(n)\otimes v$,
and the inclusion $V_{\sigma_0}\subset\C[D]$ is isometric for $\langle\ ,\ \rangle_\pi\otimes\langle\ ,\ \rangle_\sigma$; so
$B_n(F_0\otimes v,F_0\otimes w)=\|c_{F_0}(n)\|_\pi^2\langle v,w\rangle_\sigma$, and $t_n=\|F_0\|_\pi^2/\|c_{F_0}(n)\|_\pi^2$.
\end{proof}

\begin{remark}[the converse]\label{rem:converse}
If $S(D)$ is not irreducible, the $O(D)$-invariant Hermitian forms on $S(D)$ do not form a single ray: for $S(D)=S_1\oplus S_2$ with $S_1,S_2$ invariant and
orthogonal, $P|_{S_1}\oplus2P|_{S_2}$ is invariant and not proportional to $P$. So ``the Petersson product is fixed by the symmetry'' is equivalent to
irreducibility, and the numbers $t_n$ are attached to rigid modules. They exist slice by slice without rigidity (\S\ref{subsec:slices}).
\end{remark}

\subsection{The form \texorpdfstring{$\mathcal H$}{H}}\label{subsec:harm}
For a smooth function $f$ on $\mathbb H$ with values in $\C^m$ put $\xi_kf=2iy^k\,\overline{\partial f/\partial\bar\tau}$, the bar acting on the
coordinates; it is antilinear, and it maps functions that transform with weight $k$ and a representation $\rho$ to functions that transform with weight
$2-k$ and $\bar\rho$. By a harmonic weak Maass form, in the terminology of [BO10], we mean an element of $H^+$ in [BF04, \S3]: a weak Maass form in the sense of [BF04]
whose image under $\xi_{1/2}$ is a cusp form (a harmonic Maass form in the terminology of [Br13]); its image under $\xi_{1/2}$ is its \emph{shadow}.

\begin{cited}{Bruinier--Funke [BF04, \S3: Prop.~3.5, Cor.~3.8, Prop.~3.11]}
Let $L$ be an even lattice and $k\in\frac12\Z$, $k\ne1$. The weak Maass forms $f\in H^+_{k,L}$ transform with $\rho_L$, and $\xi_k$ maps $H^+_{k,L}$, the inverse
image of the cusp forms, onto the space $S_{2-k,L^-}$ of cusp forms for $\rho_{L^-}=\overline{\rho_L}$, with kernel $M^!_{k,L}$. For $g\in S_{2-k,L^-}$ and
$f\in H^+_{k,L}$, $(g,\xi_kf)=\sum_{\gamma\in L'/L}\sum_{n\le0}c_f^+(\gamma,n)\,c_g(\gamma,-n)$. For every Fourier polynomial $Q$ of the shape of a principal part
there is $f\in H^+_{k,L}$ whose principal part is $Q+c_0$, with $c_0$ a $T$-invariant constant.
\end{cited}
\why{the analytic theory of harmonic weak Maass forms (pairing with cusp forms, surjectivity of $\xi$, prescribed principal parts), used as stated.}

\begin{theorem}\label{thm:harm}
Let $D$ be rigid with $d_0(\sigma_0)=m_{1/2}(\sigma_0)=0$ (\S\ref{subsec:mult}).
\begin{enumerate}[label=(\alph*),leftmargin=7mm,itemsep=2pt]
\item For $F\in S(D)$ and $n_0>0$ with $c_F(n_0)\ne0$ there is exactly one harmonic weak Maass form $\Hh_F$ of weight $\frac12$ for $\rho_D^*$ with values
in $V_{\sigma_0}$ whose principal part is $q^{-n_0}\,\overline{c_F(n_0)}$, and $\xi_{1/2}\Hh_F=F/t_{n_0}$.
\item A harmonic weak Maass form of weight $\frac12$ for $\rho_D^*$ with values in $V_{\sigma_0}$ and non-zero shadow has a pole of order at least $\mu$.
If $\mu$ has multiplicity one in $\pi_{\sigma_0}$ and $c_{F_0}(\mu)\ne0$, then $F\mapsto\Hh_F$ (with $n_0=\mu$) is an antilinear $O(D)$-isomorphism from
$S(D)$ onto the space of harmonic weak Maass forms of weight $\frac12$ for $\rho_D^*$ with values in $V_{\sigma_0}$ whose principal part is $q^{-\mu}x$ for
a vector $x$, and $\xi_{1/2}\Hh_F=t_\mu^{-1}F$ for every $F\in S(D)$.
\end{enumerate}
The hypotheses of (a) and (b) hold for the twelve lattices of Theorem~A \gC\ (\S\ref{subsec:mult}).
\end{theorem}
\begin{proof}
We apply the quoted statement to an even lattice $L$ with $D(L)\cong D(-1)$, so that $\rho_L=\rho_D^*$ and $S_{3/2,L^-}=S(D)$: for $D=D(L_0)$ take $L=L_0(-1)$.
\emph{Projections.} $O(D)$ and $\iota$ act by real permutation matrices, which commute with $\Mp$ and with $\xi_{1/2}$, and $V_{\sigma_0}$ is the image of
$V^-$ under the central idempotent $e_{\sigma_0}=\frac{\deg\sigma_0}{|O(D)|}\sum_g\chi_{\sigma_0}(g)\,g$, which has real coefficients, $\chi_{\sigma_0}$ being
real (\S\ref{subsec:mult}). So $e_{\sigma_0}\frac{1-\iota}2$ commutes with $\xi_{1/2}$ and with complex conjugation of the coordinates, maps harmonic weak Maass forms
for $\rho_D^*$ to harmonic weak Maass forms with values in $V_{\sigma_0}$, keeps a principal part with values in $V_{\sigma_0}$, and $V_{\sigma_0}$ is stable under
complex conjugation of the coordinates. As $d_0(\sigma_0)=0$, no $T$-invariant vector lies in $V_{\sigma_0}$.

(a) \emph{Existence.} There is a harmonic weak Maass form $f$ for $\rho_D^*$ with principal part $q^{-n_0}\overline{c_F(n_0)}+c_0$, $c_0$ a $T$-invariant
constant; $\Hh_F:=e_{\sigma_0}\frac{1-\iota}2f$ has principal part $q^{-n_0}\overline{c_F(n_0)}$. \emph{Uniqueness.} The difference $\Delta$ of two such forms
has values in $V_{\sigma_0}$ and principal part $0$. Its shadow lies in $S(D)$, and $(G,\xi_{1/2}\Delta)=0$ for every $G\in S(D)$; so $\xi_{1/2}\Delta=0$,
and $\Delta\in M_{1/2}(\rho_D^*)$ has values in $V_{\sigma_0}$. That is the $\sigma_0$-isotypic part of $M_{1/2}(\rho_D^*)$, which is $0$ as
$m_{1/2}(\sigma_0)=0$. \emph{The shadow.} For $G\in S(D)$, by the pairing and Proposition~\ref{prop:Bn},
\[
(G,\xi_{1/2}\Hh_F)=\sum_{\gamma\in D}\overline{c_F(\gamma,n_0)}\,c_G(\gamma,n_0)=B_{n_0}(G,F)=t_{n_0}^{-1}P(G,F)=(G,F/t_{n_0}),
\]
so $\xi_{1/2}\Hh_F-F/t_{n_0}\in S(D)$ is orthogonal to $S(D)$, and it is $0$.

(b) If $\Phi$ is a harmonic weak Maass form with values in $V_{\sigma_0}$ and $\xi_{1/2}\Phi\ne0$, then
$0\ne\|\xi_{1/2}\Phi\|^2=\sum_{n>0}\sum_\gamma c^+_{\Phi}(\gamma,-n)\,c_{\xi\Phi}(\gamma,n)$, and $c_{\xi\Phi}(\gamma,n)=0$ unless $n\equiv r\pmod1$ for a
$T$-exponent $r\in(0,1]$ of $\pi_{\sigma_0}$, and then $n\ge r\ge\mu$; so $\Phi$ has a pole of order at least $\mu$. Let $\mu$ have
multiplicity one and $c_{F_0}(\mu)\ne0$. The vectors of $V_{\sigma_0}$ with $T$-exponent $\mu$ form $\pi_{\sigma_0}(\mu)\otimes\sigma_0\cong\sigma_0$, and
$F\mapsto\overline{c_F(\mu)}$ is an antilinear $O(D)$-map from $S(D)\cong\sigma_0$ to it, non-zero since $c_{F_0\otimes v}(\mu)=c_{F_0}(\mu)\otimes v$; by Schur's
lemma it is bijective. So every principal part $q^{-\mu}x$, $x\in\pi_{\sigma_0}(\mu)\otimes\sigma_0$, is $q^{-\mu}\overline{c_F(\mu)}$ for exactly one $F$, and (a)
gives the isomorphism and $\xi_{1/2}\Hh_F=F/t_\mu$.
\end{proof}
\noindent In (a) the principal part is prescribed instead of the norm of the shadow, and the constant that appears is the canonical scalar of
Proposition~\ref{prop:Bn}. Write $\Hh=\Hh_{F_0}$, with $F_0$ realized as $F_0\otimes v_0$ (\S\ref{subsec:mult}): the form with principal part
$q^{-\mu}\overline{c_{F_0}(\mu)}$ and shadow $F_0/t_\mu$. Replacing $F_0$ by $\lambda F_0$ replaces $\Hh$ by $\bar\lambda\Hh$, and $\Hh$ does not depend on $v_0$ in
the following sense. As $F\mapsto\Hh_F$ is antilinear and $O(D)$-equivariant, $v\mapsto\overline{\Hh_{F_0\otimes v}}$ (conjugation of the coordinates) is linear and
$O(D)$-equivariant from $\sigma_0$ to functions with values in $V^-$; so $\overline{\Hh_{F_0\otimes v}}=\Psi\otimes v$ for one function $\Psi$ with values in
$\pi_{\sigma_0}$, independent of $v$. In the basis $u_{r,i}$ of \S\ref{subsec:mult}, $\Psi=\sum\psi_{r,i}u_{r,i}$ and
$\Hh=\sum_{r,i}\overline{\psi_{r,i}}\;\overline{u_{r,i}}$, so the coordinate functions $h_{r,i}=\overline{\psi_{r,i}}$ of $\Hh$ do not depend on $v_0$.

\begin{remark}[the mean of $1/t_n$]\label{rem:mean}
$\sum_{n\le X}1/t_n\sim32X^{3/2}$ for every rigid $D$: unfolding $\int\|F_0\|_\pi^2y^{3/2}E(\tau,s)\,\frac{dx\,dy}{y^2}$ against the Eisenstein series
$E(\tau,s)=\sum_{\Gamma_\infty\backslash\SL_2(\Z)}\operatorname{Im}(\gamma\tau)^s$, whose residue at $s=1$ is $3/\pi$ [Bu97, \S1.6], gives
$\Gamma(s+\frac12)(4\pi)^{-s-1/2}\sum_n\|c_{F_0}(n)\|_\pi^2n^{-s-1/2}$, so $\sum_n t_n^{-1}n^{-w}$ has a simple pole at $w=\frac32$ with residue $48$, and the
Wiener--Ikehara theorem applies.
\end{remark}

\subsection{Slices without rigidity}\label{subsec:slices}
The proofs of Propositions~\ref{prop:hecke}, \ref{prop:Bn} and of Theorem~\ref{thm:harm} use of $\sigma_0$ only that it occurs in $S(D)$ with multiplicity
one. So for every irreducible $\sigma$ with $\dim S_{3/2}(\pi_\sigma)=1$ (a \emph{rigid slice}) there are a form $F_\sigma$, numbers $t^{(\sigma)}_n$, an
eigensystem, and, if $d_0(\sigma)=m_{1/2}(\sigma)=0$, a form $\Hh_\sigma$; if $\chi_\sigma$ is not real, the principal part
$\overline{c_F(n_0)}$ lies in $V_{\bar\sigma}$, and $\Hh_\sigma$ takes its values in $V_{\bar\sigma}$ (Lemma~\ref{lem:real} gives
$\dim S_{3/2}(\pi_{\bar\sigma})=\dim S_{3/2}(\pi_\sigma)$). The eigensystem is rational where $\chi_\sigma$ is rational: the projection $e_\sigma$ then has
rational coefficients, and Proposition~\ref{prop:cm}(a) applies to $e_\sigma F$ for $F$ with rational coefficients. Rigidity is the case of a single slice. If $S(D)$ has two slices, the coefficient pairing $B_n$ need not be a multiple of $P$: in the census row $72$ of [II] ($|D|=108$, $\fd=6$,
two slices of dimension $2$, whose eigenvalues at the primes computed are those of the curves $54a1$ and $54b1$), at $n=\frac7{72}$ the two slices
give $t=1.42118834\ldots$ and $0.23194890\ldots$ \gC. Among
the $278$ modules on which \S\ref{subsec:chi} computes the character, $82$ are lattice modules with $s\ge2$; of these, $12$ are rigid and $45$ are sums of two or three
rigid slices of multiplicity one \gC. What the
eight have that these do not is that the whole space is one slice.
\section{The eight forms}\label{sec:eight}

Table~\ref{tab:canon} collects, for all twelve rigid lattices of Table~\ref{tab:twelve}, the data that \S\S\ref{sec:setting}--\ref{sec:canon} attach to a
rigid discriminant form: the multiplicity space $\pi_{\sigma_0}$ and its least exponent $\mu$, the single eigensystem, and the canonical scalar
$t_\mu$, together with the dimensions of the spaces of weight $\frac52$. The four theta lattices, above the rule, are the case in which these data are known in
closed form; we treat the eight below it first (\S\S\ref{subsec:table}--\ref{subsec:coeffs}) and return to the four in \S\ref{subsec:theta4}. What the
later sections use are $\mu$, $t_\mu$ and the eigenvalues of $T(p^2)$; the two columns of weight $\frac52$ are not used later, and \S\ref{subsec:wt52}
explains how they are computed.

\begin{table}[ht]
\centering\footnotesize\setlength{\tabcolsep}{4pt}
\caption{The canonical data of the twelve rigid lattices of Table~\ref{tab:twelve}: the dimension of $\pi_{\sigma_0}$ and its least $T$-exponent $\mu$
(\S\ref{subsec:mult}); the eigenvalues of $T(p^2)$ on $S(D)$ (Proposition~\ref{prop:hecke}), which are $p+1$ for the theta lattices
(Remark~\ref{rem:theta4}) and, for the eight, the coefficients $a_p(E)$ of the curve $E$ of [Cr97] named, at the primes computed \gC\
(\S\ref{subsec:hecke8}); the canonical scalar $t_\mu$ (\S\ref{subsec:onenumber}), in closed form for the theta lattices (Remark~\ref{rem:theta4})
and numerically for the eight (\S\ref{subsec:table}); and the dimensions of $S_{5/2}(\rho_L)$ and of its $O(D)$-invariant part \gC\ (\S\ref{subsec:wt52}).
Above the rule the four theta lattices, below it the eight.}
\label{tab:canon}
\begin{tabular}{lcclccc}
\toprule
$L$ & $\dim\pi_{\sigma_0}$ & $\mu$ & $T(p^2)$ & $t_\mu$ & $\dim S_{5/2}$ & $\dim S_{5/2}^{O(D)}$\\
\midrule
$L_{16}$ & $1$ & $1/8$ & $p+1$ & $\sqrt2/12$ & $0$ & $0$\\
$L_{54}$ & $2$ & $1/12$ & $p+1$ & $\sqrt3/6$ & $2$ & $0$\\
$L_{64}$ & $1$ & $1/8$ & $p+1$ & $\sqrt2/12$ & $3$ & $0$\\
$L_{250}$ & $4$ & $1/20$ & $p+1$ & $\sqrt5/3$ & $15$ & $1$\\
\midrule
$B_{40}$ & $4$ & $1/10$ & $20a1$ & $0.2697776010367185$ & $3$ & $1$\\
$L_{72}$ & $4$ & $1/12$ & $36a1$ & $0.3865317556546428$ & $6$ & $1$\\
$B_{88}$ & $10$ & $1/22$ & $44a1$ & $1.074741854937327$ & $6$ & $2$\\
$L_{96}$ & $3$ & $5/48$ & $24a1$ & $0.2323507202689307$ & $6$ & $0$\\
$B_{120}$ & $8$ & $1/20$ & $20a1$ & $0.9345370234792852$ & $7$ & $2$\\
$B_{184}$ & $22$ & $7/92$ & $92a1$ & $0.7680627886051352$ & $10$ & $4$\\
$B_{324}$ & $36$ & $1/27$ & $162d1$ & $1.713556331100267$ & $13$ & $2$\\
$B_{648}$ & $36$ & $23/432$ & $324a1$ & $1.064925231128516$ & $24$ & $3$\\
\bottomrule
\end{tabular}
\end{table}

\subsection{The eight forms and their scalars}\label{subsec:table}
For each of the eight lattices, Theorem~\ref{thm:harm} singles out one form: $\Hh_L=\Hh_{F_0}$ is the harmonic weak Maass form of weight $\frac12$ for
$\rho_L^*$ with values in $V_{\sigma_0}$ and
\[
\text{principal part }\ q^{-\mu}\,\overline{u_{r_0,1}},\qquad \xi_{1/2}\Hh_L=\frac{F_0}{t_\mu},
\]
where $u_{r_0,1}=c_{F_0}(\mu)$ (\S\ref{subsec:mult}, \S\ref{subsec:harm}). So $\Hh_L$ is fixed by three data: the exponent $\mu$ (Table~\ref{tab:canon}); the form
$F_0$, up to a scalar that does not change the coordinates of $\Hh_L$ (\S\ref{subsec:coeffs}), which is explicit for $B_{40},L_{72},L_{96}$
(Proposition~\ref{prop:etatheta}) and given by its coefficients for the other five ([Wi], Appendix~\ref{app:tabH}); and the number $t_\mu$ (Table~\ref{tab:canon}).

\noindent The number $t_\mu=\int\|F_0\|^2_\pi\,y^{3/2}\frac{dx\,dy}{y^2}\big/\|c_{F_0}(\mu)\|^2_\pi$ is computed by quadrature over the fundamental domain, with the values
of $F_0$ taken from its exact Fourier expansion [Wi] or from Hejhal's method on $\pi_{\sigma_0}$. The $16$-digit values of Table~\ref{tab:canon} come
from the exact expansion, on which two implementations agree to $6\cdot10^{-15}$ \gC; the values through Hejhal's method agree with them to $6\cdot10^{-12}$, relative, or better \gC, and for
$L_{72},B_{40},B_{88}$ they agree to $4\cdot10^{-16}$, relative, with one third of the $38$-digit values of the scalar $t$ of [II, \S5.5]
(Remark~\ref{rem:I-II}).

\begin{remark}[the numbers of {[I]} and {[II]}]\label{rem:I-II}
The scalar $t$ of $P=t\cdot B_S$ in [I, Thm.~12.1] and [II, \S5.5], for $L_{72},B_{40},B_{88}$, is defined with the $O(D)$-invariant form $B_S$ whose Gram
matrix in the echelon basis of [Wi] is $|O(D)|^{-1}\sum_{g\in O(D)}R_g^{\mathsf T}R_g$, $R_g$ the matrix of $g$ on the coordinate columns. By
Proposition~\ref{prop:Bn}, $t_\mu/t=B_S/B_\mu$, and in those bases $B_\mu=3B_S$ exactly, so $t_\mu=\frac13t$ \gC; numerically, with the values of Table~\ref{tab:canon} and the $38$-digit values of $t$ in [II, \S5.5],
$t_\mu/t=\frac13$ to $4\cdot10^{-16}$.
\end{remark}

\subsection{The lattice \texorpdfstring{$B_{40}$}{B40}}\label{subsec:b40}
This subsection shows the correspondence of Theorem~\ref{thm:harm} at work on one lattice: first the transformation law and the Fourier
coefficients of $F_0$, then $t_\mu$, then the coefficients of the form $\Hh$ attached to $F_0$. Here every class $r$ carries one index, and we write $u_\nu=u_{r,1}$ and $f_\nu=f_{r,1}$ with $\nu=e_{r,1}\in\{\frac1{10},\frac3{20},\frac7{20},\frac9{10}\}$
(\S\ref{subsec:mult}). For a suitable $v_0$ the four vectors $u_\nu$, realized in $\C[D]$, lie in $\Z[D]$ and have four entries $\pm1$ each, and
$F_0=\sum_\nu f_\nu u_\nu$ with
\begin{align*}
f_{1/10}&=q^{1/10}(1-2q-2q^2+4q^3+2q^4-4q^6+q^8-4q^9+\cdots),\\
f_{3/20}&=q^{3/20}(1-3q+q^2+3q^3-q^4+q^5-4q^6+\cdots),\\
f_{7/20}&=q^{7/20}(1-2q-q^2+q^3+2q^4+5q^5-3q^6+\cdots),\\
f_{9/10}&=q^{9/10}(1-2q-2q^2+4q^3+q^4+2q^5-2q^6+\cdots),
\end{align*}
with integer coefficients, and $f(-1/\tau)=\tau^{3/2}\mathbf S\,f(\tau)$ with
\[
\mathbf S=\e(\tfrac18)\sqrt{2/5}\begin{pmatrix}-a&-b&-a&-b\\-b&-a&b&a\\-a&b&a&-b\\-b&a&-b&a\end{pmatrix},\qquad a=\sin\tfrac{2\pi}5,\ b=\sin\tfrac\pi5
\]
\gC. Since the $u_\nu$ are orthogonal, lying in different eigenspaces of $\rho(T)$, and $\|u_\nu\|^2=4$ for all $\nu$,
$t_{1/10}=\int_{\SL_2(\Z)\backslash\mathbb H}\sum_\nu|f_\nu|^2y^{3/2}\,\frac{dx\,dy}{y^2}=0.26977760103672$.
So the form of \S\ref{subsec:table} is the harmonic weak Maass form $\Hh=\sum_\nu\Hh_\nu u_\nu$ with principal part $q^{-1/10}u_{1/10}$ and
\[
\xi_{1/2}\Hh=\frac{F_0}{t_{1/10}};
\]
its holomorphic part begins
\begin{align*}
\Hh^+_{1/10}&=q^{-1/10}-19.6363895744582236\,q^{9/10}-73.3165389886210290\,q^{19/10}-\cdots,\\
\Hh^+_{3/20}&=-8.12594932526362020\,q^{17/20}-46.1148160192841775\,q^{37/20}-\cdots,\\
\Hh^+_{7/20}&=-12.4487254561182844\,q^{13/20}-56.4642995062279837\,q^{33/20}-\cdots,\\
\Hh^+_{9/10}&=-2.36361042554177635\,q^{1/10}-14.6858186128801016\,q^{11/10}-\cdots
\end{align*}
\gC: two solves at $45$ and $40$ digits agree to $5\cdot10^{-33}$, and the computed $\Hh$ satisfies the $S$- and $T$-relations at test points to
$2\cdot10^{-35}$. Proposition~\ref{prop:etatheta} gives the $f_\nu$ in closed form.

\subsection{Eta--theta products}\label{subsec:etatheta}
For $m\ge1$ and $a\in\Z/2m$ put $\theta_{m,a}=\sum_{n\equiv a\,(2m)}q^{n^2/4m}$, so that $\theta_{m,-a}=\theta_{m,a}$ and $\Theta_m=(\theta_{m,a})_a$ is the theta
function of the lattice $\langle2m\rangle$, a vector-valued modular form of weight $\frac12$ for $\rho_{\langle2m\rangle}$ [Bo98, Thm.~4.1]; let
$\eta=q^{1/24}\prod_{n\ge1}(1-q^n)$. Write $F_0=\sum f_{r,i}u_{r,i}$ with $f_{r,i}=q^{e_{r,i}}(1+O(q))$ in the basis of \S\ref{subsec:mult}; for the three
lattices below every class $r$ carries one index, and $f_e$ denotes the coordinate function of exponent $e$.

\begin{proposition}[three of the eight $F_0$ as eta--theta products]\label{prop:etatheta}
\begin{align*}
B_{40}:\ & f_{1/10}=\eta^2(\theta_{15,1}-\theta_{15,11}),\quad f_{3/20}=\eta^2(\theta_{15,2}-\theta_{15,8}),\\
& f_{7/20}=\eta^2(\theta_{15,4}-\theta_{15,14}),\quad f_{9/10}=\eta^2(\theta_{15,7}-\theta_{15,17});\\
L_{96}:\ & f_{5/48}=\eta^2(\theta_{12,1}-\theta_{12,7}),\quad f_{1/6}=\eta^2(\theta_{12,2}-\theta_{12,10}),\quad
 f_{29/48}=\eta^2(\theta_{12,5}-\theta_{12,13});\\
L_{72}:\ & f_{1/12}=\eta^2\theta_{3,0},\quad f_{1/6}=\eta^2\theta_{3,1},\quad f_{5/12}=\eta^2\theta_{3,2},\quad f_{5/6}=\tfrac12\,\eta^2\theta_{3,3}.
\end{align*}
For $B_{88},B_{120},B_{184},B_{324},B_{648}$ one has $\mu<\frac1{12}$, so no identity of this shape holds for them; moreover the Wronskian of the
coordinate functions of $F_0$ has zeros in $\mathbb H$. \gC\ (given the computed exponents)
\end{proposition}
\begin{proof}
Let $k=\frac32$, $\partial_k=q\frac d{dq}-\frac k{12}E_2$ the Serre derivative, $\partial^i=\partial_{k+2i-2}\circ\cdots\circ\partial_k$,
$V=\operatorname{span}\{f_{r,i}\}$ and $d=\dim V=\dim\pi_{\sigma_0}$ (the $u_{r,i}$ are a basis of $\pi_{\sigma_0}$). Number the $f_{r,i}$ as $f_1,\dots,f_d$, with
exponents $e_1,\dots,e_d$.

\emph{Step 1: the Wronskian.} Writing $\rho(\gamma)u_l=\sum_jR(\gamma)_{jl}u_j$ on $\pi_{\sigma_0}$, the law
$F_0(\gamma\tau)=\varphi(\tau)^{2k}\rho(\gamma)F_0(\tau)$, $\gamma=(M,\varphi)\in\Mp$, reads $f_j(M\tau)=\varphi^{2k}\sum_lR(\gamma)_{jl}f_l(\tau)$; the Serre derivative
preserves this and raises the weight by $2$. So $\mathbf M=(\partial^if_j)_{0\le i<d,\,1\le j\le d}$ satisfies
$\mathbf M(M\tau)=\operatorname{diag}_i(\varphi^{2(k+2i)})\,\mathbf M(\tau)\,R(\gamma)^{\mathsf T}$, and $W=\det\mathbf M=\det\bigl((q\tfrac d{dq})^if_j\bigr)$ is
holomorphic on $\mathbb H$, of weight $\sum_i(k+2i)=d(2d+1)/2$, with the character $\det R$ of finite order. As $\partial^if_j=q^{e_j}\bigl(p_i(e_j)+O(q)\bigr)$ with
$p_i$ monic of degree $i$ and the $e_j$ are distinct, $\operatorname{ord}_\infty W=\sum_je_j$.

\emph{Step 2: zeros.} By the valence formula [Se73, Ch.~VII, \S3, Thm.~3], applied to a power of $W$ with trivial character,
$\sum_P\operatorname{ord}_P(W)/|\bar\Gamma_P|+\sum_je_j=d(2d+1)/24$, the sum over $P\in\SL_2(\Z)\backslash\mathbb H$ and $\bar\Gamma_P$ the stabilizer of $P$ in
$\mathrm{PSL}_2(\Z)$. So $W$ has no zero in $\mathbb H$ iff $\sum_je_j=d(2d+1)/24$. With $(d,\sum_je_j,d(2d+1)/24)$ from Table~\ref{tab:H}: $B_{40},L_{72}$
$(4,\frac32,\frac32)$, $L_{96}$ $(3,\frac78,\frac78)$: no zeros; $B_{88}$ $(10,\frac{19}4,\frac{35}4)$, $B_{120}$ $(8,\frac{11}3,\frac{17}3)$, $B_{184}$
$(22,\frac{85}4,\frac{165}4)$, $B_{324}$ and $B_{648}$ $(36,\frac{63}2,\frac{219}2)$: zeros.

\emph{Step 3: the equation.} Let $W$ be zero-free. Expanding the $(d+1)\times(d+1)$ determinant with columns $(\partial^if)_i,(\partial^if_1)_i,\dots,(\partial^if_d)_i$
along the first column, every $f\in V$ satisfies $\partial^df+P_1\partial^{d-1}f+\dots+P_df=0$ with $P_i=\pm W_{d-i}/W$, $W_{d-i}$ the minor without row $d-i$,
of weight $d(2d+1)/2+2i$, character $\det R$ and $\operatorname{ord}_\infty\ge\sum_je_j$. So $P_i\in M_{2i}(\SL_2(\Z))$: $P_1=0$,
$(P_2,P_3,P_4)=(c_2E_4,c_3E_6,c_4E_4^2)$ for $d=4$, $(P_2,P_3)=(c_2E_4,c_3E_6)$ for $d=3$, with constants $c_i$. Inserting $f=f_j$ and comparing leading terms,
with $x_j=e_j-\frac18$, the monic polynomial $\prod_j(x-x_j)$ equals $\prod_{l<d}(x-\frac l6)+c_2\,x(x-\frac16)+c_3x+c_4$ ($d=4$), resp.\
$\prod_{l<3}(x-\frac l6)+c_2x+c_3$ ($d=3$); its $x^{d-1}$-coefficient is the zero-free condition again, and the others determine the $c_i$. So $V$ is the full
solution space of a monic equation with holomorphic coefficients determined by the multiset $\{e_j\}$.

\emph{Step 4: the eta--theta spaces.} Let $A_m=(\Z/2m,x^2/4m)$, the discriminant form of $\langle2m\rangle$. For $m=pp'$, $(p,p')=(3,5),(3,4)$, let
$\beta(a)\equiv a\ (2p)$, $\beta(a)\equiv-a\ (2p')$; it is well defined, as $\operatorname{lcm}(2p,2p')=2m$, and $\beta\in O(A_m)$, as $2p\mid\beta(a)-a$ and
$2p'\mid\beta(a)+a$. Put $U_m=\{v\in\C[A_m]:v(-a)=v(a),\ v(\beta a)=-v(a)\}$, which is $\rho_{A_m}$-stable. The vectors
$\tilde u_s=e_{p'-ps}+e_{-(p'-ps)}-e_{p'+ps}-e_{-(p'+ps)}$, $1\le s<p'$, lie in $U_m$, have disjoint supports and span it: the even part of $\C[A_m]$ has
dimension $16$, resp.\ $13$, and $\beta$ acts on it with trace $8$, resp.\ $7$ (an orbit $\{a,-a\}$ is $\beta$-stable iff $p\mid a$ or $p'\mid a$, and then
$\beta$ fixes its even vector), so $\dim U_m=(p-1)(p'-1)/2$. For $m=3$ let $U_3$ be the even part of $\C[\Z/6]$. Then $\eta^2\langle\Theta_m,\cdot\rangle|_{U_m}$ is a
vector-valued form of weight $\frac32$ whose coordinate functions are, up to constant factors, the right sides of the Proposition, with exponents
$\frac1{12}+\frac{(p'-ps)^2}{4m}$, resp.\ $\frac1{12}+\frac{j^2}{12}$:
$\{\frac3{20},\frac1{10},\frac7{20},\frac9{10}\}$, $\{\frac5{48},\frac16,\frac{29}{48}\}$, $\{\frac1{12},\frac16,\frac5{12},\frac56\}$ --- the multisets $\{e_{r,i}\}$
of $B_{40}$, $L_{96}$, $L_{72}$. Their Wronskian is zero-free by Step 2, so by Step 3 they solve the same equation as $V$, and $V$ is their span.

\emph{Step 5.} The exponents are distinct modulo $1$, so each $f_e$ is proportional to the eta--theta function with the same exponent; the leading
coefficients are $1$ on both sides, except $\theta_{3,3}=2q^{3/4}+\cdots$. For the five other lattices, $\eta^2$ times a theta function of a lattice of rank one
has exponents $\ge\frac1{12}$, while $\mu<\frac1{12}$ (Table~\ref{tab:canon}).
\end{proof}

\subsection{The coefficients of \texorpdfstring{$\mathcal H$}{H}}\label{subsec:coeffs}
In the basis $u_{r,i}$ of \S\ref{subsec:mult} (for $B_{40}$ the basis $u_\nu$ of \S\ref{subsec:b40}) write $\Hh=\sum h_{r,i}\,\overline{u_{r,i}}$, as in
\S\ref{subsec:harm}. Replacing $F_0$ by $\lambda F_0$ replaces $u_{r,i}$ by $\lambda u_{r,i}$ and $\Hh$ by $\bar\lambda\Hh$, so the $h_{r,i}$ depend neither on the
normalization of $F_0$ nor, by \S\ref{subsec:harm}, on $v_0$. With $(r_0,1)$ the index of $\mu$,
\begin{align*}
h_{r,i}=\delta_{(r,i),(r_0,1)}\,q^{-\mu}&+\sum_{n>0,\ n\equiv-r}a_{r,i}(n)\,q^n\\
&-\frac1{t_\mu}\sum_{m>0}\overline{c_{f_{r,i}}(m)}\,(4\pi m)^{-1/2}\,\Gamma(\tfrac12,4\pi my)\,q^{-m},
\end{align*}
the last sum being the non-holomorphic Eichler integral of $f_{r,i}/t_\mu$, as $\xi_{1/2}\bigl(\Gamma(\frac12,4\pi my)q^{-m}\bigr)=-(4\pi m)^{1/2}q^m$.
Table~\ref{tab:H} gives $t_\mu$ and the first four $a_{r,i}(n)$. They are computed by Hejhal's method for weight $\frac12$ and $\rho_D^*$ on $\pi_{\sigma_0}$
[BSt12], in double precision, with three choices of truncation and sample height ($M,Q,Y=30,50,0.25$; $40,60,0.20$; $36,56,0.22$), from
$F_0$ and $t_\mu$ and along two further routes. Relative to $\max(1,|a|)$, the three choices agree to $6\cdot10^{-14}$ and the routes to
$1.3\cdot10^{-11}$ on the tabulated coefficients; the computed expansions satisfy
$\Hh(-1/\tau)=\sqrt\tau\,\rho_D^*(S)\Hh(\tau)$ at four points off the sample points to $3.1\cdot10^{-14}$; the rows of $B_{40}$ agree
with \S\ref{subsec:b40} to $3\cdot10^{-15}$ \gC.

\subsection{The four theta lattices}\label{subsec:theta4}
We return to the four theta lattices above the rule of Table~\ref{tab:canon}.

\begin{remark}[the four theta lattices]\label{rem:theta4}
For $L_{16},L_{54},L_{64},L_{250}$, $s_L$, the square-free part of $|D|/2$, is $2,3,2,5$. By Theorem~A (\S\ref{subsec:classno}) $S(D)$ is spanned by the lifts, through the isotropic subgroups $U$ of $D$ with $U^\perp/U\cong(\Z/2s_L,x^2/4s_L)$, of the
vector-valued unary theta series $\vartheta_{s_L}=(\vartheta_{s_L,r})_{r\in\Z/2s_L}$ of weight $\frac32$, $\vartheta_{s_L,r}=\sum_{\ell\equiv r\,(2s_L)}\ell\,q^{\ell^2/4s_L}$ (the lift through $U$ is
$e_{\gamma+U}\mapsto\sum_{u\in U}e_{\gamma+u}$, which maps $S_{3/2}$ for $U^\perp/U$ into $S(D)$ [BEF16, Prop.~4.18], \S\ref{subsec:classno});
$\pi_{\sigma_0}\cong V^-(\Z/2s_L)$, the odd part of the Weil representation of $(\Z/2s_L,x^2/4s_L)$, and $\mu=1/4s_L$ \gC. Given these, the value of
$t_\mu$ in Table~\ref{tab:canon} follows in two lines from [BSc17]. For the Petersson product of
\S\ref{sec:canon} on $(\Z/2s_L,x^2/4s_L)$, $\|\vartheta_{s_L}\|^2=\sqrt{s_L}(s_L-1)/6$ [BSc17, Example~5.2]. The coefficient of $\vartheta_{s_L}$ at $n=1/4s_L$
is $e_1-e_{-1}$, of norm $2$. Each coordinate of $\vartheta_{s_L}$ occurs in exactly $|U|$ coordinates of its lift, so the lift through $U$ multiplies
$P$ and $B_n$ by the same factor $|U|$. As $P=t_nB_n$ on $S(D)$ (Proposition~\ref{prop:Bn}),
\[
t_{1/4s_L}=\frac{\sqrt{s_L}(s_L-1)/6}{2}=\frac{(s_L-1)\sqrt{s_L}}{12}.
\]
[BSc17] notes that the Rankin--Selberg method gives the same norm.
At every prime $p\nmid2|D|$, $T(p^2)$ acts on $S(D)$ by $p+1$. Indeed, let $A=(\Z/2s_L,x^2/4s_L)$, of signature $1$. The coefficients of
$\vartheta_{s_L}$ vanish except at $n=\ell^2/4s_L$, $\ell\in\Z$, where $2|A|n=\ell^2$. As $p$ is a unit modulo $2s_L$, which divides $|D|$, the three terms of \eqref{eq:T32} for
$A$ (whose proof in Proposition~\ref{prop:hecke} uses only $\sig A\equiv1\pmod8$) at $(r,n)$ are $p$, $(\frac{\ell^2}p)$ and ${[p\mid\ell]}$ times the
coefficient of $\vartheta_{s_L}$ at $(r,n)$, and $(\frac{\ell^2}p)+{[p\mid\ell]}=1$; at the other $n$ all three terms vanish. So $T(p^2)\vartheta_{s_L}=(p+1)\vartheta_{s_L}$. The lift through $U$
commutes with $T(p^2)$: multiplication by $p$ is an automorphism of $D$ that preserves $U$ and $U^\perp$, and $(\frac{2|D|n}p)=(\frac{2|A|n}p)$
because $|D|=|U|^2|A|$ and $p\nmid|U|$. As these lifts span $S(D)$, the claim follows.
The form $\Hh$ is $-\frac12$ times the lift of the completion $\widehat H^{(s_L)}$ of the mock modular form $H^{(s_L)}$ of umbral moonshine of lambency
$s_L$ [CDH14a], [DGO15], with $\xi_{1/2}\widehat H^{(s_L)}=-(2/t_{1/4s_L})\,\vartheta_{s_L}$. By [ORTL15, Thm.~4.1], for
$m\in\{2,3,4,5,7,9,13,25\}$ the completion $\widehat H^{(m)}$ is a harmonic weak Maass form of weight $\frac12$ for $\rho_A^*$, $A=(\Z/2m,x^2/4m)$
([ORTL15] writes $\overline\rho_L$), with principal part $-2q^{-1/4m}(e_1-e_{-1})$. As $s_L\in\{2,3,5\}$ and the coefficient of $\vartheta_{s_L}$ at $\mu$
is $e_1-e_{-1}$, $-\frac12$ times its lift has the principal part required in Theorem~\ref{thm:harm}. It also takes values in $V_{\sigma_0}$, as
Theorem~\ref{thm:harm} requires: the proof of Lemma~\ref{lem:parity} uses only the law at $Z$, so $\widehat H^{(s_L)}$ takes values in
$V^-(\Z/2s_L)$; and the lift maps $V^-(\Z/2s_L)$ into $V_{\sigma_0}$, since it maps $\vartheta_{s_L}$ into $S(D)$, whose values lie in $V_{\sigma_0}$, and the
values of $\vartheta_{s_L}$ span $V^-(\Z/2s_L)$, its components $\vartheta_{s_L,r}$, $0<r<s_L$, having the distinct leading exponents $r^2/4s_L$.
So Theorem~\ref{thm:harm} gives the constant $-(2/t_{1/4s_L})$ of the shadow. For $s_L=2$ ($L_{16}$ and $L_{64}$) this is Mathieu moonshine:
$\xi_{1/2}\widehat H^{(2)}_1=-12\sqrt2\,\eta^3$, so the constant $24$ in the completion of $H^{(2)}_1=2q^{-1/8}(-1+45q+231q^2+770q^3+2277q^4+\cdots)$ and the
value $t_{1/8}=\sqrt2/12$ are the same information, and solving for $\Hh$ with this value gives back the coefficients $45,231,770,2277,\dots$ \gC.
The lifts of $\vartheta_{s_L}$ are not the eta--theta products of \S\ref{subsec:etatheta}: there the theta functions $\theta_{m,a}$
have weight $\frac12$ and are multiplied by $\eta^2$, and by Theorem~A none of the eight lattices carries a lift of a unary theta series.
\end{remark}
\section{Hecke eigenvalues of the eight}\label{sec:arith}\label{subsec:hecke8}
For each of the eight lattices the eigenvalues $\lambda_p$ of Proposition~\ref{prop:hecke} have been computed at the primes $p\le31$, $p\nmid|D|$, that the
computed coefficients reach (six to nine primes for each lattice), from the relations \eqref{eq:T32}; they are the coefficients $a_p(E)$ of the curves
\[
\begin{array}{llll}
B_{40}:20a1,\quad & L_{72}:36a1,\quad & B_{88}:44a1,\quad & L_{96}:24a1,\\
B_{120}:20a1,\quad & B_{184}:92a1,\quad & B_{324}:162d1,\quad & B_{648}:324a1
\end{array}
\]
\gC. In the extended tables of [Cr97], each of these is the only isogeny class of conductor $<10^4$ whose $a_p$ are the
$\lambda_p$ at the primes computed, also if bad reduction at these primes is allowed \gC. The conductor is $|\det L|/2$ except for $L_{96}$ ($24$, not
$48$) and $B_{120}$ ($20$, not $60$). That the whole eigensystem is that of $E$ is Problem~\ref{prob:hecke}.

\begin{proposition}[complex multiplication]\label{prop:cm}
\begin{enumerate}[label=(\alph*),leftmargin=7mm,itemsep=1pt]
\item The eigenvalues $\lambda_p$ of a rigid $S(D)$, $D=D(L)$, are integers.
\item Let $L$ be one of the eight, and suppose that $\lambda_p=a_p(f)$ for every prime $p\nmid2|D|$, for a newform $f$ of weight $2$. Then $f=f_E$ for an
elliptic curve $E$ over $\Q$, and if $E$ has complex multiplication, then $L=L_{72}$. \gC
\end{enumerate}
\end{proposition}
\begin{proof}
\begin{cited}{McGraw [McG03, Thm.~5.6]}
For an even lattice $L$ and $k\in\frac12\Z$, the spaces $M_k(\rho_L)$ and $S_k(\rho_L)$ have bases of forms whose Fourier coefficients are integers.
\end{cited}
\why{the rational structure of the spaces of vector-valued modular forms, used as stated.}
(a) By [McG03] there is $F\ne0$ in $S(D)$ with integral coefficients; dividing it by the greatest common divisor of its coefficients, we may
assume that they have no common factor. The coefficients of \eqref{eq:T32} are integers, so $T(p^2)F=\lambda_pF$ has integral coefficients; writing $1$
as an integral combination of finitely many coefficients of $F$ gives $\lambda_p\in\Z$.

(b) The $a_p(f)=\lambda_p$, $p\nmid2|D|$, are rational by (a). Every Galois conjugate of $f$ has the same $a_p$, $p\nmid2|D|$, hence equals $f$ by
strong multiplicity one (\S\ref{subsec:newforms}(1)); so $f$ has rational coefficients, and $f=f_E$ (\S\ref{subsec:newforms}(2)). If $E$ has complex multiplication by an order in an imaginary quadratic field $K$,
then $j(E)\in\Q$, so $K$ has class number one, $K=\Q(\sqrt{-d})$ with $d\in\{1,2,3,7,11,19,43,67,163\}$. Then $a_p(E)=0$ at every prime $p$ inert in $K$:
at the primes of good reduction by Deuring's criterion, and at the primes of bad reduction because $j(E)$ is integral, so that $E$ has potentially good,
hence additive, reduction there. For each of the seven lattices other than $L_{72}$ and each of the nine fields, the table below gives a prime
$p\nmid2|D|$ inert in $K$ with $\lambda_p\ne0$ \gC:
\begin{center}\small
\begin{tabular}{llrrrrrrrrr}
\toprule
$L$ & $E$ & $-3$ & $-4$ & $-7$ & $-8$ & $-11$ & $-19$ & $-43$ & $-67$ & $-163$\\
\midrule
$B_{40}$ & $20a1$ & 17 & 3 & 3 & 7 & 7 & 3 & 3 & 3 & 3\\
$B_{88}$ & $44a1$ & 5 & 3 & 3 & 5 & 7 & 3 & 3 & 3 & 3\\
$L_{96}$ & $24a1$ & 5 & 11 & 5 & 5 & 13 & 13 & 5 & 5 & 5\\
$B_{120}$ & $20a1$ & 17 & 7 & 13 & 7 & 7 & 13 & 7 & 7 & 7\\
$B_{184}$ & $92a1$ & 17 & 3 & 3 & 7 & 7 & 3 & 3 & 3 & 3\\
$B_{324}$ & $162d1$ & 5 & 7 & 5 & 5 & 7 & 13 & 5 & 5 & 5\\
$B_{648}$ & $324a1$ & 5 & 7 & 5 & 5 & 7 & 13 & 5 & 5 & 5\\
$L_{72}$ & $36a1$ & -- & 7 & 13 & 7 & 7 & 13 & 7 & 7 & 7\\
\bottomrule
\end{tabular}
\end{center}
For $L_{72}$ the column $-3$ is empty: $36a1$ has complex multiplication by $\Z[\zeta_3]$, and $\lambda_p=0$ at the primes $p\equiv2\pmod3$ computed
\gC.
\end{proof}
\noindent Without rigidity the same form is no longer canonical: the isotropic lift through a subgroup of order $2$ carries $F_0$ and $\Hh$ of
$L_{72}$ to the sublattice $\operatorname{diag}(-8,6,6)$ of index $2$, but there they lie in a $6$-dimensional
$36a1$-part of a $14$-dimensional $S(D)$ with $\langle\chi_S,\chi_S\rangle=5$ \gC.

\section{Rigidity by counting}\label{sec:count}

Sections~\ref{sec:canon}--\ref{sec:arith} took rigidity as given and asked what it fixes. This section treats one family of discriminant forms,
$\Dp_p$ below, through the newforms of weight $2$ and level $p$, $2p$ or $4p$. Their number with Atkin--Lehner sign $w_p=+1$ is a difference of
genera of two Atkin--Lehner quotients of modular curves, and that difference has a closed formula in terms of $p$ and one class number
(Lemma~\ref{lem:three}). The rigid members of the family are the primes at which the count is $1$, equivalently at which $X_0(2p)/w_p$
has genus $0$ and $X_0(4p)/w_p$ genus $1$ (Corollary~\ref{cor:geom}, by comparison with the list of Theorem~A); Proposition~\ref{prop:allw}
then extends the count to every even weight. In this section $p$ is the prime of $\Dp_p$.

\subsection{The family}\label{subsec:fam}
Let $p$ be an odd prime and
\[
\Dp_p=(\Z/2)^3_{\mathrm{an}}\oplus(\Z/p,\,ux^2/p).
\]
The class of $u$ and the $2$-part are not free parameters: by Milgram's formula, $\sig\Dp_p\equiv1\pmod8$ forces $u$ to be a non-square modulo $p$ and the
$2$-part to be $\langle-2\rangle^3$ for $p\equiv1\pmod4$ and $\langle2\rangle^3$ for $p\equiv3\pmod4$ (the Gauss sums give $\sig\langle\pm2\rangle^3=\pm3$ and
$\sig(\Z/p,ux^2/p)\in\{0,4\}$ for $p\equiv1$, $\{2,6\}$ for $p\equiv3\pmod4$, the second value for $u$ a non-square); so $\Dp_p$ is one module for each $p$.
One has $\Dp_5=D(B_{40})$, $\Dp_{11}=D(B_{88})$, $\Dp_{23}=D(B_{184})$.

\begin{lemma}\label{lem:two}
Let $L$ be an even lattice of rank $3$ with $D(L)_2\cong(\Z/2)^3$ anisotropic. Then $L\otimes\Q_2$ is anisotropic, so $2\mid\fd(L)$
(\S\ref{subsec:lat}).
\end{lemma}
\begin{proof}
Write $L_2=L\otimes\Z_2$. As $D(L)_2=L_2'/L_2$ has exponent $2$, $2L_2'\subseteq L_2$, and $[L_2':L_2]=8=[L_2':2L_2']$; so $L_2=2L_2'$. Put $M=L_2'$ with
$(x,y)_M=2(x,y)_L$, which is integral on $M$ since $2x\in L_2$. Then $D(L)_2=M/2M$, and for $x\in M$, $q(x)=(x,x)_L/2=(x,x)_M/4\bmod1$. If $L\otimes\Q_2$ were
isotropic it would contain an isotropic vector $v\in M\setminus2M$; its class in $M/2M$ is non-zero with $q=(v,v)_M/4=0$, so $D(L)_2$ would be isotropic.
\end{proof}

\begin{corollary}\label{cor:fd}
If $D(L)\cong\Dp_p$, then $\fd(L)=2p$. If $D(L)\cong(\Z/2)^3_{\rm an}\oplus(\Z/\ell,u_\ell x^2/\ell)\oplus(\Z/\ell',u_{\ell'}x^2/\ell')$ for distinct odd primes
$\ell,\ell'$, then exactly one of $\ell,\ell'$ divides $\fd(L)$.
\end{corollary}
\begin{proof}
The ramified places of the quaternion algebra of $L\otimes\Q$ are even in number (Hilbert reciprocity), exclude $\infty$ ($L$ is indefinite) and the primes not dividing $2|D(L)|$
(\S\ref{subsec:lat}), and include $2$ by Lemma~\ref{lem:two}.
\end{proof}
\noindent Directly, $L\otimes\Z_p\cong\langle a_1,a_2,pa_3\rangle$ with units $a_i$ and $u\equiv a_3/2$, and $L\otimes\Q_p$ is anisotropic exactly when $u$ is a
non-square, which Milgram's formula forces; Hilbert symbols confirm $\fd=2p$ for $p\le97$ \gC.

\noindent\textbf{Lattices.} If $D(L)\cong\Dp_p$, the proof of Lemma~\ref{lem:two} shows that all inner products of $L$ are even: $L=G'(2)$ for an integral
symmetric $G'$ of signature $(2,1)$ and determinant $-p$ with $v^{\mathsf T}G'v\not\equiv0\pmod4$ for $v\in\{0,1\}^3\setminus\{0\}$; conversely every such
$G'$ gives $\Dp_p$. Examples: $G'=\operatorname{diag}(1,1,-p)$ for $p\equiv3\pmod4$; for $p\equiv1\pmod4$,
\[
G'=\langle-1\rangle\oplus\begin{pmatrix}p'&b\\b&(p+b^2)/p'\end{pmatrix},\qquad p'\equiv3\pmod4\ \text{prime},\ \Bigl(\frac{p'}p\Bigr)=-1,\ b^2\equiv-p\pmod{p'},
\]
where $p'$ exists by Dirichlet's theorem and $b$ because $(\frac{-p}{p'})=-(\frac{p'}p)=1$; the binary form is positive definite of determinant $p$ and
represents $p'$, so that
$G'\otimes\Z_2\cong\langle-1,p',pp'\rangle$, all three entries $\equiv3\pmod4$, and the condition mod $4$ holds (for $p=17$:
$\langle-1\rangle\oplus\bigl(\begin{smallmatrix}3&1\\1&6\end{smallmatrix}\bigr)$; for $p\equiv5\pmod8$ also
$\bigl(\begin{smallmatrix}1&0&0\\0&2&1\\0&1&(1-p)/2\end{smallmatrix}\bigr)$). For the $211$ primes $p\equiv1\pmod4$ below $3000$ the determinant, the
signature, the condition mod $4$, $D(L)_2\cong\langle-2\rangle^3$ and $u$ a non-square are also checked directly \gC.

\noindent\textbf{Atkin--Lehner signs.} Let $f$ be a newform of weight $w$ and level $N$ (\S\ref{subsec:newforms}). For $\ell\,\|\,N$ the newform $f$
is an eigenform of the Atkin--Lehner involution $W_\ell$, normalized so that $W_\ell^2=1$; its eigenvalue $w_\ell(f)=\pm1$ is the
\emph{Atkin--Lehner sign} at $\ell$, and $a_\ell(f)=-w_\ell(f)\,\ell^{w/2-1}$ [AL70, Thm.~3] (stated there in weight $2k$: $a(q)=-q^{k-1}\lambda(q)$).

\noindent\emph{Notation.} For a level $M$, a weight $w$ and $\varepsilon=\pm1$, $n^{\varepsilon}_M(w)$ is the number of newforms of weight $w$ and level $M$
(\S\ref{subsec:newforms}) with $W_p$-eigenvalue $w_p(f)=\varepsilon$, $p\,\|\,M$; for $w=2$ this is $a_p(f)=-\varepsilon$; $n^{\varepsilon}_M=n^{\varepsilon}_M(2)$.
$h(d)$ is the number of classes of primitive positive definite binary quadratic forms of discriminant $d$; so $h(-4p)$ is the class number of $\Z[\sqrt{-p}]$,
and for $p\equiv3\pmod4$, $h(-p)$ that of the maximal order. $X_0(N)/w_p$ is the quotient of $X_0(N)$ by the Atkin--Lehner involution and $g(\cdot)$ is the
genus. The letter $w$ is a weight of integral-weight forms; Atkin--Lehner signs are written $\varepsilon$ or $w_\ell(f)$.

\subsection{Fixed points of \texorpdfstring{$w_p$}{wp}}

\begin{lemma}\label{lem:fix}
Let $p>3$ be prime and $M\in\{1,2,4\}$. The number $F_{Mp}$ of fixed points of $w_p$ on $X_0(Mp)$ is
\[
F_{Mp}=\sum_{\mathcal O\ni\sqrt{-p}}h(\mathcal O)\,\rho_M(\mathcal O),
\]
the sum running over the orders of $\Q(\sqrt{-p})$ that contain $\alpha=\sqrt{-p}$, where $\rho_M(\mathcal O)$ is the number of cyclic subgroups of order $M$ of
$\mathcal O/M\mathcal O$ stable under multiplication by $\alpha$. The values are
\begin{center}
\begin{tabular}{llccc}
\toprule
$p$ & $\mathcal O$ & $\rho_1$ & $\rho_2$ & $\rho_4$\\
\midrule
$\equiv1\ (4)$ & $\Z[\sqrt{-p}]$ & $1$ & $1$ & $0$\\
$\equiv3\ (4)$ & $\Z[\frac{1+\sqrt{-p}}2]$ & $1$ & $3$ & $2\bigl(1+(\frac{-p}{2})\bigr)$\\
$\equiv3\ (4)$ & $\Z[\sqrt{-p}]$ & $1$ & $1$ & $2$\\
\bottomrule
\end{tabular}
\end{center}
so that $F_p=h(-4p)+[p\equiv3\ (4)]\,h(-p)$, and
\[
F_{2p}=\begin{cases}h(-4p)&p\equiv1\ (4),\\ 3h(-p)+h(-4p)&p\equiv3\ (4),\end{cases}\qquad
F_{4p}=\begin{cases}0&p\equiv1\ (4),\\ 2\bigl(1+(\tfrac{-p}{2})\bigr)h(-p)+2h(-4p)&p\equiv3\ (4).\end{cases}
\]
For $p\equiv3\pmod4$ this reads $F_{2p}=4h(-p)$, $F_{4p}=6h(-p)$ if $p\equiv7\pmod8$ and $F_{2p}=F_{4p}=6h(-p)$ if $p\equiv3\pmod8$.
\end{lemma}

\begin{proof}
\emph{Fixed points as endomorphisms.} A point of $Y_0(Mp)$ is a class of triples $(E,C_M,C_p)$ with $C_M$, $C_p$ cyclic of orders $M$, $p$, and
$w_p(E,C_M,C_p)=(E/C_p,(C_M+C_p)/C_p,E[p]/C_p)$. A fixed point is an isomorphism $\varphi\colon E/C_p\to E$ with $\varphi((C_M+C_p)/C_p)=C_M$ and
$\varphi(E[p]/C_p)=C_p$. Put $\alpha=\varphi\circ\pi$ with $\pi\colon E\to E/C_p$. Then $\deg\alpha=p$, $\ker\alpha=C_p$, $\alpha(C_M)=C_M$, and
$\alpha(E[p])=C_p=\ker\alpha$, so $\ker\alpha^2\supseteq E[p]$; as $\deg\alpha^2=p^2$, $\alpha^2=p\,u$ with $u\in\mathrm{Aut}(E)$, and $\alpha=u\hat\alpha$. If $u=\pm1$,
then $\alpha=-\hat\alpha$ (the case $\alpha=\hat\alpha$ would make $\alpha$ a rational integer of square $p$), so $\alpha^2=-p$ and $\alpha=\pm\sqrt{-p}\in\mathrm{End}(E)$. If
$\mathrm{End}(E)$ is $\Z[i]$ or $\Z[\zeta_3]$ and $u\ne\pm1$, then $\alpha$ and $\hat\alpha$ are associate elements of norm $p$, so $p$ ramifies in $\Z[i]$ or $\Z[\zeta_3]$,
i.e.\ $p\in\{2,3\}$. Hence for $p>3$ the fixed points are the triples $(E,C_M,\ker\sqrt{-p})$ with $\sqrt{-p}\in\mathrm{End}(E)$ and $C_M$ stable under
$\sqrt{-p}$; conversely each such triple is fixed. No cusp is fixed, since $w_p$ changes the $p$-part of the denominator of a cusp. The two elements
$\pm\sqrt{-p}$ have the same kernel and $\mathrm{Aut}(E)=\{\pm1\}$ fixes every subgroup, so the fixed points with $\mathrm{End}(E)=\mathcal O$ are the pairs ($E$ up
to isomorphism, $C_M$), $h(\mathcal O)$ choices for $E$.

\emph{The level structure.} If $\mathrm{End}(E)=\mathcal O$, then $E\cong\C/\mathfrak a$ with $\mathfrak a$ a proper, hence invertible, $\mathcal O$-ideal, and
$E[M]\cong\mathfrak a/M\mathfrak a\cong\mathcal O/M\mathcal O$ as $\mathcal O$-modules; so the number of admissible $C_M$ is $\rho_M(\mathcal O)$ for every such $E$.

\emph{The numbers $\rho_M$.} $\rho_1=1$. On a stable cyclic $\langle v\rangle$ of order $4$, $\alpha$ acts by some $\lambda\in(\Z/4)^\times$ ($\alpha$ has odd norm,
so it is invertible on $E[4]$), and $\alpha^2=-p$ gives $(\lambda^2+p)v=0$, i.e.\ $p\equiv3\pmod4$ and $\lambda=\pm1$.
\begin{itemize}[leftmargin=*,itemsep=1pt]
\item $\mathcal O=\Z[x]/(x^2+p)$, $\alpha=x$. Modulo $2$, $x^2+p\equiv(x+1)^2$, so $\alpha-1$ is nilpotent and non-zero on $\mathcal O/2\mathcal O\cong\F_2^2$: one
stable line, $\rho_2=1$. For $p\equiv1\pmod4$, $\rho_4=0$ by the remark above. For $p\equiv3\pmod4$, modulo $4$ one has $x^2\equiv1$, and
$x(a+bx)=b+ax=\lambda(a+bx)$ forces $b=\lambda a$; with $a$ odd this gives exactly $\langle1+x\rangle$ and $\langle1-x\rangle$, so $\rho_4=2$.
\item $\mathcal O_1=\Z[\omega]$, $\omega=\frac{1+\sqrt{-p}}2$, $p\equiv3\pmod4$, $\alpha=2\omega-1$. Modulo $2\mathcal O_1$, $\alpha\equiv1$, so all three subgroups of
order $2$ are stable: $\rho_2=3$. Modulo $4$, $\alpha-1=-2\bar\omega$ and $\alpha+1=2\omega$. If $p\equiv3\pmod8$, $N(\omega)=\frac{1+p}4$ is odd, $\omega$ and
$\bar\omega$ are units at $2$, and $\ker(\alpha\mp1)$ on $\mathcal O_1/4\mathcal O_1$ is $2\mathcal O_1/4\mathcal O_1$, of exponent $2$: $\rho_4=0$. If $p\equiv7\pmod8$,
then $2=\mathfrak p\mathfrak p'$ with $\omega\in\mathfrak p$, $\bar\omega\in\mathfrak p'$; as $2\omega\in2\mathfrak p\subset\mathfrak p^2$ and
$\omega\equiv1\pmod{\mathfrak p'}$, $\mathcal O_1/4\mathcal O_1\cong\mathcal O_1/\mathfrak p^2\times\mathcal O_1/\mathfrak p'^2\cong\Z/4\times\Z/4$ with
$\alpha=(-1,+1)$, whose stable cyclic subgroups of order $4$ are $\langle(1,0)\rangle$, $\langle(1,2)\rangle$, $\langle(0,1)\rangle$, $\langle(2,1)\rangle$: $\rho_4=4$.
\end{itemize}
Finally $h(\Z[\frac{1+\sqrt{-p}}2])=h(-p)$, $h(\Z[\sqrt{-p}])=h(-4p)$, and for $p\equiv3\pmod4$, $p>3$, $h(-4p)=\bigl(2-(\frac{-p}2)\bigr)h(-p)$, the class number
of the order of conductor $2$, the unit groups being $\{\pm1\}$ [Co13, Thm.~7.24].
\end{proof}
\noindent The values agree with $F_N=2-2\tr(W_p\,|\,S_2(\Gamma_0(N)))$ (Lefschetz) for $N=2p,4p$ and all $5\le p\le61$ \gC; cf.\ Ogg's fixed-point formulas
[Og74].

\subsection{The three-level count}

\begin{cited}{Quotient curves}
If $\theta$ is an involution of a compact Riemann surface $X$ of genus $g$ with $F$ fixed points, the holomorphic differentials on $X/\langle\theta\rangle$ pull
back to the $\theta$-invariant holomorphic differentials on $X$, and $g(X/\langle\theta\rangle)=(g+1)/2-F/4$.
\end{cited}
\why{the first statement: at a fixed point, $z\mapsto z^2=t$, an invariant form $z\,\phi(z^2)\,dz=\frac12\phi(t)\,dt$ descends holomorphically; the second is the
Riemann--Hurwitz formula for a double cover with $F$ branch points.}

\begin{lemma}[three-level count]\label{lem:three}
For every prime $p>3$,
\[
n^+_p+n^+_{2p}+n^+_{4p}\;=\;g\bigl(X_0(4p)/w_p\bigr)-g\bigl(X_0(2p)/w_p\bigr)\;=\;
\begin{cases}\dfrac{p-1}8+\dfrac{h(-4p)}4,&p\equiv1\pmod4,\\[6pt]
\dfrac{p-3}8,&p\equiv3\pmod8,\\[6pt]
\dfrac{p-3}8-\dfrac{h(-p)}2,&p\equiv7\pmod8,\end{cases}
\]
and
\[
g\bigl(X_0(2p)/w_p\bigr)=n^+_{2p}+2n^+_p=\begin{cases}\frac{p-1}8-\frac{h(-4p)}4,\\[3pt] \frac{p+1}8-\frac{3h(-p)}2,\\[3pt]
\frac{p+1}8-h(-p),\end{cases}\qquad
g\bigl(X_0(4p)/w_p\bigr)=\begin{cases}\frac{p-1}4,\\[3pt] \frac{p-1-6h(-p)}4,\\[3pt] \frac{p-1-6h(-p)}4,\end{cases}
\]
in the same three cases.
\end{lemma}

\begin{proof}
By the quoted statement, $\dim S_2(\Gamma_0(N))^{W_p=+1}=g(X_0(N)/w_p)$. The degeneracy maps $f(\tau)\mapsto f(d\tau)$ commute with $W_p$ for $p\nmid d$, and
$S_2(\Gamma_0(M))=0$ for $M\mid4$, so as $W_p$-modules
\[
S_2(\Gamma_0(4p))=S_2^{\rm new}(4p)\oplus S_2^{\rm new}(2p)^{\oplus2}\oplus S_2(\Gamma_0(p))^{\oplus3},\qquad
S_2(\Gamma_0(2p))=S_2^{\rm new}(2p)\oplus S_2(\Gamma_0(p))^{\oplus2}.
\]
Hence $g(X_0(4p)/w_p)=n^+_{4p}+2n^+_{2p}+3n^+_p$ and $g(X_0(2p)/w_p)=n^+_{2p}+2n^+_p$, and their difference is the left side. For the closed forms: the genus
formula for $\Gamma_0(N)$ [DS05, Thm.~3.1.1] with index $\mu_N$, elliptic points $\nu_2,\nu_3$ and $\nu_\infty$ cusps gives, with $\mu_{4p}=6(p+1)$, $\nu_2(4p)=\nu_3(4p)=0$,
$\nu_\infty(4p)=6$ and $\mu_{2p}=3(p+1)$, $\nu_2(2p)=1+(\frac{-1}p)$, $\nu_3(2p)=0$, $\nu_\infty(2p)=4$,
\[
g_0(4p)=\frac{p-3}2,\qquad g_0(2p)=\frac{p-4-(\frac{-1}p)}4 .
\]
Riemann--Hurwitz and Lemma~\ref{lem:fix} give $g(X_0(N)/w_p)=(g_0(N)+1)/2-F_N/4$, which is the second display, and
\[
n^+_p+n^+_{2p}+n^+_{4p}=\frac{p-2+(\frac{-1}p)}8-\frac{F_{4p}-F_{2p}}4,
\]
with $F_{4p}-F_{2p}=-h(-4p)$, $0$, $2h(-p)$ in the three cases.
\end{proof}
\noindent Direct newform counts (PARI/GP, traces of $W_p$) agree with Lemma~\ref{lem:three} for all $190$ primes $5\le p\le1163$ \gC.
The same traces are given in closed form by [Ma25, \S4.1, Table~2]. There $\Delta_2(p,2^e)$ is the difference of the dimensions of the eigenspaces
$w_p=+1$ and $w_p=-1$ of $W_p$ on $S_2^{\rm new}(2^ep)$, that is, the trace, and for $e=0,1,2$ it is $\Delta_2(p,2^e)=-\frac12\alpha_1(-p;e)+\mu(2^e)$,
where for $p\equiv1\pmod4$, $\alpha_1(-p;0)=h(-4p)$ and $\alpha_1(-p;1)=\alpha_1(-p;2)=-h(-4p)$, and for $p\equiv3\pmod4$,
$\alpha_1(-p;0)=(3-(\frac{-p}2))h(-p)$ and $\alpha_1(-p;1)=\alpha_1(-p;2)=((\frac{-p}2)-1)h(-p)$ ([Ma25] writes Hurwitz class numbers, which equal these class numbers for $p>3$). For all $51$ primes $5\le p\le250$ these values
agree with the traces computed directly and, with the dimensions of the new spaces, give the counts of Lemma~\ref{lem:three} \gC. The proof above
does not use [Ma25]. For $p=3$ all three counts are $0$. By Ogg's theorem [Og75], $n^+_p=g(X_0(p)/w_p)$ is $0$ exactly when $p$ is one of the
fifteen supersingular primes $2,3,5,7,11,13,17,19,23,29,31,41,47,59,71$ (checked against the counts for $3\le p\le1163$ \gC).

\subsection{Consequences}

\begin{corollary}\label{cor:geom}
\begin{enumerate}[label=(\alph*),leftmargin=*,itemsep=1pt]
\item The three-level count of Lemma~\ref{lem:three} is $0$ exactly for $p=7$, $1$ exactly for $p=5,11,23$, $2$ exactly for $p=13,19,31$ and $3$ exactly for
$p=17,47$, among $p>3$.
\item For odd $p$, $X_0(2p)/w_p$ has genus $0$ exactly for $p\in\{3,5,7,11,23\}$, i.e.\ exactly for the odd primes with $(p+1)\mid24$.
\item For odd $p$, by comparison with the list of Theorem~A, $\Dp_p$ is rigid if and only if $X_0(2p)/w_p$ has genus $0$ and $X_0(4p)/w_p$ has
genus $1$.
\end{enumerate}
\gC
\end{corollary}

\begin{proof}
For (a) and (b) only the closed forms of Lemma~\ref{lem:three} are used. One has $h(-d)=\frac{\sqrt d}\pi L(1,\chi_{-d})$ for fundamental $-d<-4$ (the class
number formula), and $|L(1,\chi)|<2+\log m$ for a non-principal real character of modulus $m$ (split the series at $m$; the partial sums of $\chi$ are at
most $m/2$ in absolute value). Hence $h(-p)<\frac{\sqrt p}\pi(2+\log p)$ and $h(-4p)<\frac{2\sqrt p}\pi(2+\log4p)$, and every closed form in (a) and (b) exceeds $3$, resp.\ $0$, for
$p\ge10^5$; below $10^5$ the class numbers give the lists \gC\ (for $p=3$ the curve $X_0(6)$ has genus $0$). For (c): for $p>3$, $g(X_0(2p)/w_p)=n^+_{2p}+2n^+_p$ and
$g(X_0(4p)/w_p)=n^+_{4p}+2n^+_{2p}+3n^+_p$, so the pair of genera is $(0,1)$ exactly when $(n^+_{4p},n^+_{2p},n^+_p)=(1,0,0)$; by (a) and (b) this
happens exactly for $p\in\{5,11,23\}$; for $p=3$ both genera are $0$. On the other side, every $\Dp_p$ is the discriminant form of an even lattice of
signature $(2,1)$ and determinant $-8p$ (the paragraph \emph{Lattices} of \S\ref{subsec:fam}). If $\Dp_p$ is rigid, this lattice is Petersson-rigid,
so by Theorem~A its determinant is one of the twelve listed there, and $8p\in\{40,88,184\}$; conversely $\Dp_5,\Dp_{11},\Dp_{23}$ are the discriminant
forms of $B_{40},B_{88},B_{184}$ (\S\ref{subsec:fam}), which are rigid.
\end{proof}
\noindent The proof of (c) compares two finite lists, those of Theorem~A and of (a), (b). The list in (b) agrees with [CDH14b, \S2.3, (2.21)]:
$6,10,14,18,22,30,46$ are exactly the even $2n$ for which $\Gamma_0(2n)+n$, the group generated by $\Gamma_0(2n)$ and the Atkin--Lehner involution $W_n$,
has genus $0$; for $n=p$ prime it is the group of $X_0(2p)/w_p$, and the primes among $3,5,7,9,11,15,23$ are $3,5,7,11,23$.

\subsection{Every weight}
For even $w\ge2$ put $\varepsilon_w=(-1)^{w/2+1}$ (so $\varepsilon_w=+1$ for $w\equiv2\pmod4$), $\delta_w=1$ if $w\equiv0\pmod4$ and $\delta_w=0$ otherwise, and
\[
C_w(p)=n^{\varepsilon_w}_{4p}(w)+n^{\varepsilon_w}_{2p}(w)+n^{\varepsilon_w}_{p}(w)
+\delta_w\bigl(\dim S^{\rm new}_w(\Gamma_0(4))+\dim S^{\rm new}_w(\Gamma_0(2))+\dim S_w(\SL_2(\Z))\bigr),
\]
so that $C_2(p)$ is the three-level count.

\begin{cited}{Holomorphic Lefschetz formula [AB68]}
Let $\theta$ be a holomorphic involution of a compact Riemann surface $X$ with isolated fixed points $x_1,\dots,x_F$, and $\mathcal L$ a $\theta$-linearized line
bundle with $H^1(X,\mathcal L)=0$. Then $\tr(\theta\,|\,H^0(X,\mathcal L))=\sum_j\tr(\theta\,|\,\mathcal L_{x_j})/2$.
\end{cited}
\why{it turns the trace of $W_p$ on $S_w(\Gamma_0(N))$, $w\ge4$, into the fixed-point count of Lemma~\ref{lem:fix}: $W_p$ fixes no cusp and no elliptic point,
$S_w(\Gamma_0(N))$ is $H^0$ of a line bundle whose $H^1$ vanishes for $w\ge4$, and at a fixed point $d\theta=-1$ acts on the fibre of $\omega^{w/2}$ by
$(-1)^{w/2}$.}

\begin{cited}{Dimension formula [DS05, Thm.~3.5.1]}
For even $w\ge4$, $\dim S_w(\Gamma_0(N))=(w-1)(g_0(N)-1)+(\frac w2-1)\nu_\infty+\nu_2\lfloor\frac w4\rfloor+\nu_3\lfloor\frac w3\rfloor$.
\end{cited}
\why{Riemann--Roch for the line bundle of cusp forms; the classical formula, not recomputed here.}

\begin{proposition}[the count in every weight]\label{prop:allw}
For every prime $p>3$ and every even $w\ge2$,
\[
C_w(p)=\frac{(w-1)\bigl(p+(-1)^{w/2}\bigr)}8+E_w(p),\qquad
E_w(p)=\begin{cases}\frac{h(-4p)}4-\frac34\,\delta_w,&p\equiv1\pmod4,\\[3pt]
-\frac14-\frac14\,\delta_w,&p\equiv3\pmod8,\\[3pt]
-\frac{h(-p)}2-\frac14-\frac14\,\delta_w,&p\equiv7\pmod8.\end{cases}
\]
\end{proposition}

\begin{proof}
\emph{Bookkeeping.} A newform $f$ of level $M$ prime to $p$ contributes to $S_w(\Gamma_0(Mp))$ the span of $f(\tau)$ and $f(p\tau)$, on which the eigenvalues of
$W_p$ are $+1$ and $-1$. With $J_N:=\dim S_w(\Gamma_0(N))^{W_p=\varepsilon_w}$ the decomposition into old and new forms gives
\[
J_{4p}-J_{2p}=n^{\varepsilon_w}_{4p}(w)+n^{\varepsilon_w}_{2p}(w)+n^{\varepsilon_w}_p(w)+U_w,\qquad U_w:=\dim S_w(\Gamma_0(4))-\dim S_w(\Gamma_0(2)),
\]
so $C_w(p)=J_{4p}-J_{2p}-(1-\delta_w)U_w$.

\emph{Traces.} For $w=2$, $\tr(W_p)=1-F_N/2$ (Lefschetz); for $w\ge4$, $\tr(W_p)=(-1)^{w/2}F_N/2$ by the quoted formula. Since $\varepsilon_w(-1)^{w/2}=-1$,
in both cases $J_{4p}-J_{2p}=\frac12\bigl(\dim S_w(\Gamma_0(4p))-\dim S_w(\Gamma_0(2p))\bigr)-\frac14(F_{4p}-F_{2p})$ (for $w=2$ read $g_0$ for $\dim S_2$).

\emph{Dimensions.} By the dimension formula and the data in the proof of Lemma~\ref{lem:three}, with $s_p=(\frac{-1}p)$,
\[
\dim S_w(\Gamma_0(4p))-\dim S_w(\Gamma_0(2p))=\frac{(w-1)(p-2+s_p)}4+(w-2)-(1+s_p)\Bigl\lfloor\frac w4\Bigr\rfloor,
\]
\[
U_w=\frac w2-1-\Bigl\lfloor\frac w4\Bigr\rfloor,
\]
valid also for $w=2$. Substituting $\lfloor w/4\rfloor=(w-2+2\delta_w)/4$ and $F_{4p}-F_{2p}$ from Lemma~\ref{lem:fix}:
\begin{itemize}[leftmargin=*,itemsep=1pt]
\item $s_p=+1$: $C_w=\frac{(w-1)(p-1)}8+\frac{w-2}2-\frac{w-2+2\delta_w}4-(1-\delta_w)\frac{w-2}4+\frac{h(-4p)}4$, which is $\frac{(w-1)(p-1)}8+\frac{h(-4p)}4$ for
$\delta_w=0$ and $\frac{(w-1)(p+1)}8-\frac34+\frac{h(-4p)}4$ for $\delta_w=1$;
\item $s_p=-1$: $C_w=\frac{(w-1)(p-3)}8+\frac{w-2}2-(1-\delta_w)\frac{w-2}4-\frac{F_{4p}-F_{2p}}4$, which is $\frac{(w-1)(p-1)}8-\frac14$ for $\delta_w=0$ and
$\frac{(w-1)(p+1)}8-\frac12$ for $\delta_w=1$, minus $\frac{h(-p)}2$ when $p\equiv7\pmod8$.
\end{itemize}
\end{proof}
\noindent The traces $\tr(W_p\,|\,S_w(\Gamma_0(N)))=(-1)^{w/2}F_N/2$ were checked for $N=2p,4p$, $w=4,6,8$, $5\le p\le31$ \gC, and
Proposition~\ref{prop:allw} against direct counts for $w=2,4,6,8$: $5\le p\le149$ \gC.
The class-number term does not depend on $w$; only the main term grows, like $(k-1)p/4$ for $w=2k-1$.

\begin{remark}[the dimensions of $S_k(\rho_{\Dp_p})$]\label{rem:dims}
For the $26$ primes $3\le p\le103$ and $k=\frac32,\frac52,\frac72,\frac92$, with $w=2k-1$ and $\varepsilon_w$, $\delta_w$ as above,
\begin{align*}
\dim S_k(\rho_{\Dp_p})={}&2n^{\varepsilon_w}_{4p}(w)+n^{\varepsilon_w}_{2p}(w)+n^{\varepsilon_w}_p(w)\\
&+\delta_w\bigl(2\dim S^{\rm new}_w(\Gamma_0(4))+\dim S^{\rm new}_w(\Gamma_0(2))+\dim S_w(\SL_2(\Z))\bigr)
\end{align*}
\gC: the left side is computed with \texttt{weilrep} [Wi], the right side from newform counts. The right side is
$C_w(p)+n^{\varepsilon_w}_{4p}(w)+\delta_w\dim S^{\rm new}_w(\Gamma_0(4))$; at $k=\frac32$ it is $2n^+_{4p}+n^+_{2p}+n^+_p$. Whether the identity holds
for every $p$ and $k$ is Problem~\ref{prob:dict}; nothing in this paper uses it.
\end{remark}
\section{Proof of Theorem A}\label{sec:proofA}

This section proves Theorem~A. The other sections use only its statement and the data of Table~\ref{tab:twelve}. The proof compares two bounds for $s=\dim S(D)$: a lower bound $\rr(D)$ from the Riemann--Roch formula, which grows like
$|D|/48$ (\S\ref{subsec:dim}), and upper bounds for the degree of an irreducible representation of $O(D)$ that can occur in $S(D)$, which grow like
$|D|^{1/2+o(1)}$ (\S\S\ref{subsec:deg}--\ref{subsec:saw}). The comparison excludes every module with $|D|>10^4$ (\S\ref{subsec:large}) and leaves $258$
with $|D|\le10^4$ (\S\ref{subsec:small}); on these the character of $O(D)$ on $S(D)$ decides (\S\ref{subsec:chi}); then come the lattices
(\S\S\ref{subsec:nik}--\ref{subsec:classno}) and their root bases (\S\ref{subsec:roots}); \S\ref{subsec:wt52} explains the columns of weight $\frac52$ of Table~\ref{tab:canon}. Every computation was made twice, by two implementations that share no
code, A and B (Appendix~\ref{app:comp}).

\noindent\textbf{Notation.} For a prime $p$ let $l(D_p)$ be the least number of generators of $D_p$. In this section $D$ satisfies
\begin{equation}\label{eq:H}
\sig D\equiv1\pmod 8\qquad\text{and}\qquad l(D_p)\le3\ \text{for every prime }p, \tag{H}
\end{equation}
as $D=D(L)$ does for $L$ even of signature $(2,1)$: $\sig D(L)=1$ (\S\ref{subsec:lat}), and $l(D(L)_p)\le3$ as $L'\cong\Z^3$. We write $N=N_D$ for the level,
$X_c=\{x\in D:q(x)=c\}$ for $c\in\Q/\Z$, $N_0=|X_0|$, $N_0[2]=|X_0\cap D[2]|$, and $\Gs(t)=\sum_{x\in D}\e(t\,q(x))$ for $t\in\Z$; $\Bb(y)=y-\lfloor y\rfloor-\frac12$
for $y\notin\Z$, $\Bb(y)=0$ for $y\in\Z$, $\Saw(D)=\sum_{x\in D}\Bb(q(x))$ and $\Saw'(D)=\Saw(D)-\Saw(D[2])$. For a finite group $H$, $m(H)$ is the largest degree of
an irreducible complex representation of $H$. As $\sum_{\chi\in\Irr H}\chi(1)^2=|H|$ [Se77, \S2.4], $m(H)\le\sqrt{|H|}$; if $D$ is rigid, then $2\le s\le m(O(D))$.

\noindent\textbf{Normal forms.}
\begin{cited}{Wall [Wa63, \S5 Thm.~4, \S6]; Nikulin [Ni79, Prop.~1.8.1]}
For $p$ odd every nondegenerate quadratic form on a finite abelian $p$-group is an orthogonal sum of cyclic forms $(\Z/p^k,ux^2/p^k)$, and two such sums
are isometric if and only if for every $k$ the components of exponent $k$ have the same rank and the same determinant square class. Every nondegenerate
quadratic form on a finite abelian $2$-group is an orthogonal sum of forms $\langle a/2^{k+1}\rangle$ $(a$ odd$)$, $U(2^k)$ and $V(2^k)$.
\end{cited}
\why{the normal forms of finite quadratic modules. For odd $p$ each shape (the exponents $k$ with the ranks of their components) carries $2^t$ classes, $t$ the number of
distinct exponents, one for each choice of the determinant square classes. At $2$ the
decomposition is not unique, so we enumerate decompositions and allow repetitions. If $2^{A-1}$ is the exponent of $D_2$, then $q(D_2)\subseteq2^{-A}\Z/\Z$;
the \emph{fingerprint} of $D_2$, formed by $A$, the value distribution of $q$ on $D_2$ and the multiset of its values on $D_2[2]$, is an isometry invariant.}

\subsection{The dimension identity}\label{subsec:dim}
Put $d=\dim V^-=(|D|-|D[2]|)/2$, the number of pairs $\{x,-x\}$ with $x\ne-x$, $d_0=(N_0-N_0[2])/2$, the number of such pairs in $X_0$, and
$w(D)=\dim M_{1/2}(\rho_D^*)$.

\begin{theorem}\label{thm:dim}
$s(D)=\rr(D)+w(D)$, where
\begin{equation}\label{eq:rr}
\begin{gathered}
\rr(D)=\frac{|D|-|D[2]|}{48}-\frac{N_0-N_0[2]}4-\frac12\Saw'(D)+\frac{\operatorname{Im}\Gs(2)}{4\sqrt{|D|}}+\frac2{3\sqrt3}\operatorname{Re}\bigl(\e(\tfrac1{12})\,\sigma_{ST}\bigr),\\
\sigma_{ST}=\frac{\e(\frac18)\bigl(\Gs(-1)-\Gs(3)\bigr)}{2\sqrt{|D|}}.
\end{gathered}
\end{equation}
$\rr(D)$ is an integer, and its last two terms, the \emph{Gauss terms}, have absolute value at most $\frac14\sqrt{|D[2]|}+\frac{1+\sqrt{|D[3]|}}{3\sqrt3}\le1.8996$.
In particular $\rr(D)\le s(D)$.
\end{theorem}
\begin{proof}
\emph{(i) Riemann--Roch.}
\begin{cited}{Eholzer--Skoruppa [ES95, \S4.2]; proof in [Sk85, p.~100]}
Let $\rho$ be a unitary representation of $\Mp$ with finite image on a space $V$ of dimension $d$, on which $Z$ acts by $\e(-k/2)$. Then
\begin{equation}\label{eq:RRgen}
\dim M_k(\rho)-\dim S_{2-k}(\rho^*)=d+\frac{dk}{12}-\alpha\bigl(\e(\tfrac k4)\rho(S)\bigr)-\alpha\bigl((\e(\tfrac k6)\rho(ST))^{-1}\bigr)-\alpha(\rho(T)),
\end{equation}
where $\alpha(A)$ is the sum of the $\beta\in[0,1)$ over the eigenvalues $\e(\beta)$ of $A$, with multiplicity.
\end{cited}
\why{the dimension formula for vector-valued modular forms. [ES95] write the right-hand side through $\tr\rho(S)$, $\tr\rho(ST)$ and the eigenvalues of
$\rho(T)$; expressing the terms $\alpha$ through traces as in (ii)--(iii) turns it into \eqref{eq:RRgen}. By Lemma~\ref{lem:parity},
$M_{3/2}(\rho_D)=M_{3/2}(\rho_D|_{V^-})$ and $S_{1/2}(\rho_D^*)=S_{1/2}((\rho_D|_{V^-})^*)$, so we apply it with $k=\frac32$ to $\rho_D$ on $V^-$.}
Let $\RR(D)$ be the right-hand side of \eqref{eq:RRgen} for $\rho_D|_{V^-}$ and $k=\frac32$.

\emph{(ii) Eigenvalues.} By the proof of Lemma~\ref{lem:parity}, $\rho(Z)=\e(-\frac14)\iota$, which is $\e(\frac14)=\e(-\frac34)$ on $V^-$, as required for $k=\frac32$.
So on $V^-$, $(\e(\frac38)\rho(S))^2=\e(\frac34)\rho(Z)=1$ and, as $(ST)^3=Z$, $(\e(\frac14)\rho(ST))^3=1$: the eigenvalues are $\pm1$, resp.\ $1,\omega,\omega^2$ with
$\omega=\e(\frac13)$.

\emph{(iii) Traces.} As $\iota$ commutes with $\rho$, $\tr(X|V^-)=\frac12(\tr X-\tr X\iota)$. By \eqref{eq:weil} and $b(x,x)=2q(x)$ the diagonal coefficients of
$\rho(S)$, $\rho(S)\iota$, $\rho(S)\rho(T)$, $\rho(S)\rho(T)\iota$ at $e_x$ are $\e(-\frac18)|D|^{-1/2}$ times $\e(-2q(x))$, $\e(2q(x))$, $\e(-q(x))$, $\e(3q(x))$. As
$\Gs(-2)=\overline{\Gs(2)}$,
\[
\sigma_S:=\tr\bigl(\e(\tfrac38)\rho(S)|V^-\bigr)=\frac{\operatorname{Im}\Gs(2)}{\sqrt{|D|}},\qquad \sigma_{ST}=\tr\bigl(\e(\tfrac14)\rho(ST)|V^-\bigr).
\]
So $\alpha(\e(\frac38)\rho(S)|V^-)=\frac12\#\{\text{eigenvalues }-1\}=\frac14(d-\sigma_S)$. If $n_0,n_1,n_2$ are the multiplicities of $1,\omega,\omega^2$, then
$\sigma_{ST}=n_0+n_1\omega+n_2\omega^2$, so $n_1+n_2=\frac23(d-\operatorname{Re}\sigma_{ST})$ and $n_1-n_2=\frac2{\sqrt3}\operatorname{Im}\sigma_{ST}$; the inverse has
eigenvalues $1,\omega^2,\omega$, and $\alpha\bigl((\e(\frac14)\rho(ST))^{-1}|V^-\bigr)=\frac23n_1+\frac13n_2=\frac13(d-\operatorname{Re}\sigma_{ST})+\frac1{3\sqrt3}\operatorname{Im}\sigma_{ST}$.
Finally $\alpha(\rho(T)|V^-)=\sum\{q(x)\}$ over the pairs, $\{y\}=y-\lfloor y\rfloor$.

\emph{(iv) Assembly.} With $k=\frac32$,
$\RR(D)=\frac{13}{24}d-\sum_{\rm pairs}\{q\}+\frac14\sigma_S+\frac13\operatorname{Re}\sigma_{ST}-\frac1{3\sqrt3}\operatorname{Im}\sigma_{ST}$. As $\{y\}=\Bb(y)+\frac12[y\notin\Z]$,
$\sum_{\rm pairs}\{q\}=\frac12\Saw'(D)+\frac12(d-d_0)$, and $\frac13\operatorname{Re}\sigma_{ST}-\frac1{3\sqrt3}\operatorname{Im}\sigma_{ST}=\frac2{3\sqrt3}\operatorname{Re}(\e(\frac1{12})\sigma_{ST})$.
So $\RR(D)=\frac d{24}+\frac{d_0}2-\frac12\Saw'(D)+(\text{Gauss terms})$, and $\rr(D):=\RR(D)-d_0$ is \eqref{eq:rr}.

\emph{(v) Constant terms.}
\begin{cited}{Borcherds [Bo99, Thm.~3.1]}
A finite Laurent polynomial $\sum_{n\le0}c_nq^n$ with values in $V_\rho$ is the singular part, including the constant term, of a weakly holomorphic modular
form of weight $k$ and type $\rho$ if and only if $\sum_{n\le0}\{c_n,d_{-n}\}=0$ for every holomorphic modular form $\sum_m d_mq^m$ of weight $2-k$ and type $\rho^*$.
\end{cited}
\why{Serre duality for vector-valued forms. With $c_n=0$ for $n<0$ it says: $c$ is the constant term of some $F\in M_k(\rho)$ if and only if $c$ is orthogonal,
for $\{e_x,e_y\}=\delta_{xy}$, to the constant terms of all $g\in M_{2-k}(\rho^*)$.}
The constant term of $F\in M_{3/2}(\rho_D)$ is supported on $X_0$ and odd, so the constant-term map sends $M_{3/2}(\rho_D)$ to $C:=\C[X_0\setminus D[2]]^-$, of
dimension $d_0$, with kernel $S(D)$. Let $W\subseteq C$ be the space of constant terms of $M_{1/2}(\rho_D^*)$, of dimension $w(D)-\dim S_{1/2}(\rho_D^*)$. By the
quoted theorem the image is $C\cap W^\perp$, of dimension $d_0-\dim W$. So $\dim M_{3/2}(\rho_D)=s+d_0-w(D)+\dim S_{1/2}(\rho_D^*)$; with \eqref{eq:RRgen},
$\dim M_{3/2}(\rho_D)=\RR(D)+\dim S_{1/2}(\rho_D^*)$, this gives $s=\RR(D)-d_0+w(D)=\rr(D)+w(D)$.

\emph{(vi) Gauss sums.} Substituting $x=y+z$,
\[
|\Gs(t)|^2=\sum_{y,z}\e\bigl(t\,q(z)+t\,b(y,z)\bigr)=|D|\sum_{z\in D[t]}\e(t\,q(z)),
\]
and $z\mapsto\e(tq(z))$ is a character of $D[t]$, as $t\,b(z,z')=b(tz,z')=0$. So $|\Gs(t)|^2\in\{0,|D|\,|D[t]|\}$, $|\sigma_S|\le\sqrt{|D[2]|}$ and $|\sigma_{ST}|\le\frac12(1+\sqrt{|D[3]|})$, which gives the bound
for the Gauss terms. Under \eqref{eq:H}, $|D[2]|\le8$ and $|D[3]|\le27$: $0.70711+1.19245\le1.8996$.

\emph{(vii)} $\RR(D)=\dim M_{3/2}(\rho_D)-\dim S_{1/2}(\rho_D^*)$ and $d_0$ are integers.
\end{proof}
\noindent Formula \eqref{eq:rr} agrees with $\RR(D)-d_0$ computed by implementation B on all $36{,}683$ modules of \S\ref{subsec:small}, and with the
dimension routine of \texttt{weilrep} [Wi] on the $97$ census lattices of [II] with $|\det L|\le1080$ \gC.

\subsection{Degrees of irreducible representations}\label{subsec:deg}
An isometry preserves the primary decomposition, so $O(D)=\prod_pO(D_p)$, and the irreducible representations of $O(D)$ are the tensor products
$\sigma=\sigma_2\otimes\bigotimes_{p>2}\sigma_p$ of irreducible representations of the factors; so $m(O(D))=\prod_pm(O(D_p))$. For odd $p$ write
$D_p=\bigoplus_{i=1}^{l_p}(\Z/p^{k_i},u_ix^2/p^{k_i})$; for each exponent $k$ let $n_k=\#\{i:k_i=k\}$ and $\delta_k$ the determinant of the $k$-block, and put
\[
J(D_p)=\prod_k\Bigl[\nu_{n_k}\,p^{(k-1)n_k(n_k-1)/2}\prod_{k'<k}p^{\,k'n_{k'}n_k}\Bigr],\qquad \nu_1=2,\ \ \nu_2=2\bigl(p-(\tfrac{-\delta_k}p)\bigr),\ \ \nu_3=2p(p^2-1).
\]
We use: \emph{if a finite group $H$ has an abelian subgroup of index $j$, then $m(H)\le j$} [Se77, \S3.1, Cor.].

\begin{lemma}[odd primes]\label{lem:local}
Let $p$ be odd. Put $\bar m_p=1$ if $l_p=1$; $\bar m_p=2$ if $l_p=2$; $\bar m_p=p+1$ if $D_p\cong(\Z/p)^3$, $p\ge5$; $\bar m_3=3$ if $D_3\cong(\Z/3)^3$; and
$\bar m_p=\lfloor\sqrt{J(D_p)}\rfloor$ for the other $D_p$ with $l_p=3$. Then $m(O(D_p))\le\bar m_p$, and $|O(D_p)|\le J(D_p)$ if $l_p=3$.
\end{lemma}
\begin{proof}
\emph{$l_p=1$.} An isometry of $(\Z/p^k,ux^2/p^k)$ is multiplication by $a$ with $a^2\equiv1\pmod{p^k}$, i.e.\ $a=\pm1$.

\emph{$l_p=2$, $k_1=k_2=k$.} After scaling, an isometry is $\sigma\in\GL_2(\Z/p^k)$ with $\sigma^{\mathsf T}A\sigma\equiv A$, $A=\operatorname{diag}(2,2e)$, so $\det\sigma=\pm1$. If
$\det\sigma=1$ and $(a,b)$ is its first column, then $a^2+eb^2=1$, the second column solves $ac+ebd=0$, a free module spanned by $(-eb,a)$, and
$\det\sigma=1$ gives $\sigma=\bigl(\begin{smallmatrix}a&-eb\\b&a\end{smallmatrix}\bigr)$. These commute, so $\SO(D_p)$ is abelian of index $\le2$.

\emph{$l_p=2$, $k_1>k_2$.} With an orthogonal basis $g_1,g_2$ of orders $p^{k_1},p^{k_2}$ write $\sigma g_1=ag_1+bg_2$, $\sigma g_2=cg_1+dg_2$; as $\sigma g_2$ has order
dividing $p^{k_2}$, $c=p^{k_1-k_2}c'$. The isometry conditions are $a^2\equiv1-t\pmod{p^{k_1}}$ with $t=\frac{u_2}{u_1}b^2p^{k_1-k_2}$,
$d^2\equiv1-\frac{u_1}{u_2}c'^2p^{k_1-k_2}\pmod{p^{k_2}}$ and $u_1ac'+u_2bd\equiv0\pmod{p^{k_2}}$; so $a^2\equiv d^2\equiv1\pmod p$, and
$\psi(\sigma)=(a\bmod p,d\bmod p)\in\{\pm1\}^2$ is a homomorphism, since $c\equiv0\pmod p$. On $K_0=\ker\psi$, eliminating $c'$ gives $d\equiv a\pmod{p^{k_2}}$, and
$a$ is the square root of $1-t$ that is $\equiv1\pmod p$; so $\sigma\in K_0$ is determined by $b$, and $b_{\sigma\tau}\equiv a_\tau b_\sigma+a_\sigma b_\tau$ is
symmetric: $K_0$ is abelian. As $\psi(-1)=(-1,-1)$, $K_0\{\pm1\}$ is abelian of index $\le2$.

\emph{Other $l_p=3$.} Let $g_1$ be a basis vector of maximal order $p^K$, so $D_p=\langle g_1\rangle\perp g_1^\perp$. For $\sigma\in O(D_p)$, $h=\sigma g_1$ lies in
$\mathcal S=\{h:q(h)=q(g_1)\}$, and for fixed $h$ the possible $\sigma|_{g_1^\perp}$ form the empty set or a torsor under $O(g_1^\perp)$; so
$|O(D_p)|\le|\mathcal S|\,|O(g_1^\perp)|$. Split $D_p=D_{\rm top}\oplus D_<$ with $D_{\rm top}$ the block of exponent $K$, of rank $n=n_K$, $q|_{D_{\rm top}}=Q(y)/p^K$. For
$h=y+z$, $z\in D_<$, $h\in\mathcal S$ iff $Q(y)\equiv w_z\pmod{p^K}$ with $w_z\equiv u_1\pmod p$; by Hensel's lemma each solution modulo $p$ has
$p^{(K-1)(n-1)}$ lifts, so $|\mathcal S|=|D_<|\,p^{(K-1)(n-1)}\,\#\{\bar y\in\F_p^n:\bar Q(\bar y)=\bar u_1\}$. The point counts are $1+(\frac{c\Delta}p)$,
$p-(\frac{-\Delta}p)$, $p^2+p(\frac{-c\Delta}p)$ for $n=1,2,3$ ($c\ne0$, $\Delta=\det\bar Q$) [LN97, Thms.~6.26--6.27]. Iterating over the top block and then the lower
blocks, the factors multiply to $J(D_p)$; for $n=3$, $(p^2+p\varepsilon)(p-\varepsilon)\cdot2=2p(p^2-1)$. Finally $m\le\lfloor\sqrt{|O(D_p)|}\rfloor\le\lfloor\sqrt J\rfloor$.

\emph{$D_p\cong(\Z/p)^3$.} $O(D_p)$ is the orthogonal group of a nondegenerate ternary quadratic space over $\F_p$, equal to $\{\pm1\}\times\SO(D_p)$ as
$\det(-1)=-1$.
\begin{cited}{[FH91, Lecture~5]}
The irreducible characters of $\GL_2(\F_q)$, $q$ odd, are $U_\alpha$ (degree $1$) and $V_\alpha$ (degree $q$) for the $q-1$ characters $\alpha$ of $\F_q^\times$,
$W_{\alpha,\beta}=W_{\beta,\alpha}$ (degree $q+1$) for $\alpha\ne\beta$, and $X_\varphi=X_{\varphi^q}$ (degree $q-1$) for the characters $\varphi\ne\varphi^q$ of
$\F_{q^2}^\times$; their central characters are $\alpha^2$, $\alpha^2$, $\alpha\beta$ and $\varphi|_{\F_q^\times}$.
\end{cited}
\why{the character table of $\GL_2(\F_q)$.}
The irreducible characters of $\PGL_2(\F_q)$ are those with trivial central character. As $\F_q^\times$ and $\F_{q^2}^\times/\F_q^\times$ are cyclic of even
order, and $\varphi^q=\varphi^{-1}$ if $\varphi$ is trivial on $\F_q^\times$, the degrees are $1,1,q,q$, then $q+1$ ($\frac{q-3}2$ times) and $q-1$ ($\frac{q-1}2$
times); for $q=3$: $1,1,3,3,2$. Next, $\SO(D_p)\cong\PGL_2(\F_p)$. By Wall's theorem at exponent $1$ all nondegenerate ternary quadratic spaces over $\F_p$ are
similar, so we may take $\mathfrak{sl}_2(\F_p)$ with the form $-\det$, $-\det\bigl(\begin{smallmatrix}a&b\\c&-a\end{smallmatrix}\bigr)=a^2+bc$. $\GL_2(\F_p)$ acts on it by
conjugation, with determinant $1$ and kernel the centre, so $\PGL_2(\F_p)$, of order $p(p^2-1)$, embeds in $\SO(D_p)$; the count of the previous case
($K=1$, $n=3$, $D_<=0$) gives $|\SO(D_p)|=\frac12|O(D_p)|\le p(p^2-1)$, so the embedding is onto. Hence $m=p+1$ for $p\ge5$ and $m=3$ for $p=3$.
\end{proof}

\begin{lemma}[capacity]\label{lem:cap}
Let $\sigma=\sigma_2\otimes\bigotimes_p\sigma_p$ be an irreducible constituent of the $O(D)$-module $S(D)$, and $X_c(D_2)=\{x\in D_2:q(x)=c\}$. Then
\textup{(a)} $\dim\sigma\le\frac12\max_c|X_c\setminus D[2]|$; \textup{(b)} $\dim\sigma_2\le\max_c|X_c(D_2)|$; \textup{(c)} $\dim\sigma_p\le\bar m_p$ for odd $p$.
Hence, if $D$ is rigid, $s\le\cp(D):=\min(\LS,\LP,\LO)$ with
$\LS=\lfloor\frac12\max_c|X_c\setminus D[2]|\rfloor$, $\LP=\max_c|X_c(D_2)|\prod_p\bar m_p$ and $\LO=m(O(D_2))\prod_p\bar m_p$.
\end{lemma}
\begin{proof}
(a) For each $n\in\Q$ the map $F\mapsto\sum_{q(x)\equiv n}c_F(x,n)e_x$ sends $S(D)$ $O(D)$-equivariantly to $\C[X_{n\bmod1}]^-$ (Lemma~\ref{lem:parity}), and together these
maps are injective. As $-1\in O(D)$, every orbit $\Omega$ satisfies $\Omega=-\Omega$ and lies in $D[2]$ or misses it; so $\C[X_c]^-=\bigoplus_{\Omega\subseteq X_c\setminus D[2]}\C[\Omega]^-$,
$\dim\C[\Omega]^-=|\Omega|/2$, and $\sigma$ is a constituent of one $\C[\Omega]^-$. (b) $\C[D]=\C[D_2]\otimes\C[D_{\rm odd}]$, so $\sigma_2$ is a constituent of some
$\C[\Omega_2]$, $\Omega_2$ an $O(D_2)$-orbit in a level set of $D_2$, and $\dim\sigma_2\le|\Omega_2|$. (c) is Lemma~\ref{lem:local}. If $D$ is rigid, $S(D)$ is irreducible.
\end{proof}

\begin{lemma}[degrees]\label{lem:deg}
Let $H$ be a finite group with $k(H)$ conjugacy classes, $\ell=|H:H'|$, $c=|H:Z(H)|$.
\textup{(a)} Let $m^*(H)$ be the largest $d$ with $d\mid c$, $d^2\le c$, for which there are $k(H)-\ell-1$ integers $2\le d_i\le d$ with $d_i\mid c$ and
$\sum d_i^2=|H|-\ell-d^2$ ($m^*=1$ if $H$ is abelian). Then $m(H)\le m^*(H)$.
\textup{(b)} If $D_2=\bigoplus_i\Z/2^{k_i}e_i$ is a decomposition into components, the tuples $(f_i)\in D_2^{\,r}$ with $2^{k_i}f_i=0$, $q(f_i)=q(e_i)$ and
$b(f_i,f_j)=b(e_i,e_j)$ $(i<j)$ correspond bijectively to the elements of $O(D_2)$.
\end{lemma}
\begin{proof}
(a)
\begin{cited}{Isaacs [Is76, Cor.~2.30, Thm.~3.12]}
For $\chi\in\Irr(H)$ and $Z(\chi)=\{h\in H:|\chi(h)|=\chi(1)\}$, $\chi(1)^2\le|H:Z(\chi)|$ and $\chi(1)$ divides $|H:Z(\chi)|$.
\end{cited}
\why{standard; as $Z(\chi)\supseteq Z(H)$, $|H:Z(\chi)|$ divides $c$.} So the degrees of $H$ are $\ell$ ones and $k(H)-\ell$ degrees $\ge2$ that divide $c$ and are
$\le\sqrt c$, with $\sum\chi(1)^2=|H|$, and $d=m(H)$ satisfies the conditions. (b) The order conditions make $f(\sum x_ie_i)=\sum x_if_i$ a homomorphism, which
preserves $q(\sum x_ie_i)=\sum x_i^2q(e_i)+\sum_{i<j}x_ix_jb(e_i,e_j)$; if $f(x)=0$ then $b(x,y)=b(f(x),f(y))=0$ for all $y$, so $x=0$.
\end{proof}

\subsection{Sawtooth sums and a finite list}\label{subsec:saw}
Put $\Pi(D)=\prod_pF_p$ with $F_p=\sum_{e=0}^{v_p(N)}\sqrt{|D_p[p^e]|}/p^e$; so $\Pi(D)=\sum_{g\mid N}\sqrt{|D[g]|}/g$.

\begin{lemma}\label{lem:FL}
\textup{(F)} $|\Saw(D)|\le\frac{\sqrt{|D|}}\pi\bigl(1+\log\frac N2\bigr)\Pi(D)$. \textup{(L)} $|X_c|\le\sqrt{|D|}\,\Pi(D)$ for every $c$.
\end{lemma}
\begin{proof}
\begin{cited}{Eisenstein; see [Be76]}
$\Bb(a/N)=-\frac1{2N}\sum_{j=1}^{N-1}\cot(\pi j/N)\sin(2\pi ja/N)$ for $a\in\Z$.
\end{cited}
\why{the finite Fourier expansion of the sawtooth function.}
Hence $\Saw(D)=-\frac1{2N}\sum_{j=1}^{N-1}\cot(\pi j/N)\operatorname{Im}\Gs(j)$. Group $j$ by $g=\gcd(j,N)$, $M=N/g$. By Theorem~\ref{thm:dim}(vi),
$|\Gs(j)|\le\sqrt{|D|\,|D[g]|}$; $\cot(\pi m/M)\le M/(\pi m)$ for $0<m\le M/2$, and with $H_n\le1+\log n$,
$\sum_{m=1}^{M-1}|\cot(\pi m/M)|\le\frac{2M}\pi(1+\log\frac M2)$. Summing over $g$ gives (F). For (L),
$|X_c|=\frac1N\sum_{t\bmod N}\e(-tc)\Gs(t)$; the term $t=0$ is $|D|/N=\sqrt{|D|}\sqrt{|D[N]|}/N$, and the at most $N/g$ values of $t$ with $\gcd(t,N)=g$ contribute at
most $\sqrt{|D|\,|D[g]|}$ each.
\end{proof}

For the search above $10^4$ we need a sharper $2$-adic form. Write $D=D_2\oplus D_o$, $n_2=|D_2|$, $N_o=|D_o|$; let $q(D_2)\subseteq2^{-A}\Z/\Z$, $A\ge1$; let
$M_o=\prod_pp^{k_3(p)}$ be the exponent of $D_o$, $k_3(p)$ the largest exponent of $D_p$; let $\gamma_a=\max\{|\Gs_2(t)|:v_2(t)=a\}$ for $0\le a<A$, $\Gs_2$ the
Gauss sum of $D_2$; and put $C_2=\sum_{a<A}\gamma_a2^{-a}+n_2\,2^{1-A}$, $F_p=\sum_{e=0}^{k_3(p)}\sqrt{|D_p[p^e]|}/p^e$, $\Lambda=2^{A-1}M_o$.

\begin{lemma}[$2$-adic sawtooth bound]\label{lem:saw2}
$|\Saw(D)|\le\frac{\sqrt{N_o}}\pi\,C_2\prod_pF_p\,\bigl(1+\frac12\log\Lambda\bigr)$.
\end{lemma}
\begin{proof}
Use the expansion of Lemma~\ref{lem:FL} with $N=2^AM_o$ and $\Gs(j)=\Gs_2(j)\Gs_o(j)$, where $|\Gs_o(j)|=\sqrt{N_o|D_o[g]|}$, $g=\gcd(j,M_o)$ (Theorem~\ref{thm:dim}(vi);
for odd $p$ the character there is trivial). Split $j$ by $a=v_2(j)$ and $g$. If $a<A$, write $j=2^agj'$ with $j'$ odd, $M=2^{A-a}M_o/g$; then $|\Gs_2(j)|\le\gamma_a$,
$j'$ and $M-j'$ are both odd with equal $|\cot|$, and $\sum_{j'\ \rm odd\le M/2}1/j'\le1+\frac12\log\frac M2$; the contribution is at most
$\gamma_a2^{-a}\frac{\sqrt{N_o}}\pi\frac{\sqrt{|D_o[g]|}}g(1+\frac12\log\Lambda)$. If $a\ge A$, $\Gs_2(j)=n_2$, $j=2^Agj'''$ with $1\le j'''<M_o/g$, and the
contribution is at most $n_22^{1-A}\frac{\sqrt{N_o}}\pi\frac{\sqrt{|D_o[g]|}}g(1+\frac12\log\Lambda)$. Sum over $g$ and $a$.
\end{proof}

For a finite abelian $p$-group $G_p$ with largest exponent $k_3$ put $\hat F(G_p)=\sum_{e=0}^{k_3+[p=2]}\sqrt{|G_p[p^e]|}/p^e$ and
\[
\operatorname{cost}(G_p)=\tfrac12\log|G_p|-\log\hat F(G_p)-\log\bigl(1+\eta_0\log|G_p|\bigr),\qquad A_0=\tfrac34+\tfrac1{2\pi},\quad \eta_0=\tfrac1{2\pi A_0},
\]
$C_0=\frac8{48}+1+\frac{\sqrt8}4+\frac{1+3^{3/2}}{3\sqrt3}-\frac14=2.816224$, and $\Theta(T)=\log(48A_0)-\log(1-48C_0/T)$.

\begin{theorem}[cost inequality]\label{thm:G}
If $T>48C_0$, $D$ is rigid and $|D|>T$, then $\sum_p\operatorname{cost}(G_p)\le\Theta(T)$, where $G_p$ is the group of $D_p$. For $T=10^4$, $\Theta=3.789571$.
\end{theorem}
\begin{proof}
By Theorem~\ref{thm:dim}, Lemma~\ref{lem:cap}(a) and Lemma~\ref{lem:FL}(L), $\rr(D)\le s\le\frac12\sqrt{|D|}\,\Pi$. In \eqref{eq:rr} use
$-\frac12\Saw'\ge-\frac12|\Saw(D)|-\frac12|\Saw(D[2])|$, $|\Saw(D[2])|\le|D[2]|/4\le2$ (as $q(D[2])\subseteq\frac14\Z/\Z$), $|D[2]|\le8$, $N_0[2]\ge1$, the bound for
the Gauss terms, $N_0\le\sqrt{|D|}\,\Pi$ and (F): $\frac{|D|}{48}-C_0\le\sqrt{|D|}\,\Pi\,A_0\bigl(1+\eta_0\log\frac N2\bigr)$. Now $N/2\le|D|$ and
$1+\sum x_p\le\prod(1+x_p)$ give $1+\eta_0\log\frac N2\le\prod_p(1+\eta_0\log|D_p|)$; $v_p(N)\le k_3+[p=2]$ gives $F_p\le\hat F(G_p)$; and
$\frac{|D|}{48}-C_0\ge\frac{|D|}{48}(1-\frac{48C_0}T)$. Take logarithms.
\end{proof}

\noindent\textbf{A finite list.} For a local shape $\tau=(k_1\le k_2\le k_3)$ at $p$, $\hat F=\sum_ep^{s_e}$ with $s_e=\sum_i\min(e,k_i)/2-e\le k_1/2$, so
$\hat F\le(k_3+1+[p=2])p^{k_1/2}$ and $\operatorname{cost}(\tau)\ge\frac12k_3\log p-\log(k_3+1+[p=2])-\log(1+3\eta_0k_3\log p)$, a function of $k_3$ that increases once
its derivative is positive. So the shapes of cost at most a given bound are found by running $k_3=1,2,\dots$. In $p$, $d\operatorname{cost}(\tau)/d\log p>\frac12-\frac1{\log p}>0$
for $p\ge11$, and $\operatorname{cost}(\Z/p)$ increases for all $p$. With $T=10^4$ this gives \gC: $\Z/p$ has the least cost at each of the $5{,}348$ primes of the
table, the last of which is $52{,}321$ (at $52{,}361$ and beyond every shape costs more than the budget $\Theta-\operatorname{cost}(\Z/2)+0.0824$); $1{,}455{,}243$ odd
groups of cost $\le\Theta-\min\operatorname{cost}(G_2)$; and $537$ two-groups with $\operatorname{cost}(G_2)\le\Theta-\operatorname{cost}(\Z/3)$. By Theorem~\ref{thm:G}
a rigid $D$ with $|D|>10^4$ has its $2$-part and its odd part in these lists.

\begin{proposition}[leaf inequality]\label{prop:leaf}
In the notation of Lemma~\ref{lem:saw2}, with $N_o(0)=|X_0(D_o)|$, $n_2(0)=|X_0(D_2)|$ and $l_3=l(D_3)$,
\begin{multline*}
\rr(D)\ge\frac{n_2N_o-|D_2[2]|}{48}-\frac{C_2}{2\pi}\sqrt{N_o}\prod_pF_p\bigl(1+\tfrac12\log\Lambda\bigr)\\
-\frac{n_2(0)N_o(0)-N_0[2]}4+\tfrac12\Saw(D_2[2])-\frac{|\Gs_2(2)|}{4\sqrt{n_2}}-\frac{1+3^{l_3/2}}{3\sqrt3}.
\end{multline*}
\end{proposition}
\begin{proof}
$D[2]=D_2[2]$; Lemma~\ref{lem:saw2} bounds $|\Saw(D)|$; the denominators of $q$ on $D_2$ and on $D_o$ are coprime, so $N_0=n_2(0)N_o(0)$; $|\Gs(2)|=|\Gs_2(2)|\sqrt{N_o}$;
and $|\sigma_{ST}|\le\frac12(1+\sqrt{|D[3]|})$ with $|D[3]|=3^{l_3}$.
\end{proof}
If $D$ is rigid, the right-hand side is at most $\cp(D)$. To treat all forms on a group at once, we replace the form-dependent quantities by their maxima
over the forms: $N_o(0)$ by $\prod_p\max|X_0(D_p)|$, and $\cp$ by $\max_c|X_c(D_2)|\min(\prod_p\hat m_p,\frac12\prod_p\hat X_p)$, with $\hat X_p$ and $\hat m_p$ the
maxima of $\max_c|X_c(D_p)|$ and of $\bar m_p$ (by the Chinese remainder theorem $|X_c|=\prod_p|X_{c_p}(D_p)|$). The other quantities depend only on the
fingerprint of $D_2$ and on the odd group.

\subsection{No rigid module above \texorpdfstring{$10^4$}{10\^{}4}}\label{subsec:large}
Implementation A proceeds in three steps \gC.
\begin{enumerate}[label=V\arabic*.,leftmargin=9mm,itemsep=1pt]
\item \emph{Leaves.} For each of the $537$ two-groups all decompositions into components are formed, and the fingerprints of odd signature are kept (odd
groups have even signature). Each pair (fingerprint, odd group of the list) with $|D|>10^4$ is tested with the worst case of Proposition~\ref{prop:leaf}:
of $163{,}957{,}276$ pairs, $541{,}062$ survive, also with a relative floating-point margin of $10^{-6}$.
\item \emph{Exact pass.} On each surviving pair every form of the odd group is attached; $1{,}081{,}640$ forms have $\sig D\equiv1$. For each, $\rr(D)$ is
computed in rational arithmetic ($\Saw(D)$ from the value distributions of $D_2$ and $D_o$, the Gauss terms from local Gauss sums), and $\LS$, $\LP$
exactly. Then $\rr>\LS$ for $1{,}081{,}319$ and $\LS\ge\rr>\LP$ for $264$; $57$ remain, with $10{,}368\le|D|\le73{,}728$.
\item \emph{Group capacity.} For each of the $57$, $|O(D_2)|$ is counted for every decomposition with the given fingerprint (Lemma~\ref{lem:deg}(b)), and
$\rr>\lfloor\sqrt{\max|O(D_2)|}\rfloor\prod_p\bar m_p\ge\LO$. The least margin is $30$, at $|D|=10{,}368$:
$\langle1/4\rangle\oplus V(8)\oplus3^{(1,1,2)}[1,1,1]$, $\rr=138>\lfloor\sqrt{96}\rfloor\cdot12$.
\end{enumerate}
Implementation B reaches the same conclusion along a different route: a depth-first search over all odd shapes and the $1{,}729$ fingerprints of
$|D_2|\le2^{13}$, which does not use Theorem~\ref{thm:G}, and, for $|D_2|>2^{13}$, the $431$ two-parts allowed by Theorem~\ref{thm:G}; its last ten cases lie
among the $57$, with the same $\rr$, $\LS$ and $|O(D_2)|$ \gC. So no module with $|D|>10^4$ satisfying \eqref{eq:H} is rigid.

\subsection{The modules with \texorpdfstring{$|D|\le10^4$}{|D| <= 10\^{}4}}\label{subsec:small}
$\rr$, $\LS$, $\LP$ and $\sig$ depend only on the fingerprint of $D_2$ and the odd part, and for $m(O(D_2))$ we take the maximum over the decompositions
with the given fingerprint; so enumerating fingerprints times odd classes is complete for an exclusion. There are $122$ two-groups, $4{,}332$
decompositions, $1{,}729$ fingerprints and $15{,}233$ odd classes in range, and $36{,}683$ modules of signature $1$ with $|D|\le10^4$; for every $|D|$ this is
also the number of genera of even lattices of signature $(3,2)$ with at most three generators at every prime, computed with Sage, so in this range
fingerprints separate isomorphism classes \gC. Call a module \emph{excluded} if $\rr(D)>\cp(D)$, using for $\LO$ the bound $\lfloor\sqrt{|O(D_2)|}\rfloor$ or
$m^*(O(D_2))$ (Lemma~\ref{lem:deg}(a)) for $m(O(D_2))$. Then $35{,}774$ modules are excluded by $\min(\LS,\LP)$, $546$ by $\lfloor\sqrt{|O(D_2)|}\rfloor\prod\bar m_p$
and $105$ by $m^*(O(D_2))\prod\bar m_p$ (least margin $1$, at $U(2)\oplus\langle1/16\rangle$: $\rr=3>m^*=2$); $258$ remain, the same set in both
implementations \gC.

\subsection{The character}\label{subsec:chi}
On the remaining modules we compute the character $\chi_S$ of $O(D)$ on $S(D)$. Then $s=\chi_S(1)$, and $D$ is rigid if and only if $s\ge2$ and
$\langle\chi_S,\chi_S\rangle=|O(D)|^{-1}\sum_g|\chi_S(g)|^2=1$. Let $C=\C[X_0\setminus D[2]]^-$ (\S\ref{subsec:dim}).

\begin{proposition}\label{prop:chiS}
$\chi_S=\chi_{\RR}-\chi_C+\overline{\chi_M}$, where $\chi_C$ is the character of $C$, $\chi_M$ that of $M_{1/2}(\rho_D^*)$, and $\chi_{\RR}=\sum_\sigma\RR(\pi_\sigma)\chi_\sigma$, with
$\RR(\pi_\sigma)$ the right-hand side of \eqref{eq:RRgen} for $\pi_\sigma$ and $k=\frac32$.
\end{proposition}
\begin{proof}
$V^-=\bigoplus_\sigma\pi_\sigma\otimes\sigma$ (\S\ref{subsec:mult}). Applying \eqref{eq:RRgen} to each $\pi_\sigma$, with
$M_{3/2}(V^-)=\bigoplus M_{3/2}(\pi_\sigma)\otimes\sigma$ and $S_{1/2}((V^-)^*)=\bigoplus S_{1/2}(\pi_\sigma^*)\otimes\sigma^*$, gives $\chi_{\RR}=\chi_{M_{3/2}}-\overline{\chi_{S_{1/2}}}$.
The constant-term map $M_{3/2}(\rho_D)\to C$ is equivariant with kernel $S(D)$ and image $C\cap W^\perp$ (Theorem~\ref{thm:dim}(v)); as $C/(C\cap W^\perp)\cong W^\vee$, the image
has character $\chi_C-\overline{\chi_W}$, $\chi_W=\chi_M-\chi_{S_{1/2}}$. Hence $\chi_S=\chi_{M_{3/2}}-\chi_C+\overline{\chi_W}=\chi_{\RR}-\chi_C+\overline{\chi_M}$.
\end{proof}
\noindent So $S(D)\cong\bigoplus_\sigma S_{3/2}(\pi_\sigma)\otimes\sigma$ with $\dim S_{3/2}(\pi_\sigma)=\RR(\pi_\sigma)-d_0(\sigma)+m_{1/2}(\sigma)$, where $d_0(\sigma)$ and
$m_{1/2}(\sigma)$ are the multiplicities of $\sigma$ in $C$ and in the representation with character $\overline{\chi_M}$ (for real $\chi_\sigma$, in $M_{1/2}(\rho_D^*)$).

\noindent\textbf{The term $\chi_{\RR}$.} For an operator $A$ of finite order $n$ commuting with $O(D)$, $\sum_\sigma\alpha(A|_{\pi_\sigma})\chi_\sigma(g)=\sum_\beta\beta\,\tr(gP_\beta(A))$ with
$P_\beta(A)=n^{-1}\sum_{j<n}\e(-j\beta)A^j$. With $Y_1=\e(\frac38)\rho(S)$, $Y_2=(\e(\frac14)\rho(ST))^{-1}$ and $a(g)=\tr(g|V^-)$, traces on $V^-$,
\begin{multline}\label{eq:chiRR}
\chi_{\RR}(g)=\tfrac98a(g)-\tfrac14\bigl(a(g)-\tr(gY_1)\bigr)\\-\tfrac19\bigl[3a(g)+(\e(-\tfrac13)+2\e(\tfrac13))\tr(gY_2)+(\e(\tfrac13)+2\e(-\tfrac13))\tr(gY_2^2)\bigr]-\Psi(g),
\end{multline}
$\Psi(g)=\frac12\bigl(\sum_{gx=x}\{q(x)\}-\sum_{gx=-x}\{q(x)\}\bigr)$. All traces factor over the primes, as $g=\otimes g_p$, $\rho(S)=\otimes\rho_p(S)$ and $\iota=\otimes\iota_p$.
Likewise $\chi_C(g)=\frac12(\#\{x\in X_0:gx=x\}-\#\{x\in X_0:gx=-x\})$, from local fixed-point value distributions.

\noindent\textbf{The term $\chi_M$.}
\begin{cited}{Skoruppa [Sk08, \S5], from [Sk85, Thm.~5.2, p.~101] and [Sk85, Thm.~1.8, p.~22]}
$M_{1/2}(\Gamma(4m))\cong\bigoplus_{l\mid m,\ m/l\ \mathrm{squarefree}}\bigl(W(l)^\varepsilon\bigr)^*$ as $\Mp$-modules, where $W(l)$ is the Weil representation of
$D_l=(\Z/2l,x^2/4l)$, $\varepsilon:x\mapsto-x$, and $\Gamma(4m)^*$ is the principal congruence subgroup with the theta multiplier.
\end{cited}
\begin{cited}{[Sk08, \S3]; Scheithauer [Sch09, \S4, Prop.~4.3]}
$W(A\perp B)\cong W(A)\otimes W(B)$ and $W(-A)\cong W(A)^*$. If $\sig A$ is odd, the level $N_A$ of $A$ is divisible by $4$ and $\Gamma(N_A)^*$ acts
trivially on $W(A)$. If $A$ has even signature and its level divides $N'$, then $\Gamma(N')$ acts trivially on $W(A)$.
\end{cited}
\why{the vector-valued Serre--Stark theorem and basic properties of Weil representations, in the notation of this paper: $W(A)$ is the Weil
representation $\rho_A$ of \eqref{eq:weil}, and [Sk08] writes the condition on the signature as $\sigma(A)^4\ne1$ with $\sigma(A)=\e(-\sig A/8)$, the factor
in \eqref{eq:weil}; so it says that $\sig A$ is odd. We use this statement, attributed in [Sk08, \S3] to work in preparation, only for $\rho_D^*$ itself, and check it by replacing $m$ with $2m$, $3m$, $4m$ below; the normalization of [Sch09] is the dual of ours, with
the same kernel.}
Let $m=N/4$. As $\rho_D^*=W(D(-1))$, $\sig D(-1)=-1$ is odd and $D(-1)$ has level $N$, $\rho_D^*$ is trivial on $\Gamma(4m)^*$, and $F\mapsto(\lambda\mapsto\lambda\circ F)$ identifies $M_{1/2}(\rho_D^*)$, $O(D)$-equivariantly, with
$\Hom_{\Mp}(\rho_D,M_{1/2}(\Gamma(4m)))$. So $M_{1/2}(\rho_D^*)\cong\bigoplus_l\Inv\bigl(W(E^{(l)})\bigr)^{\varepsilon\times1}$ with $E^{(l)}=D_l(-1)\oplus D(-1)$, $O(D)$ acting on the
second summand. As $\sig E^{(l)}=-2$, every $E^{(l)}_p$ has even signature, $W(E^{(l)}_p)$ factors through $\SL_2(\Z/p^{k_p})$, and as $\SL_2(\Z)\to\prod_p\SL_2(\Z/p^{k_p})$ is
onto,
\begin{equation}\label{eq:chiM}
\chi_M(g)=\sum_{l\mid m,\ m/l\ \mathrm{squarefree}}\ \sum_{(e_p):\ \prod e_p=1}\ \prod_p\tr\bigl(g_p\,\big|\,\Inv W(E^{(l)}_p)^{\varepsilon_p=e_p}\bigr).
\end{equation}

\begin{lemma}\label{lem:inv}
Let $E$ have even signature, $I=q_E^{-1}(0)$, and $S_{II}$ the compression of $\rho_E(S)$ to $\C[I]$. Then $\Inv W(E)=\ker(2-S_{II}-S_{II}^*)$ on $\C[I]$.
\end{lemma}
\begin{proof}
$\rho(T)v=v$ iff $\operatorname{supp}v\subseteq I$. For $v\in\C[I]$, $\langle v,(2-S_{II}-S_{II}^*)v\rangle=\|\rho(S)v-v\|^2$, as $\rho(S)$ is unitary; so the operator is
$\ge0$ and its kernel is $\{v\in\C[I]:\rho(S)v=v\}$.
\end{proof}
\noindent\textbf{Computation} \gC. The $278$ candidates are the $258$ of \S\ref{subsec:small} and $20$ more, each with a component $(\Z/p)^3$, $p$ odd,
which the bound $\bar m_p=p+1$ of Lemma~\ref{lem:local} excludes and implementation B keeps as candidates. Conjugacy classes of $O(D)=\prod_pO(D_p)$ are products of local classes, the local groups
enumerated by Lemma~\ref{lem:deg}(b); the largest $|O(D)|$ is $11{,}520$. Kernels in Lemma~\ref{lem:inv} are taken with threshold $10^{-8}$ after checking
that no eigenvalue lies in $(10^{-8},10^{-4})$. Checks: $\chi_{\RR}(1)-\chi_C(1)$ equals \eqref{eq:rr}; replacing $m$ by $2m,3m,4m$ leaves $\chi_M$ unchanged on
$40$ random modules; $s$, $\rr$, $m_{1/2}$ and $\langle\chi_S,\chi_S\rangle$ deviate from integers by at most $5.7\cdot10^{-14}$, and $|\operatorname{Im}\chi_S|\le8.1\cdot10^{-15}$.
\emph{Result}: $s=0$ for $45$ modules, $s=1$ for $60$, and exactly $23$ are rigid, the same sets in both implementations. Among the $258$, twenty are
discriminant forms of lattices with $|D|>648$ (the largest of order $4608$); they satisfy $\rr\le\cp$, and the character excludes them. So the bound
$|\det L|\le648$ comes from this computation, not from the inequality $\rr\le\cp$. As a third check, the basis of \texttt{weilrep} [Wi] confirms $s$ and
the character on the $23$ rigid modules. On the $23$, $\chi_S$ is rational-valued and $d_0(\sigma_0)=m_{1/2}(\sigma_0)=0$.

\subsection{Lattices}\label{subsec:nik}
\begin{cited}{Nikulin [Ni79, Thm.~1.10.1]}
An even lattice with invariants $(t_+,t_-,q)$ exists if and only if \textup{(1)} $t_+-t_-\equiv\sig q\pmod 8$; \textup{(2)} $t_\pm\ge0$, $t_++t_-\ge l(A_q)$;
\textup{(3)} $(-1)^{t_-}|A_q|\equiv\operatorname{discr}K(q_p)\bmod(\Z_p^\times)^2$ for every odd $p$ with $t_++t_-=l(A_{q_p})$; \textup{(4)}
$|A_q|\equiv\pm\operatorname{discr}K(q_2)\bmod(\Z_2^\times)^2$ if $t_++t_-=l(A_{q_2})$ and $q_2\ne q_\theta^{(2)}(2)\oplus q_2'$. Here $K(q_p)$ is the $p$-adic lattice
of rank $l(A_{q_p})$ with discriminant form $q_p$.
\end{cited}
\why{the existence theorem for even lattices with given discriminant form; we use only the necessity of (3) and (4), and (3) is proved directly below.
(1) and (2) hold here.}

\begin{proposition}\label{prop:oddnik}
Let $L$ be even of signature $(2,1)$ with $D(L)\cong D$, $p$ odd, $l(D_p)=3$, $D_p=p^{(k_1,k_2,k_3)}[u_1,u_2,u_3]$. Then
$\bigl(\frac{-|D|/p^{k_1+k_2+k_3}}p\bigr)=\bigl(\frac{2u_1u_2u_3}p\bigr)$.
\end{proposition}
\begin{proof}
$L\otimes\Z_p\cong\perp_i\langle p^{a_i}\varepsilon_i\rangle$ with $\varepsilon_i\in\Z_p^\times$, and all $a_i\ge1$ as $D_p$ needs three generators. The discriminant form of
$\langle p^a\varepsilon\rangle$ is $(\Z/p^a,(2\varepsilon)^{-1}x^2/p^a)$, so by Wall's theorem $a_i=k_i$ and $\prod\varepsilon_i\equiv2\prod u_i$ modulo squares. As
$\det L=-|D|\equiv p^{\sum a_i}\prod\varepsilon_i$ modulo $(\Q_p^\times)^2$, the claim follows.
\end{proof}
\noindent\textbf{Computation} \gC. Conditions (3) and (4) split the $36{,}683$ modules into $32{,}980$ lattice modules and $3{,}703$ others; implementation B
finds the same split twice, with its own implementation of (3), (4), and by listing the rank-three $\Z_p$-lattices with discriminant form $D_p$ by their Jordan
decompositions. The $23$ rigid modules split $12/11$: six violate Proposition~\ref{prop:oddnik} at $p=3$, five violate only (4). For the existence of the
twelve, which does not need Nikulin's theorem: the models of Table~\ref{tab:twelve} are even, of signature $(2,1)$ and determinant $-|D|$, and their
discriminant forms have the value distributions and $2$-torsion values of the twelve modules; as fingerprints separate classes in this range, they are
isometric to them \gC.

\subsection{Class number one, theta series, conclusion}\label{subsec:classno}
\begin{cited}{Conway--Sloane [CS99, Ch.~15, Thm.~14] (Eichler) and [CS99, Ch.~15, \S7, p.~378]}
For indefinite forms of dimension at least $3$, a spinor genus contains exactly one class. An indefinite form of dimension at least $3$ with $|\det f|<128$ is
alone in its genus.
\end{cited}
\why{class numbers of indefinite ternary lattices.}
\begin{cited}{O'Meara [OM00, \S\S101--102]}
Let $V$ be a regular quadratic space over a global field $F$ with $\dim V\ge3$, and $T$ the set of real places at which $V$ is anisotropic. Then
$\theta(O^+(V))=P_T$, the group of $a\in F^\times$ positive at every place of $T$, modulo $(F^\times)^2$; and the number of spinor genera in
$\operatorname{gen}L$ is $[J_F:P_TJ_F^L]$, $J_F^L=\{i:i_p\in\theta(O^+(L_p))\ \forall p\}$.
\end{cited}
\why{strong approximation for indefinite lattices. For $F=\Q$ and $V$ indefinite, $T=\emptyset$, also when $L\otimes\Q$ is anisotropic, as for $B_{184}$ and $B_{324}$.
As $\dim V=3$ is odd, $-1$ is improper and acts trivially on lattices, so proper and improper spinor genera coincide; and as $\Q$ has class number one,
the genus is one spinor genus if every $u\in\prod_p\Z_p^\times$ lies in $\Q^\times J^L$.}
\emph{$|\det L|<128$.} $L_{16},B_{40},L_{54},L_{64},L_{72},B_{88},L_{96},B_{120}$ are alone in their genera.

\emph{$L_{250}$.} Let $M=U\oplus\langle2\rangle$, so $h(M)=1$ and $L_{250}=M(5)$. If $L\in\operatorname{gen}M(5)$, then $L\otimes\Z_5\cong M(5)\otimes\Z_5$ has all products in
$5\Z_5$, and $5$ is a unit at $\ell\ne5$; so $L(\frac15)$ is even, integral and in $\operatorname{gen}M$. Thus $L\mapsto L(\frac15)$ is a bijection
$\operatorname{gen}M(5)\to\operatorname{gen}M$ on isometry classes, and $h(L_{250})=1$.

\emph{Preparation.} If $p\nmid2\det L$, $L_p$ is unimodular of rank $3$; a unimodular plane in it represents every unit square class, and the reflections in
such vectors are integral, so $\theta(O^+(L_p))\supseteq\Z_p^\times$. If $v\in L_p$ has $2(x,v)/(v,v)\in\Z_p$ for all $x$, the reflection $\tau_v$ lies in $O(L_p)$, and
$\theta(\tau_v\tau_w)=Q(v)Q(w)$ modulo squares, $Q(v)=(v,v)/2$.

\emph{$B_{324}=A_2\oplus\langle-108\rangle$, $B_{648}=A_2\oplus\langle-216\rangle$.} At $2$: $Q(x,y)=x^2-xy+y^2$ takes the unit values $1,3,13,7$ at
$(1,0),(1,2),(1,4),(1,3)$, and $2(x,v)/(v,v)=(x,v)/Q(v)\in\Z_2$; so $\theta_2\supseteq\Z_2^\times$. At $3$: $e_1=(1,0)$, $e_2=(1,2)$ and the generator $e_3$ of
$\langle c\rangle$ form an orthogonal $\Z_3$-basis with $Q=1,3,c/2$, and reflections in an orthogonal basis are integral, so $\theta_3\ni Q(e_i)Q(e_j)$. For $c=-216$,
$Q(e_2)Q(e_3)=-324=-18^2$ is the non-square unit class, so $\theta_3\supseteq\Z_3^\times$. For $c=-108$ the products are $3$, $-54\sim3$, $-162\sim1$, so
$\theta_3\supseteq\{1,3\}(\Q_3^\times)^2$; if $u_3$ is not a square, $-u\in J^L$, so $u\in\Q^\times J^L$. Both genera are one spinor genus, hence one class.

\emph{$B_{184}=\langle2\rangle\oplus\langle2\rangle\oplus\langle-46\rangle$.} At $23$, $\langle2\rangle\oplus\langle2\rangle$ is a unimodular plane, so $\theta_{23}\supseteq\Z_{23}^\times$. At
$2$, $-23\equiv1\pmod8$ gives $L_2\cong\langle2\rangle^3$; the vectors $(1,0,0)$, $(1,1,1)$, $(1,2,0)$ have $Q=1,3,5$, the reflections are integral, and
$\theta_2\ni3,5,15$: $\theta_2\supseteq\Z_2^\times$. So the genus is one spinor genus. A Sage computation confirms class number one for all twelve \gC.

\noindent\textbf{Theta series.}
\begin{cited}{Bruinier--Ehlen--Freitag [BEF16, Prop.~4.18]}
Let $U$ be a totally isotropic subgroup of a finite quadratic module $A$. The map $\C[U^\perp/U]\to\C[A]$, $e_{\gamma+U}\mapsto\sum_{u\in U}e_{\gamma+u}$, maps
$M_k(\rho_{U^\perp/U})$ into $M_k(\rho_A)$ and cusp forms to cusp forms.
\end{cited}
\why{lifting through isotropic subgroups. [BEF16] use the dual normalization (\S\ref{subsec:weil}); applied to $A(-1)$ it gives the statement in ours.}
Let $U\subseteq D(L)$ be totally isotropic with $U^\perp/U\cong(\Z/2m,x^2/4m)$, the discriminant form of $\langle2m\rangle$. The odd unary theta series
$\vartheta_m=(\vartheta_{m,r})_{r\in\Z/2m}$, $\vartheta_{m,r}=\sum_{\ell\equiv r\,(2m)}\ell\,q^{\ell^2/4m}$, lies in $S_{3/2}$ for $(\Z/2m,x^2/4m)$, and its lifts through
$U$, for the isomorphisms $U^\perp/U\cong(\Z/2m,x^2/4m)$, lie in $S(D(L))$. Let $\Theta(L)\subseteq S(D(L))$ be the span of these lifts; it is $O(D)$-stable, and $\vartheta_1=0$. If $D_2$ is anisotropic of order $8$, then $\Theta(L)=0$: the
denominators of $q$ on $D_2$ and on $D_{\rm odd}$ are coprime, so $q(\gamma)=0$ forces $\gamma\in D_{\rm odd}$, every isotropic $U$ lies in $D_{\rm odd}$, and
$U^\perp/U\supseteq D_2$ is not cyclic. This applies to $B_{40},L_{72},B_{88},B_{120},B_{184}$. For the other seven we list the isotropic $U\ne0$, lift
$\vartheta_m$ through every isomorphism $U^\perp/U\cong(\Z/2m,x^2/4m)$, and compute the exact rank of the lifts on the coefficients of exponent $\le30$, a lower
bound for $\dim\Theta(L)$: it is $s=2,3,2,5$ for $L_{16},L_{54},L_{64},L_{250}$ (with $m=2,3,2,5$); $L_{96}$ and $B_{324}$ have one isotropic $U\ne0$ and
no quotient of this form, and of the seven isotropic $U\ne0$ of $D(B_{648})$ only two have $U^\perp/U$ of this form, both with $m=1$; so $\Theta(L)=0$ for these three \gC.

\noindent\textbf{Conclusion.} By \S\ref{subsec:large} and \S\ref{subsec:small} every rigid module satisfying \eqref{eq:H} is among the $278$ candidates, and by
\S\ref{subsec:chi} exactly $23$ of them are rigid. A Petersson-rigid $L$ has $D(L)$ among these, hence among the twelve lattice modules of \S\ref{subsec:nik};
by Nikulin's corollary (\S\ref{subsec:lat}) $L$ lies in the corresponding genus, which consists of one class. With the theta series, this proves
Theorem~A. \qed

\subsection{Root bases}\label{subsec:roots}
A \emph{root} of $L$ is a primitive $r\in L$ with $(r,r)>0$ and $2(r,x)/(r,r)\in\Z$ for all $x\in L$; a \emph{root basis} is a $\Z$-basis $r_1,r_2,r_3$ of roots with
$(r_i,r_j)\le0$ for $i\ne j$. Its Cartan matrix $A=(a_{ij})$, $a_{ij}=2(r_i,r_j)/(r_i,r_i)$, is a symmetrizable generalized Cartan matrix: with $d_i=(r_i,r_i)/2$, the
Gram matrix is $G=\operatorname{diag}(d_1,d_2,d_3)A$. $A$ is \emph{connected} if its Dynkin diagram is; a connected $A$ whose symmetrization has
signature $(2,1)$ is \emph{hyperbolic} if every proper connected subdiagram is of finite or affine type [Ka90, Ch.~4], and \emph{compact} if they are all
of finite type. In rank three these subdiagrams are the vertices and the edges, and the edge $\{i,j\}$ is of finite type iff $a_{ij}a_{ji}\le3$ and of affine
type iff $a_{ij}a_{ji}=4$ [Ka90, Ch.~4]; so, for $A$ whose symmetrization has signature $(2,1)$,
\begin{equation}\label{eq:hyp}
A\ \text{hyperbolic}\iff A\ \text{connected and}\ a_{ij}a_{ji}\le4\ (i\ne j);\qquad A\ \text{compact}\iff\text{moreover}\ a_{ij}a_{ji}\le3.
\end{equation}
For connected $A$ the symmetrizing diagonal is unique up to a positive scalar; let $d^0$ be the primitive integral one and $M_A$ the lattice with Gram matrix
$\operatorname{diag}(d^0)A$, the \emph{minimal symmetrization}. The census of [II, \S3.1] lists $100$ Cartan matrices with the Gram matrices of their $M_A$: the
\emph{core}, the $44$ classes (under simultaneous permutation) of rank-three symmetrizable hyperbolic Cartan matrices whose symmetrization has signature
$(2,1)$, $11$ compact and $33$ non-compact, and the \emph{rim}, $56$ matrices one edit away. We re-enumerated the core from \eqref{eq:hyp}: $44$ classes, $11$
compact, equal to the core rows of [II] \gC.

\begin{lemma}\label{lem:rootfin}
\textup{(a)} For $k\ge1$ the root bases of $L(k)$ are those of $L$, with the same Cartan matrices. \textup{(b)} If $L$ has a root basis with connected Cartan
matrix $A$, then $L\cong M_A(k)$, $k=\gcd(d_1,d_2,d_3)$. \textup{(c)} If $G=\bigl(\begin{smallmatrix}a&x&y\\x&b&z\\y&z&c\end{smallmatrix}\bigr)$ is the Gram matrix of a root basis
of $L$, then $a,b,c\in\{2n:n\mid\exp D(L)\}$; $x,y,z\le0$ with $cx^2,\,by^2,\,az^2\le abc+|\det L|$; and $\operatorname{lcm}(\frac a2,\frac b2)\mid x$,
$\operatorname{lcm}(\frac a2,\frac c2)\mid y$, $\operatorname{lcm}(\frac b2,\frac c2)\mid z$.
\end{lemma}
\begin{proof}
(a) Primitivity does not see the form, and $2k(r,x)/k(r,r)=2(r,x)/(r,r)$. (b) $(d_i)$ symmetrizes $A$, so $(d_i)=k\,d^0$. (c) For a root $r$ with $(r,r)=2n$,
$r/n\in L'$ has order exactly $n$ in $D(L)$, by primitivity. $\det G=abc+2xyz-az^2-by^2-cx^2$ with $2xyz\le0$ gives $-|\det L|\le abc-cx^2$, and likewise.
Finally $a_{12}=2x/a$ and $a_{21}=2x/b$ are integers.
\end{proof}
\noindent If the genus of $L$ has class number one, the Gram matrices of the root bases of $L$ are exactly the integral symmetric $G$ satisfying (c) with
$\det G=\det L$ and $D(G)\cong D(L)$: such a $G$ has signature $(2,1)$, lies in the genus of $L$ by Nikulin's corollary, hence is a Gram matrix of $L$, and its
basis vectors are roots by (c). So there are finitely many, and they can be listed.

\begin{proposition}\label{prop:roots}
\textup{(a)} An even lattice of signature $(2,1)$ has a hyperbolic root basis if and only if it is $M(k)$ for a core row $M$ of the census and some $k\ge1$.
\textup{(b)} The twelve lattices of Theorem~A have $39$ root bases up to isometry and permutation, $32$ of them connected \gC. Exactly
$L_{16},L_{54},L_{64},L_{72},L_{96},L_{250}$ have a hyperbolic one; they are, in this order, $FF(2)$, $FF(3)$, $M_8(2)$, $L_{72}(1)$, $M_{12}(2)$,
$FF(5)$ (Table~\ref{tab:twelve}), and their hyperbolic root bases are compact exactly for $L_{72}$ and $L_{96}$. $B_{40}$ and $B_{88}$ are rim rows of the census, and $B_{120},B_{184},B_{324},B_{648}$ are not of
the form $M(k)$ for any census row $M$ \gC. The $39$ root bases, with their labels (the chain names of [II, \S3.2]), are listed in the file
\texttt{ROOT\_BASES.md} of the reproduction archive. \textup{(c)} Among the pairs $(M,k)$, $M$ a census row and $k^3|\det M|\le648$, exactly twelve give a
Petersson-rigid $M(k)$; for the Feingold--Frenkel row $FF$, $FF(k)$ is Petersson-rigid for $k=2,3,5$ and not for $k=1,4,6$ \gC.
\end{proposition}
\begin{proof}
(a) A hyperbolic Cartan matrix is, by \eqref{eq:hyp} and the enumeration of the core, one of the $44$ core matrices; by Lemma~\ref{lem:rootfin}(b) the lattice
is $M_A(k)$. Conversely, by Lemma~\ref{lem:rootfin}(a) the root basis of $M$ is one of $M(k)$ with the same Cartan matrix. (b) All twelve genera have class
number one (Theorem~A), so the root bases are listed as above; each listed Gram matrix $G$ is moreover certified by an explicit $U\in\GL_3(\Z)$ with
$U^{\mathsf T}G_0U=G$, $G_0$ the model of Table~\ref{tab:twelve}. Read off $A$ and $k$ and apply (a), Lemma~\ref{lem:rootfin}(b) and \eqref{eq:hyp}. For the
last four, no Cartan matrix of their root bases occurs in the census; if one of them were $M_A(k)$ for a census row with matrix $A$, it would have a root
basis with Cartan matrix $A$, by Lemma~\ref{lem:rootfin}(a). (c) By Theorem~A, $M(k)$ is Petersson-rigid if and only if $D(M(k))$ is isomorphic to
one of the twelve lattice modules, and $|\det|\le648$ for these; $158$ pairs are tested.
\end{proof}

\subsection{The dimensions at weight \texorpdfstring{$\frac52$}{5/2}}\label{subsec:wt52}
This subsection is not part of the proof of Theorem~A; with the formula \eqref{eq:RRgen} it explains the last two columns of
Table~\ref{tab:canon}.

\noindent\textbf{The columns of weight $\frac52$} \gC. By Lemma~\ref{lem:parity}, $M_{5/2}(\rho_D)$ takes values in $V^+$, as $i^{1-5}=1$, and
$S_{-1/2}(\rho_D^*)=0$.
\begin{cited}{Bruinier [Br02, Thm.~1.6]}
Let $L$ be an even lattice of signature $(2,l)$, $l\ge3$, and $\kappa=1+l/2$. For $\gamma\in L'/L$ with $q(\gamma)=0$ the Eisenstein series
$E_\gamma=\frac12\sum_{(M,\varphi)}e_\gamma|^*_\kappa(M,\varphi)$, summed over $\langle T\rangle\backslash\Mp$, converges and is a holomorphic modular form of weight $\kappa$ for
the dual Weil representation $\rho_L^*$, with constant term $e_\gamma+e_{-\gamma}$.
\end{cited}
\why{the Eisenstein series of the Weil representation. Take $l=3$ and $L$ with $D(L)\cong D(-1)$; it exists by Nikulin's theorem (\S\ref{subsec:nik}),
as $3-2\equiv\sig D\pmod8$ and $5>l(D_p)$, so that conditions (3) and (4) are void. Then $\rho_L^*=\rho_D$ and $\kappa=\frac52$.}
So the constant-term map from $M_{5/2}(\rho_D)$ onto the even functions on $X_0$ is surjective with kernel $S_{5/2}(\rho_D)$, and \eqref{eq:RRgen} for
$\rho_D|_{V^+}$, $k=\frac52$, minus the number of $\{\pm1\}$-orbits in $X_0$, is $\dim S_{5/2}(\rho_D)$; this is [BEF16, (3.4)] for $A=D(-1)$. The
invariants $W=\C[D]^{O(D)}$ form a subrepresentation in $V^+$, $S_{5/2}(\rho_D)^{O(D)}=S_{5/2}(\rho_D|_W)$, and the averaged series $\sum_{x\in\Omega}E_x$ over the
$O(D)$-orbits $\Omega\subseteq X_0$ have the orbit sums as constant terms; so the same formula on $W$, minus the number of $O(D)$-orbits in $X_0$, is
$\dim S_{5/2}(\rho_D)^{O(D)}$. The traces are computed prime by prime; the formula reproduces the $98$ values of $\dim S_{5/2}$ computed with \texttt{weilrep}
in [II].
\section{Open problems}\label{sec:conj}

\begin{enumerate}[label=\textbf{Problem \arabic*.},leftmargin=*,itemsep=3pt,ref=\arabic*]
\item\label{prob:hecke} Prove that the eigensystem of $S(D)$ is that of $E$ for each of the eight lattices. One route is a Shimura correspondence for
vector-valued forms on these discriminant forms, with control of its level, and a Sturm bound. For $B_{40},L_{72},L_{96}$ the form $F_0$ is explicit
(Proposition~\ref{prop:etatheta}), and the curves are eta quotients: $20a1$ is $\eta(2\tau)^2\eta(10\tau)^2$, $36a1$ is $\eta(6\tau)^4$, $24a1$ is
$\eta(2\tau)\eta(4\tau)\eta(6\tau)\eta(12\tau)$ [MO97, Thms.~1, 2].
\item\label{prob:dict} Decide whether the identity of Remark~\ref{rem:dims} holds for every odd prime $p$ and every $k\in\frac32+\Z$, $k\ge\frac32$, and
whether the eigensystems of $S_{3/2}(\rho_{\Dp_p})$ are exactly those of the newforms of weight $2$ and level $p$, $2p$ or $4p$ with $w_p=+1$. If they
are, the number of Lemma~\ref{lem:three} is the number of distinct eigensystems of $S_{3/2}(\rho_{\Dp_p})$ (\S\ref{subsec:newforms}(1)).
\item\label{prob:rat} Prove a criterion of Bruinier--Ono type for the rationality of coefficient vectors of the forms $\Hh$, and explain the
cases it cannot cover. Among $1608$ coefficient vectors tested across the eight forms, $91$ agree numerically with vectors of rationals (one
implementation; the rule is stated below). Of these, $47$ occur where a quadratic twist of the curve $E$ of Table~\ref{tab:canon} has analytic
rank $3$ or where a local factor vanishes. The remaining $44$, all from $B_{324}$, include coefficients that agree with rationals although the
derivative at $s=1$ of the corresponding twisted $L$-function does not vanish. The rule, the data, and what is proved about them follow.

The coordinates $h_{r,i}$ of $\Hh$ do not depend on the normalization of $F_0$ (\S\ref{subsec:coeffs}), so whether a
coefficient $a_{r,i}(n)$ is rational is a property of $\Hh$. For the eight lattices write $2|D|n=|\Delta_0|f^2$ with $\Delta_0<0$ a fundamental
discriminant, and $E_{\Delta_0}$ for the twist of the curve $E$ of Table~\ref{tab:canon} by $\chi_{\Delta_0}$. We say that the vector $(a_{r,i}(n))_i$
\emph{agrees with a vector of rationals} if for every $i$ there are integers $u,v$ with $1\le v\le60$, $|u|\le10^5$ and $|a-u/v|\le50\,e$, where
$a$ is the value of $a_{r,i}(n)$ computed with $(M,Q,Y)=(40,60,0.20)$ (\S\ref{subsec:coeffs}), $a'$ the value computed with $(36,56,0.22)$,
and $e=|a-a'|+10^{-13}|a|$. Among the first sixteen exponents $n$ of every class $r$ (twelve for $B_{324},B_{648}$), $1608$ pairs $(r,n)$ in all,
exactly $91$ agree, where numbers unrelated to fractions would give $0.34$ agreements on average \gC. The count depends on the rule: with $e$
the spread of the three choices of \S\ref{subsec:coeffs}, $8$ further pairs, of $B_{88}$ and $B_{120}$, agree, none of them of the kinds
below \gC. At $f=1$
the $91$ pairs are exactly the $14$ pairs at which $E_{\Delta_0}$ has analytic rank $3$, i.e.\ $L'(E,\chi_{\Delta_0},1)=0$; every twist met at $f=1$ has root
number $-1$ and analytic rank $1$ or $3$. At $f\ge2$ they are $5$ further pairs of analytic rank $3$; $21$ pairs, of $B_{88},B_{184},B_{648}$, at which a
factor $a_p(E)-\chi_{\Delta_0}(p)$ with $p\nmid2|D|$, $p\,\|\,f$ vanishes; the $7$ pairs of $B_{120}$, at which $3\,\|\,f$ and $\chi_{\Delta_0}(3)=+1$,
so that the factor $1-\chi_{\Delta_0}(3)$ vanishes; and $44$ pairs of $B_{324}$
with $\Delta_0\equiv1\pmod8$ \gC. These counts were made by one implementation only.

At the $44$ pairs $f$ is even and $\Delta_0$ is odd, so that $\Delta_0\equiv1\pmod8$ means $\chi_{\Delta_0}(2)=+1$; as $E=162d1$ has split
multiplicative reduction at $2$, $a_2(E)=+1$, these are the pairs at which $\chi_{\Delta_0}(2)=a_2(E)$. They are also the pairs at which
$E_{\Delta_0}$ has root number $+1$: for every fundamental discriminant $\Delta<0$ with $3\mid\Delta$ and $|\Delta|<4000$, $E_\Delta$ has root
number $+1$ exactly when $\Delta\equiv1\pmod8$, and the data do not separate the two descriptions \gC. A root number $+1$ does not make
$L'(E,\chi_{\Delta_0},1)$ vanish. At the $28$ discriminants of the $44$ pairs, $E_{\Delta_0}$ has additive reduction at $3$, and
$L(E,\chi_{\Delta_0},s)=\sum_n a_n(E)\chi_{\Delta_0}(n)n^{-s}$ is the $L$-function of $E_{\Delta_0}$ \gC. Let $N_{\Delta_0}$ be the conductor of
$E_{\Delta_0}$ and $\Lambda(s)=(\sqrt{N_{\Delta_0}}/2\pi)^s\,\Gamma(s)\,L(E,\chi_{\Delta_0},s)$. As the root number is $+1$, the functional equation reads
$\Lambda(s)=\Lambda(2-s)$, so that $\Lambda'(1)=0$; differentiating the product at $s=1$ and using $\Gamma(1)=1$, $\Gamma'(1)=-\gamma_{\rm E}$, $\gamma_{\rm E}$
Euler's constant, gives $0=(\sqrt{N_{\Delta_0}}/2\pi)\bigl((\log\frac{\sqrt{N_{\Delta_0}}}{2\pi}-\gamma_{\rm E})L(E,\chi_{\Delta_0},1)+L'(E,\chi_{\Delta_0},1)\bigr)$, that is,
\[
L'(E,\chi_{\Delta_0},1)=-\Bigl(\log\frac{\sqrt{N_{\Delta_0}}}{2\pi}-\gamma_{\rm E}\Bigr)\,L(E,\chi_{\Delta_0},1).
\]
The first factor vanishes only if $N_{\Delta_0}=(2\pi e^{\gamma_{\rm E}})^2=125.23\ldots$, which is not an integer; so
$L'(E,\chi_{\Delta_0},1)\ne0$ whenever $L(E,\chi_{\Delta_0},1)\ne0$. This is the case at $23$ of the $28$ discriminants; at the other five
$L(E,\chi_{\Delta_0},1)$ vanishes to the working precision. At $n=\frac{26}{27}$ ($\Delta_0=-39$, $f=4$), for instance, the coefficient is $-\frac54$,
while $L(E,\chi_{-39},1)=7.31446$ and $L'(E,\chi_{-39},1)=-19.70263$ \gC.

For the other $47$ pairs the model is [BO10, Thm.~1.1(2)]. It assumes that $p$ is
prime and that $G\in S_2^{\rm new}(\Gamma_0(p))$ is a newform with root number $\epsilon(G)=-1$; for a fundamental discriminant $\Delta>0$ with
$(\frac\Delta p)=1$ it says that the coefficient $c^+_g(\Delta)$ of the harmonic weak Maass form $f_g$ of [BO10], $g$ the form of weight $\frac32$
corresponding to $G$, is algebraic exactly when $L'(G,\chi_\Delta,1)=0$; and
the proof runs through twisted Heegner divisors on $X_0(p)$ and the Gross--Zagier formula. Here the conductors of the curves $E$ of
Table~\ref{tab:canon}, $20,36,44,24,20,92,162,324$, are not prime; the discriminants $\Delta_0$ are negative, of the sign opposite to that of
$\Delta$ in [BO10, Thm.~1.1(2)]; seven of the eight lattices have compact Shimura curves; and the pattern at $f\ge2$ has no counterpart there. The $44$ pairs lie outside a statement of this type, since at $23$ of their
$28$ discriminants $L'(E,\chi_{\Delta_0},1)\ne0$.
\end{enumerate}

\appendix
\section{The eight harmonic weak Maass forms}\label{app:tabH}
Each of the eight forms is
\begin{align*}
\Hh_L&=\sum_{r,i}h_{r,i}\,\overline{u_{r,i}},\\
h_{r,i}&=\delta_{(r,i),(r_0,1)}\,q^{-\mu}+\sum_{n>0,\ n\equiv-r}a_{r,i}(n)\,q^n
-\frac1{t_\mu}\sum_{m>0}\overline{c_{f_{r,i}}(m)}\,(4\pi m)^{-1/2}\,\Gamma(\tfrac12,4\pi my)\,q^{-m}
\end{align*}
(\S\ref{subsec:coeffs}). It is determined by $\mu$, by $F_0=\sum f_{r,i}u_{r,i}$ and by $t_\mu$ (Theorem~\ref{thm:harm}); the numbers
$a_{r,i}(n)$, which the theorem does not give in closed form, write it out. The coordinate functions $f_{r,i}$ are
$\eta^2$ times theta functions for $B_{40},L_{72},L_{96}$ (Proposition~\ref{prop:etatheta}) and are computed from the basis of [Wi] for the other five.
Table~\ref{tab:H} lists, for each lattice, $t_\mu$ and the first four
coefficients $a_{r,i}(n)$: one row per basis vector $u_{r,i}=c_{F_0}(e_{r,i})$ (\S\ref{subsec:mult}), with its class $r$, its exponent $e_{r,i}$, and the first
exponent $n_0\equiv-r$ of $h^+_{r,i}$; the row of $\mu$ carries the pole $q^{-\mu}$. $t_\mu$ is given to $16$ digits, from the exact form of [Wi] (\S\ref{subsec:table}). The other entries are given to ten significant digits, trailing zeros
omitted, with the accuracy stated in \S\ref{subsec:coeffs}; for $B_{40}$ the basis is the $u_\nu$ of \S\ref{subsec:b40} \gC.
\begin{center}\footnotesize\setlength{\tabcolsep}{3.2pt}
\begin{longtable}{lllllll}
\caption{The forms $\Hh$: $t_\mu$ and $h^+_{r,i}=[(r,i)=(r_0,1)]\,q^{-\mu}+\sum_na_{r,i}(n)q^n$.}\label{tab:H}\\
\toprule
$r$ & $e_{r,i}$ & $n_0$ & $a_{r,i}(n_0)$ & $a_{r,i}(n_0+1)$ & $a_{r,i}(n_0+2)$ & $a_{r,i}(n_0+3)$\\\midrule
\endfirsthead
\toprule
$r$ & $e_{r,i}$ & $n_0$ & $a_{r,i}(n_0)$ & $a_{r,i}(n_0+1)$ & $a_{r,i}(n_0+2)$ & $a_{r,i}(n_0+3)$\\\midrule
\endhead
\multicolumn{7}{l}{\rule{0pt}{3.2ex}$B_{40}$: $|D|=40$, $\mu=1/10$, $\dim\pi_{\sigma_0}=4$, $t_\mu=0.2697776010367185$, $E=20a1$}\\
$1/10$ & $1/10$ & $9/10$ & $q^{-1/10}+{}$ $-19.63638957$ & $-73.31653899$ & $-216.9125062$ & $-549.9328637$\\
$3/20$ & $3/20$ & $17/20$ &  $-8.125949325$ & $-46.11481602$ & $-123.7119644$ & $-336.299995$\\
$7/20$ & $7/20$ & $13/20$ &  $-12.44872546$ & $-56.46429951$ & $-163.3671504$ & $-448.8898221$\\
$9/10$ & $9/10$ & $1/10$ &  $-2.363610426$ & $-14.68581861$ & $-58.95925057$ & $-164.4086628$\\
\multicolumn{7}{l}{\rule{0pt}{3.2ex}$L_{72}$: $|D|=72$, $\mu=1/12$, $\dim\pi_{\sigma_0}=4$, $t_\mu=0.3865317556546428$, $E=36a1$}\\
$1/12$ & $1/12$ & $11/12$ & $q^{-1/12}+{}$ $-9.630077549$ & $-30.36395081$ & $-84.55539928$ & $-190.9266434$\\
$1/6$ & $1/6$ & $5/6$ &  $-7.381529822$ & $-31.30892897$ & $-73.45207234$ & $-183.6644197$\\
$5/12$ & $5/12$ & $7/12$ &  $-5.860861162$ & $-21.70826588$ & $-60.84434929$ & $-147.3139886$\\
$5/6$ & $5/6$ & $1/6$ &  $-0.9507704938$ & $-7.147744202$ & $-19.82974721$ & $-53.31176766$\\
\multicolumn{7}{l}{\rule{0pt}{3.2ex}$B_{88}$: $|D|=88$, $\mu=1/22$, $\dim\pi_{\sigma_0}=10$, $t_\mu=1.074741854937327$, $E=44a1$}\\
$1/22$ & $1/22$ & $21/22$ & $q^{-1/22}+{}$ $-1.789754489$ & $-5.519098082$ & $-9.743653398$ & $-16.32578358$\\
$3/22$ & $3/22$ & $19/22$ &  $-1.118023406$ & $-1.777394543$ & $-5$ & $-8.291244835$\\
$7/44$ & $7/44$ & $37/44$ &  $-1.746980832$ & $-4.407910683$ & $-7.254079333$ & $-15.5110784$\\
$5/22$ & $5/22$ & $17/22$ &  $-2.887782837$ & $-5.450113965$ & $-12.38094114$ & $-21.94660691$\\
$9/22$ & $9/22$ & $13/22$ &  $-1.312472292$ & $-3.433298728$ & $-7.089918378$ & $-12.81307631$\\
$19/44$ & $19/44$ & $25/44$ &  $-0.8630430822$ & $-3.710210247$ & $-6.013733553$ & $-11.54674645$\\
$15/22$ & $15/22$ & $7/22$ &  $-0.7423668293$ & $-1.558058788$ & $-3.590365586$ & $-6.506999944$\\
$35/44$ & $35/44$ & $9/44$ &  $-0.5525359019$ & $-1.308488265$ & $-3.170719794$ & $-4.765233702$\\
$39/44$ & $39/44$ & $5/44$ &  $0.09013222184$ & $-1.001941206$ & $-0.9810332988$ & $-3.742932123$\\
$43/44$ & $43/44$ & $1/44$ &  $-0.1644019264$ & $-1$ & $-1.517133066$ & $-3.822566445$\\
\multicolumn{7}{l}{\rule{0pt}{3.2ex}$L_{96}$: $|D|=96$, $\mu=5/48$, $\dim\pi_{\sigma_0}=3$, $t_\mu=0.2323507202689307$, $E=24a1$}\\
$5/48$ & $5/48$ & $43/48$ & $q^{-5/48}+{}$ $-15.28666928$ & $-72.18640022$ & $-198.3108753$ & $-551.0183413$\\
$1/6$ & $1/6$ & $5/6$ &  $-21.66883673$ & $-88.66376426$ & $-276.2177992$ & $-715.8810815$\\
$29/48$ & $29/48$ & $19/48$ &  $-5.968026494$ & $-37.47375173$ & $-119.5715356$ & $-340.5516631$\\
\multicolumn{7}{l}{\rule{0pt}{3.2ex}$B_{120}$: $|D|=120$, $\mu=1/20$, $\dim\pi_{\sigma_0}=8$, $t_\mu=0.9345370234792852$, $E=20a1$}\\
$1/20$ & $1/20$ & $19/20$ & $q^{-1/20}+{}$ $-1.855982182$ & $-6.183758181$ & $-11.3294629$ & $-18.83085243$\\
$7/60$ & $7/60$ & $53/60$ &  $-3.683575181$ & $-8.043769421$ & $-16.93079338$ & $-32.25275609$\\
$3/10$ & $3/10$ & $7/10$ &  $-0.7601873584$ & $-2.356643796$ & $-4.106064929$ & $-9.983997292$\\
$11/30$ & $11/30$ & $19/30$ &  $-1.170865253$ & $-3.78895826$ & $-6.887501338$ & $-14.1608502$\\
$9/20$ & $9/20$ & $11/20$ &  $-2.232149753$ & $-6.882664604$ & $-14$ & $-25.57966463$\\
$7/10$ & $7/10$ & $3/10$ &  $-1$ & $-2.266784066$ & $-6.316199027$ & $-10.21903023$\\
$43/60$ & $43/60$ & $17/60$ &  $-1.062974663$ & $-3.149997499$ & $-7.121519269$ & $-13.37419586$\\
$29/30$ & $29/30$ & $1/30$ &  $-0.1590973936$ & $-1.647834296$ & $-2.319520794$ & $-6.682166966$\\
\multicolumn{7}{l}{\rule{0pt}{3.2ex}$B_{184}$: $|D|=184$, $\mu=7/92$, $\dim\pi_{\sigma_0}=22$, $t_\mu=0.7680627886051352$, $E=92a1$}\\
$7/92$ & $7/92$ & $85/92$ & $q^{-7/92}+{}$ $-5.497211142$ & $-14.17680014$ & $-38.01829011$ & $-82.52958502$\\
$5/46$ & $5/46$ & $41/46$ &  $0.1091844223$ & $-3.207312209$ & $-4.496725451$ & $-16.25439784$\\
$11/92$ & $11/92$ & $81/92$ &  $-4.516217952$ & $-17.82831054$ & $-43.26620731$ & $-95.24679511$\\
$7/46$ & $7/46$ & $39/46$ &  $-4.54987544$ & $-18.73664416$ & $-41.26869887$ & $-98.28259467$\\
$15/92$ & $15/92$ & $77/92$ &  $-4.7832628$ & $-12.20831361$ & $-30$ & $-68.31915493$\\
$19/92$ & $111/92$ & $73/92$ &  $-4.559191931$ & $-15.13336393$ & $-35.73237762$ & $-80.95902485$\\
$11/46$ & $11/46$ & $35/46$ &  $-1.97705521$ & $-7.681086217$ & $-21.8276568$ & $-47.10756132$\\
$15/46$ & $15/46$ & $31/46$ &  $-3.348677561$ & $-12.98020795$ & $-26.89382561$ & $-71.38461538$\\
$17/46$ & $17/46$ & $29/46$ &  $-3.284951981$ & $-13.84306606$ & $-34.45847806$ & $-81.06531898$\\
$19/46$ & $65/46$ & $27/46$ &  $-3.594889899$ & $-8.807968316$ & $-25.36175953$ & $-56.74798919$\\
$21/46$ & $21/46$ & $25/46$ &  $-3.80865173$ & $-10.15144603$ & $-30$ & $-64.92990807$\\
$43/92$ & $135/92$ & $49/92$ &  $-2.580302186$ & $-7.640926795$ & $-19.83676539$ & $-47.6350118$\\
$51/92$ & $51/92$ & $41/92$ &  $-2.850546109$ & $-10.4867168$ & $-28.58053846$ & $-66.16925994$\\
$63/92$ & $155/92$ & $29/92$ &  $0.3990274054$ & $-2.668548948$ & $-8.313733372$ & $-17.34293107$\\
$33/46$ & $33/46$ & $13/46$ &  $-0.2093386548$ & $-7.060717078$ & $-14.93580926$ & $-41.30352191$\\
$67/92$ & $67/92$ & $25/92$ &  $-0.8708076866$ & $-3.450039317$ & $-10.63933952$ & $-24.83059552$\\
$37/46$ & $37/46$ & $9/46$ &  $-2$ & $-8.913122188$ & $-20.89641772$ & $-55.82638139$\\
$79/92$ & $355/92$ & $13/92$ &  $-0.1750589759$ & $-0.7271073398$ & $-3.980951481$ & $-7.675192296$\\
$83/92$ & $83/92$ & $9/92$ &  $-2.097307711$ & $-6.760709046$ & $-21.62918424$ & $-53.57859981$\\
$43/46$ & $135/46$ & $3/46$ &  $-0.784669697$ & $-0.7328251363$ & $-5.939603654$ & $-14.79388941$\\
$45/46$ & $91/46$ & $1/46$ &  $-1.043258651$ & $-1.735609347$ & $-9.395596112$ & $-17.80872354$\\
$91/92$ & $91/92$ & $1/92$ &  $-0.6459615671$ & $-3.780635918$ & $-10.4825843$ & $-25.35709155$\\
\multicolumn{7}{l}{\rule{0pt}{3.2ex}$B_{324}$: $|D|=324$, $\mu=1/27$, $\dim\pi_{\sigma_0}=36$, $t_\mu=1.713556331100267$, $E=162d1$}\\
$1/27$ & $1/27$ & $26/27$ & $q^{-1/27}+{}$ $-5/4$ & $-3$ & $-3.813836332$ & $-7.357753989$\\
$5/54$ & $5/54$ & $49/54$ &  $-1.024982278$ & $-1.488902318$ & $-3.5$ & $-5.084149898$\\
$23/216$ & $23/216$ & $193/216$ &  $-1.257491092$ & $-2.616776526$ & $-5.061793053$ & $-7.878613845$\\
$4/27$ & $4/27$ & $23/27$ &  $0.2080505658$ & $-1.373990076$ & $-1.769635806$ & $-3.625$\\
$1/6$ & $1/6$ & $5/6$ &  $-1.5$ & $-2$ & $-4.565572653$ & $-7.411121827$\\
$1/6$ & $7/6$ & $5/6$ &  $-0.75$ & $-0.7754338648$ & $-2$ & $-2.705560914$\\
$11/54$ & $65/54$ & $43/54$ &  $-0.4669173319$ & $-0.7347967888$ & $-1.72183182$ & $-2$\\
$47/216$ & $263/216$ & $169/216$ &  $-0.46654037$ & $-1.001514941$ & $-1.462612132$ & $-2.668186814$\\
$7/27$ & $34/27$ & $20/27$ &  $-0.3723273362$ & $-0.7076537334$ & $-1.75$ & $-1.722817336$\\
$7/24$ & $7/24$ & $17/24$ &  $-1$ & $-2.36019393$ & $-4.330641112$ & $-7$\\
$7/24$ & $55/24$ & $17/24$ &  $-0.06557265273$ & $-0.36019393$ & $-1$ & $-1.474962857$\\
$17/54$ & $17/54$ & $37/54$ &  $-0.5$ & $-2.057482223$ & $-2.181403654$ & $-4.755542877$\\
$71/216$ & $71/216$ & $145/216$ &  $-1.318596346$ & $-1.87762505$ & $-3.736985448$ & $-6.91269579$\\
$1/3$ & $1/3$ & $2/3$ &  $-0.4058291983$ & $-1.418448443$ & $-1.797085401$ & $-4.127261432$\\
$10/27$ & $10/27$ & $17/27$ &  $0.140680299$ & $0.3091078515$ & $0.3145252606$ & $-0.01249113916$\\
$23/54$ & $23/54$ & $31/54$ &  $-0.4863856789$ & $-1.25$ & $-1.75702874$ & $-3.878745546$\\
$95/216$ & $95/216$ & $121/216$ &  $-0.3577250172$ & $-1.613709414$ & $-2.467560367$ & $-3.771342271$\\
$13/27$ & $13/27$ & $14/27$ &  $-0.7598471496$ & $-2.247304529$ & $-4.429659851$ & $-6.564605499$\\
$29/54$ & $29/54$ & $25/54$ &  $-0.7479801517$ & $-1.712667871$ & $-3$ & $-5.891728151$\\
$119/216$ & $119/216$ & $97/216$ &  $-0.7956096337$ & $-1.905118711$ & $-3.997172834$ & $-6.349543427$\\
$16/27$ & $16/27$ & $11/27$ &  $-0.381784297$ & $-1.980825195$ & $-2.274660698$ & $-4.395974717$\\
$35/54$ & $89/54$ & $19/54$ &  $0.0383496105$ & $-1.310110315$ & $-1.222511017$ & $-3$\\
$143/216$ & $143/216$ & $73/216$ &  $-0.8101103146$ & $-1.721182308$ & $-3$ & $-4.654608907$\\
$2/3$ & $2/3$ & $1/3$ &  $-0.4772474181$ & $-2$ & $-2.154964707$ & $-5.25$\\
$2/3$ & $5/3$ & $1/3$ &  $-0.5227525819$ & $0.2386237091$ & $-1$ & $-0.75$\\
$19/27$ & $19/27$ & $8/27$ &  $-0.3043718987$ & $-0.8890133285$ & $-1.256807161$ & $-3.056986624$\\
$41/54$ & $95/54$ & $13/54$ &  $-0.25$ & $0.03503244722$ & $-0.4288625086$ & $-0.9039388599$\\
$167/216$ & $167/216$ & $49/216$ &  $-0.5249822783$ & $-0.4468263599$ & $-1.21896504$ & $-2.448461231$\\
$22/27$ & $49/27$ & $5/27$ &  $-0.2482182241$ & $-0.2318760338$ & $-0.7643856814$ & $-0.766541334$\\
$5/6$ & $5/6$ & $1/6$ &  $-0.1883416035$ & $-0.8267685669$ & $-1.5$ & $-3.346116537$\\
$47/54$ & $209/54$ & $7/54$ &  $-0.2401528504$ & $0.25$ & $-0.224846505$ & $-0.03345962997$\\
$191/216$ & $191/216$ & $25/216$ &  $-0.2520198483$ & $-0.145744228$ & $-1.067740444$ & $-1.853595213$\\
$25/27$ & $25/27$ & $2/27$ &  $-0.3912562026$ & $-0.5828710812$ & $-1.689961787$ & $-2.059418399$\\
$23/24$ & $23/24$ & $1/24$ &  $-0.3116583965$ & $-0.7506732828$ & $-1.491672574$ & $-2.603370105$\\
$53/54$ & $53/54$ & $1/54$ &  $-0.2825124052$ & $-0.7444388248$ & $-1.25$ & $-2.740223711$\\
$215/216$ & $215/216$ & $1/216$ &  $0.06502481036$ & $-0.04622136916$ & $-0.2382351995$ & $-0.2290662823$\\
\multicolumn{7}{l}{\rule{0pt}{3.2ex}$B_{648}$: $|D|=648$, $\mu=23/432$, $\dim\pi_{\sigma_0}=36$, $t_\mu=1.064925231128516$, $E=324a1$}\\
$23/432$ & $23/432$ & $409/432$ & $q^{-23/432}+{}$ $-1.657058442$ & $-5.266295056$ & $-9.991012754$ & $-21.51155967$\\
$11/108$ & $11/108$ & $97/108$ &  $-2.950696391$ & $-7.972087441$ & $-15.68402946$ & $-29.33817807$\\
$47/432$ & $47/432$ & $385/432$ &  $-2.371198123$ & $-5.042404293$ & $-10.55273866$ & $-20.34823033$\\
$7/48$ & $7/48$ & $41/48$ &  $-1.66027871$ & $-4.765747622$ & $-8.701638287$ & $-18.38289298$\\
$71/432$ & $503/432$ & $361/432$ &  $-1.015705374$ & $-2.281815891$ & $-5.793641239$ & $-9.535851512$\\
$23/108$ & $23/108$ & $85/108$ &  $-1.517306224$ & $-5.180174609$ & $-10.97011125$ & $-19.82852922$\\
$95/432$ & $959/432$ & $337/432$ &  $0.064070929$ & $-1.312939443$ & $0.05129120175$ & $-3$\\
$1/4$ & $1/4$ & $3/4$ &  $-2.666666667$ & $-6.427293972$ & $-14$ & $-26.60253097$\\
$1/4$ & $5/4$ & $3/4$ &  $-0.8865601402$ & $-2$ & $-5.096614643$ & $-9.205061933$\\
$119/432$ & $551/432$ & $313/432$ &  $-0.8159705394$ & $-1.414208108$ & $-4.057761332$ & $-6.235740183$\\
$5/16$ & $5/16$ & $11/16$ &  $-2.193229286$ & $-7.560213053$ & $-15.91196247$ & $-30.68656657$\\
$5/16$ & $21/16$ & $11/16$ &  $-1$ & $-2.595697591$ & $-6.911962467$ & $-11.34328328$\\
$35/108$ & $35/108$ & $73/108$ &  $-2.083560176$ & $-4.554245893$ & $-11.02248453$ & $-21.40972236$\\
$143/432$ & $143/432$ & $289/432$ &  $-0.9550309347$ & $-1.904336328$ & $-4.240361264$ & $-8.217277639$\\
$167/432$ & $167/432$ & $265/432$ &  $-0.7995827908$ & $-3.81005927$ & $-7.162504383$ & $-14.56390138$\\
$47/108$ & $47/108$ & $61/108$ &  $-1.604249407$ & $-2.667741474$ & $-7.29713168$ & $-13.68559906$\\
$191/432$ & $191/432$ & $241/432$ &  $-0.3001880383$ & $-0.7388633884$ & $-2.290155881$ & $-4.218515485$\\
$215/432$ & $215/432$ & $217/432$ &  $-0.02241412428$ & $-2.87211206$ & $-3.67077789$ & $-10.06638915$\\
$59/108$ & $59/108$ & $49/108$ &  $-1.617051549$ & $-4.808675755$ & $-11.3997914$ & $-21.38574248$\\
$239/432$ & $239/432$ & $193/432$ &  $-1.639650782$ & $-2.969421036$ & $-7.725983026$ & $-14.36535818$\\
$7/12$ & $7/12$ & $5/12$ &  $-1.807592899$ & $-4.33157293$ & $-9.927204931$ & $-18.83013935$\\
$263/432$ & $263/432$ & $169/432$ &  $-0.6645170525$ & $-4.131061196$ & $-6.387928894$ & $-15.92173829$\\
$31/48$ & $31/48$ & $17/48$ &  $-0.33685414$ & $-3.20229543$ & $-5.917140748$ & $-14.53228455$\\
$71/108$ & $287/108$ & $37/108$ &  $-0.54957576$ & $-0.1753519697$ & $-1.622021761$ & $-2.984294626$\\
$287/432$ & $287/432$ & $145/432$ &  $-1.350703939$ & $-3.536167546$ & $-6.672135701$ & $-15.73410892$\\
$311/432$ & $311/432$ & $121/432$ &  $0.2979493968$ & $-0.8791988815$ & $0$ & $-2.529189814$\\
$83/108$ & $83/108$ & $25/108$ &  $-0.5277443365$ & $-1.655840367$ & $-4.849905981$ & $-8.567183717$\\
$335/432$ & $335/432$ & $97/432$ &  $-0.901392782$ & $-2.841794709$ & $-4.073465915$ & $-10.80737788$\\
$13/16$ & $13/16$ & $3/16$ &  $-1.22687972$ & $-2.870172447$ & $-7$ & $-12.88771529$\\
$13/16$ & $29/16$ & $3/16$ &  $-0.5602130529$ & $-1$ & $-2.642920517$ & $-4.112284713$\\
$359/432$ & $791/432$ & $73/432$ &  $-0.1671203513$ & $-1.5662033$ & $-2.702839274$ & $-6.04542268$\\
$95/108$ & $203/108$ & $13/108$ &  $-0.1758142071$ & $-1.148974698$ & $-2.173092285$ & $-3.467964536$\\
$383/432$ & $383/432$ & $49/432$ &  $-0.7658969025$ & $-2.535644255$ & $-5.138935801$ & $-9.599101838$\\
$407/432$ & $839/432$ & $25/432$ &  $-0.05548867306$ & $-0.5948184237$ & $-1.258241265$ & $-3.144997327$\\
$107/108$ & $107/108$ & $1/108$ &  $-0.3298402103$ & $-1.698459841$ & $-3.994396469$ & $-7.564837159$\\
$431/432$ & $863/432$ & $1/432$ &  $-0.1596804207$ & $-0.680259277$ & $-1.748568285$ & $-2.895337038$\\
\bottomrule
\end{longtable}
\end{center}

\section{Computational record}\label{app:comp}
Every statement marked \gC\ rests on a computation listed below, and where a statement was computed in more than one way, all the scripts are
listed; the README of the folder \texttt{v100} of the reproduction archive names the log that each script writes. The proof of Theorem~A and the computations of the forms were made by two implementations that share no
code; their scripts are in the folders \texttt{implementation\_A} and \texttt{implementation\_B} of the archive, and its file \texttt{STATEMENTS.tsv}
lists, row by row, the scripts and logs of both. Below, ``rows'' are rows of \texttt{STATEMENTS.tsv}, where the scripts of implementation B are
named, and A, R, F are the folders \texttt{w1015\_classification\_proof}, \texttt{w1015\_rootbases\_wt52}, \texttt{further\_scripts} of
\texttt{implementation\_A}. Python is run with \texttt{-I}; \texttt{weilrep} [Wi] under Sage~10; PARI/GP~2.15.
\begin{center}\scriptsize
\begin{longtable}{>{\raggedright\arraybackslash}p{56mm}>{\raggedright\arraybackslash}p{96mm}}
\toprule
statement & scripts and logs\\\midrule
\endhead
Tables~\ref{tab:twelve} and~\ref{tab:canon}: every entry against the records of the scripts below, as tabulated in draft v0.9.8 (folder \texttt{paper\_v098}) & \texttt{check\_tables12.py}\\
\S\ref{subsec:weil}: $D(B_{40})$ and $D(L_{72})$ are the modules $2_5^{+3}5^{+1}$ and $2_5^{+3}3^{+2}$ of [BEF16]; \S\ref{subsec:lat}:
$L_{72}\cong\mathcal O^0(-1)$, an explicit change of basis & \texttt{bef16\_symbols.py}; \texttt{errthum\_l72.py}\\
\S\ref{subsec:mult}, the twelve: $\pi_{\sigma_0}$, $\mu$, $c_{F_0}(\mu)\ne0$, the span of the coefficients; \S\ref{subsec:hecke}: the relations \eqref{eq:T32} & F: \texttt{t\_twisted.py}, \texttt{compare\_h12.py}, \texttt{run23.py}, \texttt{hecke\_b648.py}; rows 69, 71, 86\\
\S\ref{subsec:slices}: row $72$, the two slices, their eigenvalues, $t$ at $n=\frac7{72}$ & \texttt{w32\_basis.py}, \texttt{w32\_slices.py}, \texttt{w32\_slice\_hecke.py}, \texttt{t\_slices.py}, \texttt{c72\_curves.gp}\\
\S\ref{subsec:slices}: counts of slices among the $278$ & \texttt{mf\_list.py}\\
\S\ref{subsec:table}, Remark~\ref{rem:I-II}: the values of $t_\mu$ (quadrature; Hejhal's method), $B_\mu=3B_S$, $t_\mu/t$ & F: \texttt{compare\_h12.py} ($t_\mu$ to $16$ digits, two implementations), \texttt{compare.py} (Hejhal's method; $t_\mu/t$),
\texttt{beta\_exact.py} ($B_\mu=3B_S$)\\
\S\ref{subsec:b40}: the matrix $\mathbf S$ and the coefficients of $\Hh$ & F: \texttt{run23.py}, \texttt{b40\_check.py}; rows 78--79\\
Prop.~\ref{prop:etatheta}: exponent sums & \texttt{tabH\_sums.py} (on Table~\ref{tab:H})\\
\S\ref{subsec:coeffs}, Table~\ref{tab:H}: the forms $\Hh$ & F: \texttt{harm\_a.py}, \texttt{compare\_h12.py}, \texttt{verify\_table.py}; rows 87, 92\\
Remark~\ref{rem:theta4}: $\pi_{\sigma_0}\cong V^-(\Z/2s_L)$; the coefficients $45,231,770,2277$ & F: \texttt{pi\_iso.py}, \texttt{theta\_norm.py}; rows 63--68, 70\\
\S\ref{sec:arith}: the eigenvalues and the curves; the uniqueness of the isogeny classes & F: \texttt{run23.py}, \texttt{compare.py}; rows 3, 57, 71--73\\
Prop.~\ref{prop:cm}: certificate; the lift to $\operatorname{diag}(-8,6,6)$ & \texttt{cm\_certificate.py}, \texttt{lift\_test.py}; \texttt{cm\_check.py} ($j$-invariants of the seven curves)\\
\S\ref{subsec:fam}: $\fd=2p$ for $p\le97$; the lattices $G'(2)$ for the $211$ primes $p\equiv1\pmod4$ below $3000$; Lemma~\ref{lem:two} on $18952$ Gram matrices & \texttt{qdisc8.gp}, \texttt{qdisc\_dp.gp} (every odd $p\le199$), \texttt{classical\_checks.gp}\\
Lemma~\ref{lem:fix}: $F_N$ by Lefschetz, $5\le p\le61$ & \texttt{fixpts.gp}\\
Lemma~\ref{lem:three}: $190$ primes $p\le1163$; Ogg's list & \texttt{hier\_counts.gp}, \texttt{run\_w2.gp}, \texttt{ogg\_alpha.py}; \texttt{profile\_b.log} ($p\le47$), \texttt{ogg\_check.log} ($p\le109$); comparison \texttt{cmp\_profile\_b.py}\\
after Lemma~\ref{lem:three}: the traces of [Ma25] for $e=0,1,2$ and the counts, $51$ primes $p\le250$ & \texttt{martin\_table2.gp}\\
Cor.~\ref{cor:geom}: class numbers below $10^5$ & \texttt{tiers\_classno.gp}\\
Prop.~\ref{prop:allw}: traces; counts $p\le149$, $w\le8$ & \texttt{lefschetz\_w.gp}, \texttt{run\_all\_w.gp}, \texttt{check\_allw.py}; \texttt{profile\_b.log} ($52$ rows, $p\le47$)\\
Remark~\ref{rem:dims}: $\dim S_k(\rho_{\Dp_p})$ and the newform counts, $3\le p\le103$, $k\le\frac92$ & \texttt{dp\_dims.py} (left side); \texttt{hier\_counts.gp}, \texttt{run\_all\_w.gp} (right side); comparison \texttt{cmp\_dims.py};
\texttt{profile\_b.log} ($p\le47$), comparison \texttt{cmp\_dims\_b.py}\\
\S\ref{subsec:dim}: \eqref{eq:rr} against $\RR-d_0$ ($36{,}683$ modules) and \texttt{weilrep} ($97$ census lattices) & A: \texttt{vsearch/v\_validate\_wp0.py}; row 11\\
\S\ref{subsec:saw}: the list for $T=10^4$ & A: \texttt{vsearch/v\_list.py}; rows 14--16\\
\S\ref{subsec:large}: V1--V3 and the route of implementation B & A: \texttt{vsearch/v\_leaf.py}, \texttt{v\_exact\_run.py}, \texttt{v\_cand\_ordO\_all.py}; rows 17--20\\
\S\ref{subsec:small}: the modules with $|D|\le10^4$, the count of genera, the exclusion & A: \texttt{vsearch/v\_census.py}, \texttt{v\_validate\_wp0.py}, \texttt{v\_partA.py}, \texttt{v\_partA\_mstar.py}; rows 21--22\\
\S\ref{subsec:chi}: the character on the $278$ candidates, the checks, the twenty lattice modules above $648$, $d_0(\sigma_0)=m_{1/2}(\sigma_0)=0$ &
A: \texttt{chis/run278.py}, \texttt{cmp278.py}, \texttt{test1.py}, \texttt{test\_mfac.py}; F: \texttt{c3\_big.py}, \texttt{t\_twisted.py}; rows 23--26, 90\\
\S\ref{subsec:nik}: Nikulin's conditions; the models of the twelve & A: \texttt{chis/nikulin.py}, \texttt{gramcheck.py}, \texttt{gramcheck2.py}; R: \texttt{gram\_alt.py}; rows 28--30\\
\S\ref{subsec:classno}: class number one (Sage); the theta lifts & F: \texttt{theta\_dim.py}; rows 1, 31, 59\\
\S\ref{subsec:roots}: the core; the root bases of the twelve and their certificates; the census rows; the scaled census & R: \texttt{core44.py}, \texttt{rb.py},
\texttt{run12.py}, \texttt{verify.py}, \texttt{membership.py}, \texttt{chains.py}, \texttt{scaled.py}, \texttt{check\_p2.py}; rows 6, 10, 32--37; the list \texttt{ROOT\_BASES.md}: \texttt{w1116/root\_bases.py}\\
Table~\ref{tab:canon}: the columns of weight $\frac52$ & R: \texttt{thmC\_table.py} (\texttt{weil.py}, \texttt{invdim.py}), \texttt{calib52.py}, \texttt{calib\_so.py}; rows 38--41\\
Problem~\ref{prob:rat}: the rational coefficients of $\Hh$, the twists, the expected number of accidental agreements & implementation B:
\texttt{w1019\_review/hm/rat\_test.py}, \texttt{rat\_risk.py} (logs \texttt{rat\_test\_*.log}, \texttt{rat\_summary.log}, \texttt{rat\_risk.log}); row 88.
The $44$ pairs of $B_{324}$: root numbers, reduction at $3$, $L$ and $L'$ at $s=1$ (two implementations): F: \texttt{lprime\_a.py};
implementation B: \texttt{w1020\_review/lprime\_check.py}; row 89\\
\bottomrule
\end{longtable}
\end{center}
The weight-$\frac52$ computations made with these (the level-$23$ block of $B_{184}$, the $3$-adic descent for $B_{324}$ and $B_{648}$, the $\eta$--$\theta$ test at
weight $\frac52$, the pure census rows of [II]) are archived with these scripts and are not used here.

\section{Quoted theorems}\label{app:cite}
\begin{center}\footnotesize
\begin{longtable}{>{\raggedright\arraybackslash}p{4mm}>{\raggedright\arraybackslash}p{106mm}>{\raggedright\arraybackslash}p{38mm}}
\toprule
& statement & used in\\\midrule
\endhead
1 & Milgram's formula [MH73, App.~4] & \S\S\ref{subsec:fqm}, \ref{subsec:fam}\\
2 & Genus determined by signature and discriminant form [Ni79, Cor.~1.9.4] & \S\S\ref{subsec:lat}, \ref{subsec:classno}, \ref{subsec:roots}\\
3 & Even Clifford algebra and ramification [Vo21, Prop.~5.1.2, 5.3.23]; compactness of Shimura curves [Vo21, Prop.~38.1.3, Main Thm.~38.4.3] & \S\S\ref{subsec:lat}, \ref{subsec:litrel}\\
4 & Newforms and strong multiplicity one [AL70]; $a_\ell=-w_\ell\ell^{w/2-1}$ for $\ell\,\|\,N$ [AL70, Thm.~3] & \S\S\ref{subsec:newforms}, \ref{subsec:fam}\\
5 & Eichler--Shimura; modularity [BCDT01] & \S\ref{subsec:newforms}\\
6 & Hecke operators in odd signature [BS10, Thm.~4.10, Rem.~4.11] & \S\ref{subsec:hecke}\\
7 & Schur's lemma [Se77, \S2.2]; $\sum\chi(1)^2=|H|$ [Se77, \S2.4]; abelian subgroups of index $j$ [Se77, \S3.1] & \S\S\ref{subsec:hecke}, \ref{subsec:onenumber}, \ref{sec:proofA}\\
8 & Harmonic weak Maass forms [BF04, \S3: Prop.~3.5, Cor.~3.8, Prop.~3.11] & Thm.~\ref{thm:harm}\\
9 & Residue $3/\pi$ of $E(\tau,s)$ [Bu97, \S1.6]; Wiener--Ikehara & Remark~\ref{rem:mean}\\
10 & Theta function of $\langle2m\rangle$ [Bo98, Thm.~4.1]; valence formula [Se73, Ch.~VII, \S3, Thm.~3] & Prop.~\ref{prop:etatheta}\\
11 & Integral bases [McG03, Thm.~5.6] & Prop.~\ref{prop:cm}\\
12 & Deuring ($a_p=0$ at inert primes of good reduction); integrality of CM $j$-invariants, hence potentially good and additive bad reduction; the nine imaginary quadratic fields of class number one & Prop.~\ref{prop:cm}\\
13 & Bruinier--Ono [BO10, Thm.~1.1(2), 1.2]; [Br13, Thm.~3.2] & \S\ref{subsec:litrel}, Problem~\ref{prob:rat} (compared; not used)\\
14 & Quotient curves (Riemann--Hurwitz) & Lemma~\ref{lem:three}\\
15 & Holomorphic Lefschetz formula [AB68] & Prop.~\ref{prop:allw}\\
16 & Genus and dimension formulas for $\Gamma_0(N)$ [DS05, Thm.~3.1.1, Thm.~3.5.1] & Lemma~\ref{lem:three}, Prop.~\ref{prop:allw}\\
17 & Class numbers of orders [Co13, Thm.~7.24] & Lemma~\ref{lem:fix}\\
18 & Fixed points of Atkin--Lehner involutions [Og74] & after Lemma~\ref{lem:fix} (compared; not used)\\
19 & $X_0(p)/w_p$ has genus $0$ exactly for the supersingular primes [Og75] & after Lemma~\ref{lem:three}\\
20 & Hilbert reciprocity; quadratic reciprocity; Dirichlet's theorem on primes in arithmetic progressions; the class number formula
$h(-d)=\frac{\sqrt d}\pi L(1,\chi_{-d})$ & Cor.~\ref{cor:fd}, \S\ref{subsec:fam}, Cor.~\ref{cor:geom}\\
21 & The dual normalization of the Weil representation; anisotropic modules with $S_{3/2,A}=0$ [BEF16, (8)--(9), Table~4] & \S\ref{subsec:weil} (compared; not used)\\
22 & A maximal order of $(\frac{5,6}{\Q})$ and its lattice of trace zero [Err11, \S3.1, Prop.~3.1, (7.1)] & \S\ref{subsec:lat}\\
23 & Traces of Atkin--Lehner operators on new spaces [Ma25, \S4.1, Table~2] & after Lemma~\ref{lem:three} (compared; not used)\\
24 & Petersson norm of the unary theta series of weight $\frac32$ [BSc17, Example~5.2] & Remark~\ref{rem:theta4}\\
25 & The completed umbral forms of pure A-type lambency as harmonic weak Maass forms for the Weil representation [ORTL15, Thm.~4.1] & Remark~\ref{rem:theta4}\\
26 & The even $2n$ for which $\Gamma_0(2n)+n$ has genus $0$ [CDH14b, \S2.3, (2.21)] & Cor.~\ref{cor:geom}(b) (compared; not used)\\
27 & Normal forms of finite quadratic modules [Wa63, \S5 Thm.~4, \S6]; [Ni79, Prop.~1.8.1] & \S\ref{sec:proofA}\\
28 & Riemann--Roch for vector-valued modular forms [ES95, \S4.2], [Sk85, p.~100] & Thm.~\ref{thm:dim}, Prop.~\ref{prop:chiS}, \S\ref{subsec:wt52}\\
29 & Serre duality for vector-valued forms [Bo99, Thm.~3.1] & Thm.~\ref{thm:dim}, Prop.~\ref{prop:chiS}\\
30 & Point counts of quadrics over $\F_p$ [LN97, Thms.~6.26--6.27]; characters of $\GL_2(\F_q)$ [FH91, Lecture~5] & Lemma~\ref{lem:local}\\
31 & Degrees of characters [Is76, Cor.~2.30, Thm.~3.12] & Lemma~\ref{lem:deg}\\
32 & Fourier expansion of the sawtooth function [Be76] & Lemma~\ref{lem:FL}\\
33 & $M_{1/2}(\Gamma(4m))$ as an $\Mp$-module [Sk08, \S5], [Sk85, Thm.~5.2, Thm.~1.8]; Weil representations [Sk08, \S3], [Sch09, Prop.~4.3] & \S\ref{subsec:chi}\\
34 & Existence of even lattices [Ni79, Thm.~1.10.1] & \S\ref{subsec:nik}\\
35 & Class numbers of indefinite forms [CS99, Ch.~15, Thm.~14, \S7]; spinor norms and spinor genera [OM00, \S\S101--102] & \S\ref{subsec:classno}\\
36 & Lifting through isotropic subgroups [BEF16, Prop.~4.18] & \S\ref{subsec:classno}, Remark~\ref{rem:theta4}\\
37 & Hyperbolic and compact generalized Cartan matrices [Ka90, Ch.~4] & \S\ref{subsec:roots}\\
38 & Eisenstein series of the Weil representation [Br02, Thm.~1.6]; the dimension formula [BEF16, (3.4)] & \S\ref{subsec:wt52}\\
\bottomrule
\end{longtable}
\end{center}

\section*{Acknowledgements}
\addcontentsline{toc}{section}{Acknowledgements}
The scripts of Appendix~\ref{app:comp} and drafts of the text were written with the help of AI language models (Claude, Anthropic). The author takes
full responsibility for the content.

\section*{References}
\addcontentsline{toc}{section}{References}
{\small
\begin{list}{}{\setlength{\labelwidth}{15mm}\setlength{\leftmargin}{17mm}\setlength{\labelsep}{2mm}\setlength{\itemsep}{1.5pt}}
\item[{[I]}] E.~Cho, \emph{Forced shadows of an obstructed hyperbolic Kac--Moody denominator}, arXiv:2608.19706.
\item[{[II]}] E.~Cho, \emph{Petersson-rigid lattices in a census of 100 rank-three root bases}, arXiv:2609.02843.
\item[{[AB68]}] M.~F.~Atiyah, R.~Bott, \emph{A Lefschetz fixed point formula for elliptic complexes.\ II.\ Applications}, Ann.\ of Math.\ (2) \textbf{88} (1968), 451--491.
\item[{[AGOR15]}] C.~Alfes, M.~Griffin, K.~Ono, L.~Rolen, \emph{Weierstrass mock modular forms and elliptic curves}, Res.\ Number Theory \textbf{1} (2015),
doi:10.1007/s40993-015-0026-2.
\item[{[AL70]}] A.~O.~L.~Atkin, J.~Lehner, \emph{Hecke operators on $\Gamma_0(m)$}, Math.\ Ann.\ \textbf{185} (1970), 134--160.
\item[{[BCDT01]}] C.~Breuil, B.~Conrad, F.~Diamond, R.~Taylor, \emph{On the modularity of elliptic curves over $\Q$: wild $3$-adic exercises}, J.\ Amer.\ Math.\
Soc.\ \textbf{14} (2001), 843--939.
\item[{[BEF16]}] J.~H.~Bruinier, S.~Ehlen, E.~Freitag, \emph{Lattices with many Borcherds products}, Math.\ Comp.\ \textbf{85} (2016), 1953--1981.
\item[{[Be76]}] B.~C.~Berndt, \emph{Classical theorems on quadratic residues}, Enseign.\ Math.\ (2) \textbf{22} (1976), 261--304.
\item[{[BF04]}] J.~H.~Bruinier, J.~Funke, \emph{On two geometric theta lifts}, Duke Math.\ J.\ \textbf{125} (2004), 45--90.
\item[{[BO10]}] J.~H.~Bruinier, K.~Ono, \emph{Heegner divisors, $L$-functions and harmonic weak Maass forms}, Ann.\ of Math.\ (2) \textbf{172} (2010), 2135--2181.
\item[{[Bo98]}] R.~E.~Borcherds, \emph{Automorphic forms with singularities on Grassmannians}, Invent.\ Math.\ \textbf{132} (1998), 491--562.
\item[{[Bo99]}] R.~E.~Borcherds, \emph{The Gross--Kohnen--Zagier theorem in higher dimensions}, Duke Math.\ J.\ \textbf{97} (1999), 219--233; correction, ibid.\ \textbf{105} (2000), 183--184.
\item[{[Br02]}] J.~H.~Bruinier, \emph{Borcherds Products on O(2,l) and Chern Classes of Heegner Divisors}, Lecture Notes in Math.\ 1780, Springer, 2002.
\item[{[Br13]}] J.~H.~Bruinier, \emph{Harmonic Maass forms and periods}, Math.\ Ann.\ \textbf{357} (2013), 1363--1387.
\item[{[BS10]}] J.~H.~Bruinier, O.~Stein, \emph{The Weil representation and Hecke operators for vector valued modular forms}, Math.\ Z.\ \textbf{264} (2010), 249--270.
\item[{[BSc17]}] J.~H.~Bruinier, M.~Schwagenscheidt, \emph{Algebraic formulas for the coefficients of mock theta functions and Weyl vectors of
Borcherds products}, J.\ Algebra \textbf{478} (2017), 38--57; arXiv:1607.02701.
\item[{[BSt12]}] J.~H.~Bruinier, F.~Str\"omberg, \emph{Computation of harmonic weak Maass forms}, Exp.\ Math.\ \textbf{21} (2012), 117--131.
\item[{[Bu97]}] D.~Bump, \emph{Automorphic Forms and Representations}, Cambridge Stud.\ Adv.\ Math.\ 55, Cambridge Univ.\ Press, 1997.
\item[{[CDH14a]}] M.~C.~N.~Cheng, J.~F.~R.~Duncan, J.~A.~Harvey, \emph{Umbral moonshine}, Commun.\ Number Theory Phys.\ \textbf{8} (2014), 101--242; arXiv:1204.2779.
\item[{[CDH14b]}] M.~C.~N.~Cheng, J.~F.~R.~Duncan, J.~A.~Harvey, \emph{Umbral moonshine and the Niemeier lattices}, Res.\ Math.\ Sci.\ \textbf{1}
(2014), Art.~3; arXiv:1307.5793.
\item[{[Co13]}] D.~A.~Cox, \emph{Primes of the Form $x^2+ny^2$: Fermat, Class Field Theory, and Complex Multiplication}, 2nd ed.,
Wiley, Hoboken, NJ, 2013.
\item[{[CS99]}] J.~H.~Conway, N.~J.~A.~Sloane, \emph{Sphere Packings, Lattices and Groups}, 3rd ed., Grundlehren 290, Springer, 1999.
\item[{[Cr97]}] J.~E.~Cremona, \emph{Algorithms for Modular Elliptic Curves}, 2nd ed., Cambridge Univ.\ Press, 1997; extended tables: \texttt{github.com/JohnCremona/ecdata}.
\item[{[DGO15]}] J.~F.~R.~Duncan, M.~J.~Griffin, K.~Ono, \emph{Proof of the umbral moonshine conjecture}, Res.\ Math.\ Sci.\ \textbf{2} (2015), Art.~26.
\item[{[DS05]}] F.~Diamond, J.~Shurman, \emph{A First Course in Modular Forms}, Grad.\ Texts in Math.\ 228, Springer, 2005.
\item[{[Err11]}] E.~Errthum, \emph{Singular moduli of Shimura curves}, Canad.\ J.\ Math.\ \textbf{63} (2011), 826--861.
\item[{[ES95]}] W.~Eholzer, N.-P.~Skoruppa, \emph{Modular invariance and uniqueness of conformal characters}, Comm.\ Math.\ Phys.\ \textbf{174} (1995), 117--136.
\item[{[FF83]}] A.~J.~Feingold, I.~B.~Frenkel, \emph{A hyperbolic Kac--Moody algebra and the theory of Siegel modular forms of genus 2}, Math.\ Ann.\ \textbf{263} (1983), 87--144.
\item[{[FH91]}] W.~Fulton, J.~Harris, \emph{Representation Theory: A First Course}, GTM 129, Springer, 1991.
\item[{[GN96]}] V.~A.~Gritsenko, V.~V.~Nikulin, \emph{Igusa modular forms and ``the simplest'' Lorentzian Kac--Moody algebras}, Sb.\ Math.\ \textbf{187} (1996), no.~11, 1601--1641.
\item[{[Is76]}] I.~M.~Isaacs, \emph{Character Theory of Finite Groups}, Academic Press, 1976.
\item[{[Ka90]}] V.~G.~Kac, \emph{Infinite Dimensional Lie Algebras}, 3rd ed., Cambridge Univ.\ Press, 1990.
\item[{[KL13]}] H.~H.~Kim, K.-H.~Lee, \emph{Root multiplicities of hyperbolic Kac--Moody algebras and Fourier coefficients of modular forms}, Ramanujan J.\
\textbf{32} (2013), 329--352.
\item[{[KP84]}] V.~G.~Kac, D.~H.~Peterson, \emph{Infinite-dimensional Lie algebras, theta functions and modular forms}, Adv.\ Math.\ \textbf{53} (1984), 125--264.
\item[{[LMFDB]}] The LMFDB Collaboration, \emph{The $L$-functions and modular forms database},\\ \url{https://www.lmfdb.org}.
\item[{[LN97]}] R.~Lidl, H.~Niederreiter, \emph{Finite Fields}, 2nd ed., Encyclopedia Math.\ Appl.\ 20, Cambridge Univ.\ Press, 1997.
\item[{[Ma18]}] K.~Martin, \emph{Refined dimensions of cusp forms, and equidistribution and bias of signs}, J.\ Number Theory \textbf{188} (2018), 1--17.
\item[{[Ma25]}] K.~Martin, \emph{Distribution of local signs of modular forms and murmurations of Fourier coefficients}, Mathematika
\textbf{71} (2025), no.~3, Paper No.~e70028, 37~pp.
\item[{[Mac72]}] I.~G.~Macdonald, \emph{Affine root systems and Dedekind's $\eta$-function}, Invent.\ Math.\ \textbf{15} (1972), 91--143.
\item[{[McG03]}] W.~J.~McGraw, \emph{The rationality of vector valued modular forms associated with the Weil representation}, Math.\ Ann.\ \textbf{326} (2003),
105--122.
\item[{[MH73]}] J.~Milnor, D.~Husemoller, \emph{Symmetric Bilinear Forms}, Ergebnisse 73, Springer, 1973.
\item[{[MO97]}] Y.~Martin, K.~Ono, \emph{Eta-quotients and elliptic curves}, Proc.\ Amer.\ Math.\ Soc.\ \textbf{125} (1997), 3169--3176.
\item[{[Ni79]}] V.~V.~Nikulin, \emph{Integral symmetric bilinear forms and some of their applications}, Izv.\ Akad.\ Nauk SSSR Ser.\ Mat.\ \textbf{43} (1979), 111--177;
English transl.\ Math.\ USSR-Izv.\ \textbf{14} (1980), 103--167.
\item[{[Og74]}] A.~P.~Ogg, \emph{Hyperelliptic modular curves}, Bull.\ Soc.\ Math.\ France \textbf{102} (1974), 449--462.
\item[{[Og75]}] A.~P.~Ogg, \emph{Automorphismes de courbes modulaires}, S\'eminaire Delange--Pisot--Poitou (Th\'eorie des nombres) \textbf{16} (1974/75), no.~1,
exp.~no.~7.
\item[{[OM00]}] O.~T.~O'Meara, \emph{Introduction to Quadratic Forms}, Classics in Mathematics, Springer, 2000.
\item[{[ORTL15]}] K.~Ono, L.~Rolen, S.~Trebat-Leder, \emph{Classical and umbral moonshine: connections and $p$-adic properties}, J.\ Ramanujan
Math.\ Soc.\ \textbf{30} (2015), 135--160; arXiv:1403.3712.
\item[{[Sch09]}] N.~R.~Scheithauer, \emph{The Weil representation of $\mathrm{SL}_2(\Z)$ and some applications}, Int.\ Math.\ Res.\ Not.\ IMRN \textbf{2009}, no.~8, 1488--1545.
\item[{[Se73]}] J.-P.~Serre, \emph{A Course in Arithmetic}, Grad.\ Texts in Math.\ 7, Springer, 1973.
\item[{[Se77]}] J.-P.~Serre, \emph{Linear Representations of Finite Groups}, Grad.\ Texts in Math.\ 42, Springer, 1977.
\item[{[Sk85]}] N.-P.~Skoruppa, \emph{\"Uber den Zusammenhang zwischen Jacobiformen und Modulformen halbganzen Gewichts}, Bonner Math.\ Schriften 159, 1985.
\item[{[Sk08]}] N.-P.~Skoruppa, \emph{Jacobi forms of critical weight and Weil representations}, in: Modular Forms on Schiermonnikoog, Cambridge Univ.\ Press, 2008; arXiv:0707.0718.
\item[{[Vo21]}] J.~Voight, \emph{Quaternion Algebras}, Grad.\ Texts in Math.\ 288, Springer, 2021.
\item[{[Wi]}] B.~Williams, \texttt{weilrep}, a SageMath package for vector-valued modular forms (software).
\item[{[Wa63]}] C.~T.~C.~Wall, \emph{Quadratic forms on finite groups, and related topics}, Topology \textbf{2} (1963), 281--298.
\item[{[Ya73]}] M.~Yamauchi, \emph{On the traces of Hecke operators for a normalizer of $\Gamma_0(N)$}, J.\ Math.\ Kyoto Univ.\ \textbf{13} (1973), no.~2, 403--411.
\end{list}
}

\end{document}